\documentclass[11pt]{article}
\usepackage[T1]{fontenc}
\usepackage[utf8]{inputenc}
\usepackage{lmodern}
\usepackage{amsmath,amssymb,amsthm,mathtools,mathrsfs,microtype,booktabs}
\usepackage[margin=1in]{geometry}
\usepackage{needspace,graphicx,tikz}
\usetikzlibrary{arrows.meta}
\usepackage[hidelinks]{hyperref}
\newtheorem{theorem}{Theorem}[section]
\newtheorem{lemma}[theorem]{Lemma}
\newtheorem{proposition}[theorem]{Proposition}
\newtheorem{corollary}[theorem]{Corollary}
\theoremstyle{definition}\theoremstyle{plain}
\theoremstyle{remark}\newtheorem{remark}[theorem]{Remark}
\newcommand{\D}{\mathbb D}
\newcommand{\Ds}{\mathbb D^*}
\newcommand{\E}{\mathbb E}
\newcommand{\cL}{c_{\mathrm L}}
\newcommand{\cm}{c_{\mathrm m}}

\newcommand{\dd}{\,d^2z}
\DeclareMathOperator{\Rea}{Re}
\DeclareMathOperator{\Ima}{Im}
\title{Relative entropy and K\"ahler geometry of SLE weldings}
\author{Gunhee Cho\thanks{Department of Mathematics, Texas State University,
San Marcos, TX, USA. wvx17@txstate.edu}
\and Vivian Olsiewski Healey\thanks{Department of Mathematics, Texas State
University, San Marcos, TX, USA. healey@txstate.edu}
\and Jihun Yum\thanks{Department of Mathematics Education and Institute of
Pure and Applied Mathematics, Jeonbuk National University, Jeonju, Republic
of Korea. jihun0224@jbnu.ac.kr}}
\date{November 24, 2026}
\begin{document}
\maketitle
\begin{abstract}
We study the rotation-randomized SLE welding law \(P\), \(0<\kappa\le4\),
under composition by deterministic circle maps. For every Weil--Petersson
homeomorphism \(a\) we prove the exact identity
\[
 \operatorname{Ent}(P_a\,|\,P)=H_\kappa(a):=\frac{c_{\mathrm L}}{24}\,I^L(a)+2K(a),
\]
where \(I^L\) is the Loewner energy, \(K\) the log ratio of the exterior
and interior conformal radii, and \(c_{\mathrm L}=13+24/\kappa+3\kappa/2\). The
function \(H_\kappa\) is not new: by Takhtajan--Teo it is a global Kähler
potential of the right-invariant Kähler metric
\(2g_{\rm VK}+\frac{c_{\mathrm L}}6g_{\rm WP}\) on the Weil--Petersson universal
Teichm\"uller curve, which extends the Kirillov--Yuriev metric
\(g_{c_{\mathrm L},1}\) of \(S^1\backslash\operatorname{Diff}_+(S^1)\), and by a recent
preprint of Cai (independent of our work) this metric has the Fourier form of the score covariance
of the welding law. What is novel is the identity of this potential: it is
the relative entropy of the welding laws. Consequently, the relative
entropy \(h\mapsto\operatorname{Ent}(P_{h^{-1}}|P)\) is a global Kähler potential of
the Fisher information metric of the welding family, and this metric is
the Hessian of the relative entropy in every direction. A Kähler metric
determines its potentials only up to pluriharmonic functions, and the
Fisher metric, an infinitesimal quantity, does not determine the
entropy; the identity requires a global computation, which we carry out
by a mean-value argument along radial holomorphic discs in the universal
Teichm\"uller curve. Quasi-invariance of \(P\) was known for \(\kappa<4\)
(Fan--Sung) and, with explicit densities, for analytic maps
(Baverez--Jego); the identity gives it for all Weil--Petersson maps also
at \(\kappa=4\), as obtained concurrently, by different methods, by S.~Fan.
The probabilistic input, a horizontal integration-by-parts formula for
the welding law contained in the preprints of Baverez--Jego and of Cai,
is proved here from published results on the SLE loop measure, together
with, at \(\kappa=4\), the removability of SLE\(_4\) loops
(Kavvadias--Miller--Schoug).
\end{abstract}

\noindent\textbf{MSC 2020:} 60J67, 30F60, 32G15, 53C55, 30C62, 81T40, 94A17.\\
\noindent\textbf{Keywords:} SLE loop, conformal welding, relative entropy,
Loewner energy, Kähler potential, Weil--Petersson Teichm\"uller
space, Kirillov--Yuriev metric, quasi-invariance, Fisher information.

\begingroup\scriptsize\tableofcontents\endgroup

\section{Introduction}\label{sec:introduction}\label{sec:prior}

A conformal welding records a Jordan curve by the correspondence between
its two boundary coordinates: if \(f\) and \(g\) map the unit disc \(\D\) and
its exterior \(\Ds\) conformally onto the two sides of the curve, the
welding is the circle homeomorphism \(g^{-1}\circ f|_{S^1}\). For an SLE\(_\kappa\)
loop, \(0<\kappa\le4\), the welding is a random element of
\(\operatorname{Homeo}_+(S^1)\), and deterministic circle maps act on it by
composition. This paper computes the cost of that action exactly. The
relative entropy between the welding law and its image under a
Weil--Petersson map \(a\) equals a conformal action \(H_\kappa(a)\), built from
the Loewner energy and a capacity term (Theorem~\ref{thm:A}). The function
\(H_\kappa\) has been known for twenty years as a Kähler potential in
Teichm\"uller theory; what we prove is that it is a relative entropy.
As a consequence, the relative entropy is a global Kähler potential of
the Fisher information metric of the family of deformed laws
(Theorem~\ref{thm:Akahler}). The theorems below are stated so that they
can be read without the rest of this introduction.

The complex-geometric background is the following. The Weil--Petersson
homeomorphisms modulo M\"obius maps form the Weil--Petersson universal
Teichm\"uller space \(T_0(1)\), a complex Hilbert manifold on which
Takhtajan and Teo \cite{TTI,TT} showed that the Weil--Petersson metric is
Kähler--Einstein with negative Ricci curvature, and that the universal
Liouville action \(S_1\) is a global Kähler potential for it; on the
universal Teichm\"uller curve over it, the capacity term \(K\) is a global
Kähler potential of the Velling--Kirillov metric. Wang \cite{WangInv}
proved that \(S_1\) is \(\pi\) times the Loewner energy \(I^L\), the
large-deviation rate function of SLE\(_\kappa\) as \(\kappa\to0\), so that a
quantity from probability is the potential of a Kähler--Einstein metric.
Why this should be so, and what the Kähler structure means for SLE at a
fixed \(\kappa>0\), is not yet fully understood; see the survey \cite{Wang}. 
The same function 
has since appeared at fixed \(\kappa\) in several roles. For
Carfagnini and Wang \cite{CarfagniniWang} the Loewner energy is, up to the
central charge, the Onsager--Machlup functional of the loop measure. For
Baverez and Jego \cite{BJweld}, \(H_\kappa\) is the action of a formal density
of the welding law. For Cai \cite{CaiResponse} its increments are
log-densities of finite-mode deformations, and its Hessian at the round
circle is the covariance of the score. For Fan \cite{FanCM} a renormalized
relative action built from it is the log-density for all Weil--Petersson
deformations. In the present paper, \(H_\kappa\) is the relative entropy
itself. So the new result is not a new potential but the identity of a
known one, and with it the Kähler metric becomes the Fisher information
metric of the welding laws with the entropy as its potential; the
Einstein condition then singles out \(\kappa=8/3\).

Two facts about the composition action were known. Fan and Sung \cite{FS}
proved, by the welding of quantum disks, that for \(\kappa<4\) the welding
law is quasi-invariant under Weil--Petersson maps, without a formula for
the density. Baverez and Jego \cite[Theorem~1.2, Corollary~1.3]{BJweld}
proved an integration-by-parts formula for the welding law and deduced
explicit densities for analytic maps; formally the welding law has
density \(e^{-H_\kappa}\) with respect to a non-existent invariant measure,
where \(H_\kappa\) is the action of Theorem~\ref{thm:A} \cite[(1.10),
(1.14)]{BJweld}. In a Gaussian space, the analogous statement is the
Cameron--Martin theorem, and there the relative entropy of a translate is
half the squared Cameron--Martin norm of the shift. Theorem~\ref{thm:A} is
the corresponding statement for SLE weldings: the entropy is the value of
the formal log-density at the deforming map.

\begin{theorem}[Entropy--action identity]\label{thm:A}\label{thm:introEntropy}
Let \(0<\kappa\le4\). Let \(\eta\) be a loop drawn from the shape law of
the SLE\(_\kappa\) loop measure, that is, the loop measure restricted to
loops separating \(0\) from \(\infty\) and normalized to interior conformal
radius one (Section~\ref{sec:setting}). Let \(f:\D\to\operatorname{int}\eta\)
and \(g^+:\Ds\to\operatorname{ext}\eta\) be the conformal maps with \(f(0)=0\),
\(f'(0)>0\), \(g^+(\infty)=\infty\), \((g^+)'(\infty)>0\), let \(\alpha\) be an
independent uniform angle, and let \(P=P_\kappa\) be the law of the welding
\[
 \sigma=(g^+\circ R_\alpha)^{-1}\circ f|_{S^1}\in\operatorname{Homeo}_+(S^1),
 \qquad R_\alpha(z)=e^{i\alpha}z .
\]
For a circle homeomorphism \(a\), let \(P_a\) be the law of \(a\circ\sigma\).
Let \(a\) be a Weil--Petersson homeomorphism: \(a\) is absolutely continuous
with \(\log a'\in H^{1/2}(S^1)\); equivalently, \(a=g_a^{-1}\circ f_a|_{S^1}\) is
the welding of a Jordan curve \(\gamma_a\) of finite Loewner energy
\(I^L(\gamma_a)\), with \(f_a(0)=0\), \(g_a(\infty)=\infty\). Then
\[
 \operatorname{Ent}(P_a\,|\,P)=\operatorname{Ent}(P\,|\,P_a)=H_\kappa(a)
 :=\frac{\cL}{24}\,I^L(\gamma_a)+2K(a),
\]
where
\[
 \cL=13+\frac{24}{\kappa}+\frac{3\kappa}{2},\qquad
 K(a)=\log\left|\frac{g_a'(\infty)}{f_a'(0)}\right|,\qquad
 \pi I^L(\gamma_a)=\int_\D\Bigl|\frac{f_a''}{f_a'}\Bigr|^2\dd
 +\int_{\Ds}\Bigl|\frac{g_a''}{g_a'}\Bigr|^2\dd-4\pi K(a),
\]
and \(\operatorname{Ent}(Q|R)=\int\log(dQ/dR)\,dQ\) is the relative entropy
(\(+\infty\) unless \(Q\ll R\)). The same identity holds for the law of
\(\sigma\circ a\). In particular \(P_a\) and \(P\) are mutually absolutely
continuous, also at \(\kappa=4\).
\end{theorem}

Here \(\pi I^L\) is the universal Liouville action \(S_1\) of
Takhtajan--Teo, which equals \(\pi\) times the Loewner energy by
\cite{WangInv}, and \(K\) is the log ratio of the exterior and interior
conformal radii of \(\gamma_a\), the Velling--Kirillov potential; both are
unchanged when \((f_a,g_a)\) is replaced by \((cf_a,cg_a)\), so \(H_\kappa\) is
well defined, and \(H_\kappa(\mathrm{id})=0\). Three special cases show what
the identity says. For the disc automorphisms
\(m_t(z)=(z+\bar t)/(1+tz)\) the curve is a circle, the Loewner energy
vanishes, and \(\operatorname{Ent}(P_{m_t}|P)=-2\log(1-|t|^2)\) for every
\(\kappa\). As \(\kappa\to0\), \(\kappa\operatorname{Ent}(P_a|P)\to I^L(\gamma_a)\),
in agreement with the large deviations of SLE\(_{0+}\), whose rate
function is the Loewner energy \cite{WangJEMS,PeltolaWang}. And to second
order along the flow of \(\cos n\theta\,\partial_\theta\) the entropy is
\(\frac{s^2}{4}C_n\), \(C_n=2n+\frac{\cL}{12}(n^3-n)\), which is the Fisher
information computed in Theorem~\ref{thm:Akahler}. Theorem~\ref{thm:A} is
Theorem~\ref{thm:entropyWP}.

The proof is a mean-value argument. For maps \(a\) generated by finitely
many Fourier modes the density \(dP_a/dP\) is bounded, and its logarithm is
a renormalized increment of \(H_\kappa\): on the random welding \(\sigma\)
itself the action is infinite, since SLE loops are not Weil--Petersson, and
the increment is defined as a limit along the radial approximations
\(f(r\,\cdot)\) of the loop. The pointwise formula for the log-density
therefore holds only on smooth weldings, a \(P\)-null set, and its
expectation is not visible from it. We show that the regularized increment
is harmonic in the complex parameter \(\zeta\), \(|\zeta|<1\), of the family
\(\zeta^{-1}f(\zeta\,\cdot)\), which is a holomorphic disc in the universal
Teichm\"uller curve through the round circle (\(\zeta=0\)); rotation
invariance of \(P\) turns the expectation at radius \(r\) into a circle
average, and the mean-value property identifies it with the value at
\(\zeta=0\), which is \(H_\kappa(a)\) (Proposition~\ref{prop:entropyFin}). A
two-sided approximation, using lower semicontinuity of the entropy and
the Donsker--Varadhan inequality, extends the identity to all
Weil--Petersson maps.

The relation to earlier and concurrent work is as follows. The identity
gives a proof of the quasi-invariance theorem of Fan--Sung \cite{FS} by a
different route (from \cite{FS} we use only the deterministic Lemma~4.1),
valid also at \(\kappa=4\), and it makes rigorous, at the level of entropy,
the formal density of \cite{BJweld}. The finite-mode densities, written as
radial limits of action increments, and the radial approximations
themselves coincide with those contained in the preprint of Cai
\cite[(1.28)--(1.35)]{CaiResponse}; relative entropy is not considered
there, and our proofs do not use that preprint (Remark~\ref{rem:entCai}).
Independently and concurrently, S.~Fan \cite{FanCM} (v1 posted 18 November
2026, v2 19 November 2026; we cite the numbering of v1) proved a
Cameron--Martin theorem for the same law by different methods: for
\(0<\kappa\le4\) and Weil--Petersson maps \(\psi,\phi\) the law of
\(\psi^{-1}\circ\sigma\circ\phi\) is equivalent to \(P\), with density
\(e^{-J_{\psi,\phi}}\) given by a relative action regularized through Grunsky
operators, and if one of the maps is not Weil--Petersson the laws are
mutually singular \cite[Theorems~1.4, 1.9, 1.10]{FanCM}. Quasi-invariance at
\(\kappa=4\) is thus obtained in both works, and \cite{FanCM} has in
addition the singularity statement and the two-sided action. Relative
entropy is not considered in \cite{FanCM}. In Fan's notation,
Theorem~\ref{thm:A} says that the random relative action has mean equal to
its value at the round circle, which does not follow from the density
formula, for the reason given above
(Proposition~\ref{prop:FanComparison}, Remark~\ref{rem:FanTests}).
Conditionally on Fan's singularity theorem, Theorem~\ref{thm:A} holds for
every \(a\in\operatorname{Homeo}_+(S^1)\), both sides being infinite off the
Weil--Petersson class (Proposition~\ref{prop:FanComparison}(iv)). Both works
take \(\cL\) from the same normalization of the restriction covariance of
the loop measure, so the agreement of the constants is a consistency
check and not an independent determination (Remark~\ref{rem:cLtest}).

The second theorem interprets this identity in Kähler geometry. The
deformed laws form a statistical model indexed by the Weil--Petersson
maps modulo rotations, and such a model carries a Fisher information
metric. That a Kähler metric can arise as a Fisher information metric is
known in finite dimensions: the Bergman metric of a bounded domain is the
Fisher metric of a natural statistical model \cite{ChoYum}. The parameter
space is the Weil--Petersson universal
Teichm\"uller curve of Takhtajan--Teo \cite{TTI,TT}, a complex Hilbert
manifold with two natural right-invariant Kähler metrics, the
Velling--Kirillov metric \(g_{\rm VK}\), with global potential \(K\), and the
pull-back of the Weil--Petersson metric \(g_{\rm WP}\), with global potential
the universal Liouville action. Mixed derivatives are linear, so
\(H_\kappa=2K+\frac{\cL}{24}I^L\) is a global Kähler potential of
\(2g_{\rm VK}+\frac{\cL}6\pi^*g_{\rm WP}\); this is geometry, known before, and
involves no probability. On the smooth part
\(S^1\backslash\operatorname{Diff}_+(S^1)\) the combinations of the two metrics are the
Kirillov--Yuriev metrics \(g_{c,h}\) \cite{KirillovYuriev}, labelled by a
central charge \(c\) and a weight \(h\); they are the invariant Kähler
metrics of this homogeneous space, which is a coadjoint orbit of the
Virasoro group. Their curvature is where the number \(26\) enters: Bowick
and Rajeev \cite{BowickRajeev} computed the regularized Ricci form of every
\(g_{c,h}\) to be \(-\omega_{26,1}\), so that \(g_{c,h}\) is Einstein in this sense
exactly when \(c=26h\), and on the base the Weil--Petersson metric
\(g_{c,0}\) is Kähler--Einstein for every \(c\) \cite{TTI}. On the
probabilistic side, Cai \cite[(1.24), (1.36)]{CaiResponse} computed the
covariance of the score of the welding law at the round circle, which is
\(n+\frac{\cL}{24}(n^3-n)\) on \(\cos n\theta\) and has the Fourier form of
\(g_{\cL,1}\). Combining the two, one could predict which geometric
potential belongs to the Fisher metric. The relative entropy does not
appear in either statement, and Theorem~\ref{thm:Akahler} is what
Theorem~\ref{thm:A} adds to them: the potential is the entropy. Whether
the Fisher metric is the Einstein member of the family \(g_{c,h}\) is
answered after Corollary~\ref{cor:introCons}.

\begin{theorem}[Fisher metric of the welding experiment]\label{thm:Akahler}
Let \(0<\kappa\le4\), and let \(P\), \(P_a\), \(\cL\) and \(H_\kappa\) be as in
Theorem~\ref{thm:A}. Let
\(\mathcal M=\operatorname{Rot}(S^1)\backslash\operatorname{WP}(S^1)\) be the space of
Weil--Petersson homeomorphisms modulo post-composition by rotations,
with the complex Hilbert manifold structure of the Weil--Petersson
Teichm\"uller curve \(\mathcal T_0(1)\) of Takhtajan--Teo
(Section~\ref{sec:markedcomparison}), and for \(h\in\mathcal M\) let
\(E_h=P_{h^{-1}}\) be the law of \(h^{-1}\circ\sigma\). For a real vector
field \(v\,\partial_\theta\) on the circle write
\[
 v=c_0+\sum_{n\ge1}\bigl(it_ne^{in\theta}-i\bar t_ne^{-in\theta}\bigr),\qquad
 \|v\|_F^2=2\sum_{n\ge1}C_n|t_n|^2,\qquad
 C_n=2n+\frac{\cL}{12}(n^3-n).
\]
\begin{enumerate}
\item \emph{Finite-mode submodels.} For every \(h\in\mathcal M\) and every flow
\(\Phi_s\) of a real trigonometric polynomial field \(v\), the laws
\(E_{\Phi_s\circ h}\) have densities with respect to \(E_h\) bounded above and
below; the family is differentiable in quadratic mean at \(s=0\), i.e.\
\(\sqrt{dE_{\Phi_s\circ h}/dE_h}=1+\frac s2\ell+o_{L^2(E_h)}(s)\), with Fisher
information \(\E_{E_h}[\ell^2]=\|v\|_F^2\); and
\(\operatorname{Ent}(E_{\Phi_s\circ h}|E_h)=\frac{s^2}2\|v\|_F^2+O(s^3)\).
\item \emph{Unique continuous extension.} These Fisher informations, as
quadratic forms on the velocities \(X=D\mathsf R_h(v)\) of
\(s\mapsto\Phi_s\circ h\), extend uniquely to a continuous Riemannian metric \(I\)
on \(\mathcal M\): \(I_h(X,X)=\|X\circ h^{-1}\|_F^2\), where \(X\circ h^{-1}=(D\mathsf R_h)^{-1}X\) is
the tangent vector transported to the identity by the right translation
\(\mathsf R_h(k)=k\circ h\). The metric \(I\) is the right-invariant Kähler metric
\(2g_{\rm VK}+\frac{\cL}{6}\pi^*g_{\rm WP}\), where \(\pi:\mathcal M\to
\operatorname{M\ddot ob}(S^1)\backslash\operatorname{WP}(S^1)\) is the projection and, at
the identity, \(g_{\rm VK}(v,v)=2\sum n|t_n|^2\) and
\(\pi^*g_{\rm WP}(v,v)=\sum(n^3-n)|t_n|^2\). It is the right-invariant
extension of the Kirillov--Yuriev metric \(g_{\cL,1}\) of
\(S^1\backslash\operatorname{Diff}_+(S^1)\), and coincides with it at the identity:
\(I(u_n,u_n)=C_n/2\) for \(u_n=\cos n\theta\).
\item \emph{Entropy.} The relative entropy
\(V(h)=\operatorname{Ent}(E_h|P)=H_\kappa(h)\) is a global Kähler potential of
\(I\), that is, \(i\partial\bar\partial V=I(J\cdot,\cdot)\) on \(\mathcal M\); and for
every \(h\in\mathcal M\), \(X\in T_h\mathcal M\) and every \(C^2\) path \(k_s\) with
\(k_0=h\), \(\dot k_0=X\),
\(\operatorname{Ent}(E_{k_s}|E_h)=\frac{s^2}2I_h(X,X)+o(s^2)\).
\end{enumerate}
\end{theorem}

Parts (i)--(iii) are parts (ii)--(iv) of Theorem~\ref{thm:globalKahler},
with Lemma~\ref{lem:conventions}. We regard Theorem~\ref{thm:Akahler} as an
interpretation of Theorem~\ref{thm:A} rather than as an independent
result in geometry, and it is important to separate what in it was
known. Two statements should be compared:
\[
 \text{(A)}\quad i\partial\bar\partial H_\kappa=\omega_I ,\qquad\qquad
 \text{(B)}\quad \operatorname{Ent}(E_h\,|\,P)=H_\kappa(h)\quad(h\in\mathcal M).
\]
Statement (A) says that a geometric function has a certain metric as
its complex Hessian. It is due to Takhtajan--Teo \cite{TTI,TT}, together
with the right invariance of the metric
\(2g_{\rm VK}+\frac{\cL}{6}\pi^*g_{\rm WP}\); and that this metric is the one given
by the score covariance at the round circle is in Cai's preprint, where
Kähler compatibility is also proved at a general Weil--Petersson base
point, at the tangent level and in period-mapping coordinates
\cite[Section~2.3, Propositions~2.10--2.11]{CaiResponse}. That a
Kirillov--Yuriev form appears there is to be expected:
\(2hn+\frac{c}{12}(n^3-n)\) is the Shapovalov norm of \(L_{-n}v\) for a
highest-weight vector \(v\) of weight \(h\) of the Virasoro algebra with
central charge \(c\), and measures whose integration by parts realizes
such forms go back to the unitarizing measures of Airault--Malliavin
\cite{AiraultMalliavin}; see also \cite{GQW} and the Shapovalov form of
\cite{BJloop}. Statement (B) says that the entropy between actual
probability laws, computed, is that geometric function. It is
Theorem~\ref{thm:A}, and it is new.

Statement (B) does not follow from (A), nor from (A) together with the
Fisher metric. Two functions with the same complex Hessian need not be
equal: on \(\mathbb C\), \(|z|^2\) and \(|z|^2+\Rea z^2\) are potentials of the
same Kähler metric. So even if one knew in advance that the relative
entropy is some potential of \(I\), it would be determined only up to a
pluriharmonic function; and the Fisher metric gives less than that, since
it is an infinitesimal quantity, the covariance of the score at one
point, and many divergences have the same Hessian. Concretely, along the
flow \(\psi_s\) of a trigonometric field \(v\), with score density \(\rho_v\)
(Section~\ref{sec:entropy}),
\[
 \frac{d}{ds}\operatorname{Ent}(P^{\mathsf R}_{\psi_s}\,|\,P)=\E\bigl[\rho_v(\sigma\circ\psi_s)\bigr],
 \qquad P^{\mathsf R}_{\psi_s}=\operatorname{Law}(\sigma\circ\psi_s),
\]
the mean of a stress coefficient of the deformed welding. It vanishes at
\(s=0\), and its derivative there is the Fisher information. The Fisher
metric therefore gives the second-order Taylor coefficient of the entropy
at each base point, while (B) requires the first moment
\(\E[\rho_v(\sigma\circ\psi_s)]\) for every \(s\), and the entropy on the whole
Weil--Petersson class. This first moment is not a covariance at the round
circle. The finite-mode densities of \cite{CaiResponse} and the densities of
\cite{FanCM} describe finite transformations, but they are random
variables, and no expectation of their logarithm is computed there. That
expectation is what the mean-value argument provides.

Accordingly, the parts of Theorem~\ref{thm:Akahler} depend on
Theorem~\ref{thm:A} to different degrees. The score covariance at the
round circle, the finite-mode Fisher informations of part (i), the
identification with the Kirillov--Yuriev metric, the Kähler property of
the metric, and the fact that \(H_\kappa\) is a global Kähler potential of it
do not need Theorem~\ref{thm:A}; they use the finite-flow likelihoods of
Section~\ref{sec:finiteinput} and the Takhtajan--Teo geometry. That the
\emph{relative entropy} is a global Kähler potential needs
Theorem~\ref{thm:A}, and so does, in our proof, the expansion of the
entropy in every Weil--Petersson direction in part (iii). In short, what
was known is that \(i\partial\bar\partial\) of \(H_\kappa\) is a Kähler metric; what is
proved here is that \(H_\kappa\) is the relative entropy, and hence that
\(i\partial\bar\partial\) of the relative entropy is the Fisher metric. Among all
Kähler potentials of \(I\), the statistics of the welding family select
this one: it is defined on all of \(\mathcal M\), without regularization, is
nonnegative, and vanishes only at the round circle.

For a fixed finite Fourier field, differentiability in quadratic mean
follows from bounded response derivatives by Taylor expansion, and we do
not count it, or the Fourier covariance, as new. Fisher information along
general Weil--Petersson paths, at the critical \(H^{3/2}\) regularity, is
not treated here (Remark~\ref{rem:criticalPaths}). Cai's reversible
response diffusion \cite{CaiDiffusion,CaiWeak} concerns a stochastic time
evolution with invariant law \(P\); our parameter paths are deterministic
changes of the law.

\begin{corollary}\label{cor:introCons}
In the setting of Theorems~\ref{thm:A} and~\ref{thm:Akahler}:
\begin{enumerate}
\item For Weil--Petersson maps \(a,b\), with \(\|Q-R\|_{\rm TV}=\sup_B|Q(B)-R(B)|\),
\[
 \operatorname{Ent}(P_a|P_b)=H_\kappa(b^{-1}\circ a),\qquad
 \|P_a-P_b\|_{\rm TV}\le\sqrt{H_\kappa(b^{-1}\circ a)/2}.
\]
\item For the disc automorphisms \(m_t(z)=(z+\bar t)/(1+tz)\), \(|t|<1\) (for
which \(I^L=0\) and \(K(m_t)=-\log(1-|t|^2)\)),
\(\operatorname{Ent}(P_{m_t}|P)=-2\log(1-|t|^2)\) for every \(\kappa\).
\item \(\kappa\operatorname{Ent}(P_a|P)\to I^L(\gamma_a)\) as \(\kappa\to0\), for every
Weil--Petersson map \(a\).
\item \(I=g_{26,1}-c_{\rm m}\,g_{1,0}\), where
\(c_{\rm m}=26-\cL=(6-\kappa)(3\kappa-8)/(2\kappa)\) is the central charge of
SLE\(_\kappa\), \(g_{c,h}\) denotes the right-invariant metric with value
\(hn+\frac{c}{24}(n^3-n)\) on \((u_n,u_n)\) at the identity, and \(g_{1,0}\) is
the mixed Hessian of \(\pi^*I^L/24\). In particular \(I=g_{26,1}\) if and only
if \(\kappa=8/3\).
\end{enumerate}
\end{corollary}

Corollary~\ref{cor:introCons} collects Corollaries~\ref{cor:entropyConsequences}
and~\ref{cor:KE}. Part (iv) relates the result to the Onsager--Machlup
functional of Carfagnini and Wang \cite[Theorem~1.1]{CarfagniniWang}, who
showed that for fixed \(\kappa\le4\) the SLE\(_\kappa\) loop measure has
Onsager--Machlup functional \(\frac{c_{\rm m}}{24}I^L\), with the matter
charge \(c_{\rm m}\). The Fisher potential \(H_\kappa=\frac{\cL}{24}I^L+2K\) carries
the Liouville charge \(\cL\) instead, and the two charges sum to \(26\), as in
the matter--Liouville--ghost factorization of \cite{BJweld}. The two
functionals answer different questions: the Onsager--Machlup functional
compares the masses of small neighbourhoods of two loops under one
measure, the entropy compares two laws of the welding. Part (iv) also answers the Einstein question raised before
Theorem~\ref{thm:Akahler} (Remark~\ref{rem:regRicci}). The horizontal part \(\frac{\cL}{6}\pi^*g_{\rm WP}\) of
the Fisher metric is a multiple of the Kähler--Einstein metric of
Takhtajan--Teo, with Einstein constant \(-26/\cL\). For the regularized
Ricci form \(-\omega_{26,1}\) of Bowick--Rajeev \cite{BowickRajeev}, the Fisher
metric itself is Kähler--Einstein exactly when \(\cm=0\), that is at
\(\kappa=8/3\), where the loop measure satisfies exact conformal restriction;
for other \(\kappa\) it satisfies a twisted Kähler--Einstein equation whose
twist has the Onsager--Machlup functional as potential. These are
readings of known curvature computations through the entropy, not new
curvature theorems.

Theorem~\ref{thm:A} also gives the action a meaning at fixed \(\kappa\)
(Corollary~\ref{cor:entropyFurther}). By the Chernoff--Stein lemma,
\(H_\kappa(a)\) is the optimal exponential rate at which \(n\) independent
weldings distinguish \(P\) from \(P_a\). The Loewner energy and the capacity
term are recovered separately from the entropies at two values of
\(\kappa\), since \(H_\kappa\) is affine in \(\cL\). The law \(P\) is preserved by a
Weil--Petersson map only if it is a rotation; the entropy cost of a
deformation decreases with \(\kappa\) and is smallest for SLE\(_4\); and the
Donsker--Varadhan formula turns Theorem~\ref{thm:A} into lower bounds for
\(\frac{\cL}{24}I^L+2K\) by expectations of SLE observables.

The probabilistic input of both theorems is a horizontal
integration-by-parts formula. Let \(\mathfrak t_n\) be the interior stress
coefficients, \(\mathcal C\) the cylinder core and \(D_v\) the right-flow
derivative of Section~\ref{sec:finiteinput}. Then
\begin{equation}\tag{W}\label{eq:W}
 \int D_{\cos n\theta}F\,dP=\int F\,\Ima\mathfrak t_n\,dP,\qquad
 \int D_{\sin n\theta}F\,dP=\int F\,\Rea\mathfrak t_n\,dP
 \qquad(F\in\mathcal C,\ n\ge1).
\end{equation}
This formula is not new. For trigonometric fields it is the real form of
the adjoint formula of Baverez--Jego \cite[Theorem~1.2]{BJweld},
transported by the inversion \(\mathsf J(h)=h^{-1}\), and Cai states it on the
same cylinder core \cite[Proposition~4.4]{CaiResponse}. Both are unrefereed
preprints, and since every result here rests on \eqref{eq:W} we give a
proof from published results in Section~\ref{app:Wproof}
(Theorem~\ref{thm:Wproof}), independent of both: a Ward identity for Zhan's
SLE loop measure \cite{ZhanLoop}, obtained from the conformal restriction
covariance and the first-variation formula of Gordina--Qian--Wang
\cite{GQW}, is transferred to welding coordinates, where the change of
basis contributes the ghost charge \(26\) and turns the matter charge
\(c_{\rm m}\) into \(\cL=26-c_{\rm m}\). The constant of the restriction covariance
used in the literature is verified in Lemma~\ref{lem:normalization}, and
the first-variation formula is proved for all Jordan curves
(Lemmas~\ref{lem:GQWanalytic} and~\ref{lem:GQWext}). No result of this
paper depends on \cite{CaiResponse}; the coefficient responses are proved
in Appendix~\ref{app:coeffResp}. From the preprint \cite{BJweld} we use only
the definition of \(P\). At \(\kappa=4\) we use in addition the conformal
removability of SLE\(_4\) \cite{KMS}, also a preprint. Readings of the
conventions of \cite{BJweld,GQW,ZhanLoop,TTI} are collected in
Section~\ref{sec:conventionsLit}.

\begin{remark}[Scope of the claims]\label{rem:scope}
We do not assert likelihood regularity beyond finite-mode flows, or a
\(\kappa\)-score. All identities concern the deformation parameter at
fixed \(\kappa\); for different \(\kappa\) the laws are mutually singular,
since the curves have almost sure dimension \(1+\kappa/8\)
\cite{Beffara}, and for \(0<\kappa\le4\) the curve is conformally removable
(by \cite{RohdeSchramm,JonesSmirnov} for \(\kappa<4\) and \cite{KMS} for \(\kappa=4\)), hence determined
by its welding up to M\"obius maps, so the laws for different \(\kappa\) are
carried by disjoint sets of weldings. The asymptotic statistics of the
welding experiment are not treated here.
\end{remark}

Among further related work, loops of finite Loewner energy are exactly
the Weil--Petersson quasicircles \cite{WangInv}, with formulas for the
energy in terms of the welding alone in \cite{FVW}; see the survey
\cite{Wang}. The loop measures are constructed in \cite{ZhanLoop} and
characterized in \cite{BJloop}, and the welding of quantum disks
\cite{AHSloop} underlies \cite{FS}. The Kähler geometry of
\(S^1\backslash\operatorname{Diff}_+(S^1)\) and of the Weil--Petersson Teichm\"uller
space is developed in
\cite{KirillovYuriev,BowickRajeev,TT,TTI,GordinaLescot}; we also use the
deterministic flow and chart results of \cite{Shen,ShenTang,WHS}.

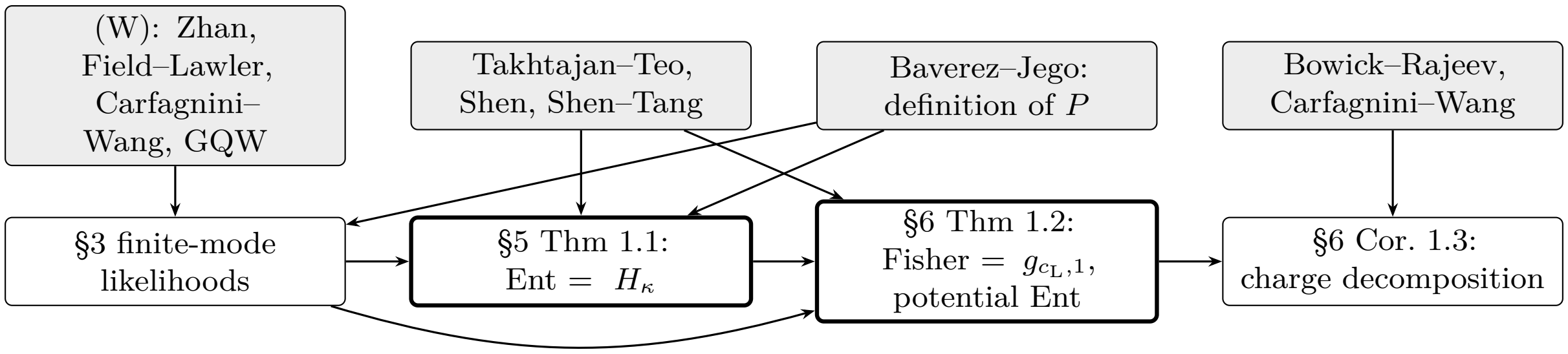
\begin{figure}[ht]
\centering
\begin{tikzpicture}[
  every node/.style={font=\footnotesize,align=center,text width=3.3cm,
                     minimum height=9mm,inner sep=2.5pt},
  inp/.style={draw,rounded corners=2pt,fill=black!7},
  res/.style={draw,rounded corners=2pt},
  main/.style={draw,very thick,rounded corners=2pt},
  ar/.style={-{Stealth[length=4pt]},semithick}]
\def\dx{4.15}
\node[inp] (cai)  at (0,0)        {\eqref{eq:W}: Zhan, Field--Lawler,\\Carfagnini--Wang, GQW};
\node[inp] (tt)   at (\dx,0)      {Takhtajan--Teo,\\Shen, Shen--Tang};
\node[inp] (bj)   at (2*\dx,0)    {Baverez--Jego:\\definition of \(P\)};
\node[inp] (cw)   at (3*\dx,0)    {Bowick--Rajeev,\\Carfagnini--Wang};
\node[res] (fin)  at (0,-1.8)     {\S\ref{sec:finiteinput} finite-mode\\likelihoods};
\node[main] (ent) at (\dx,-1.8)   {\S\ref{sec:entropy} Thm~\ref{thm:A}:\\Ent \(=H_\kappa\)};
\node[main] (kah) at (2*\dx,-1.8) {\S\ref{sec:global} Thm~\ref{thm:Akahler}:\\Fisher \(=g_{\cL,1}\), potential Ent};
\node[res] (ke)  at (3*\dx,-1.8) {\S\ref{sec:KE} Cor.~\ref{cor:introCons}:\\charge decomposition};
\draw[ar] (cai)--(fin); \draw[ar] (bj)--(fin); \draw[ar] (fin)--(ent); \draw[ar] (tt)--(ent);
\draw[ar] (bj)--(ent); \draw[ar] (ent)--(kah); \draw[ar] (kah)--(ke);
\draw[ar] (cw)--(ke);
\draw[ar] (tt)--(kah); \draw[ar] (fin) to[bend right=16] (kah);
\end{tikzpicture}
\caption{Dependence of the main results. Shaded boxes are external
inputs; thick boxes are the main theorems.}\label{fig:roadmap}
\end{figure}

The paper is organized as follows; Figure~\ref{fig:roadmap} shows the
dependencies. Section~\ref{sec:setting} fixes the welding law and
Section~\ref{sec:finiteinput} derives the finite-flow likelihoods from
\eqref{eq:W}, which is proved in Section~\ref{app:Wproof}.
Section~\ref{sec:entropy} proves Theorem~\ref{thm:A} and compares it with
\cite{FanCM}. Section~\ref{sec:global} proves Theorem~\ref{thm:Akahler} and
the charge decomposition, and Section~\ref{sec:open} collects open
problems. The appendices contain the M\"obius density and coefficient
moments, the metric calibration, the constant of the restriction
covariance, and the first-variation formula of \cite{GQW} with the score
covariance (Appendices~\ref{app:mobius}--\ref{app:caiInput}). The deterministic
lemmas on which the mean-value argument of Section~\ref{sec:entropy} rests,
in particular the continuity, the holomorphic dependence on the disc
parameter and the derivative of \(H_\kappa\) along flows, are proved in
Appendix~\ref{sec:entropyDet}, and the coefficient responses in
Appendix~\ref{app:coeffResp}. The paper is self-contained in this sense:
there is no supplementary material apart from the scripts used for
numerical and computer-algebra checks.

Table~\ref{tab:notation} lists the main symbols.
Where a letter has more than one use, the table lists each use with the
places where it occurs.

\begin{table}[ht]
\centering\footnotesize
\renewcommand{\arraystretch}{1.0}\setlength{\tabcolsep}{4pt}
\begin{tabular}{@{}p{3.6cm}p{9.6cm}p{2.3cm}@{}}
\toprule
Symbol & Meaning & Where\\
\midrule
\(X\), \(L_a\), \(\mathsf R_a\), \(\mathsf J\), \(R_\alpha\)
 & \(\operatorname{Homeo}_+(S^1)\); \(x\mapsto a\circ x\); \(x\mapsto x\circ a\); \(x\mapsto x^{-1}\); rotation \(z\mapsto e^{i\alpha}z\) & \S\ref{sec:setting}\\
\(P\), \(P_a\), \(P^{\mathsf R}_\Phi\), \(E_h\)
 & welding law; \((L_a)_*P\); \((\mathsf R_\Phi)_*P\); \(P_{h^{-1}}\), \(h\in\mathcal M\) & \S\S\ref{sec:setting}, \ref{sec:entropy}, \ref{sec:global}\\
\(P_p\) & Schwarzian coefficient \([w^{p-2}]\mathcal S(f^{-1})\); never the law \(P\) & \S\ref{app:Wproof}, App.~\ref{app:caiInput}\\
\(f,g\); \(G\) & interior/exterior uniformizers; reciprocal exterior map \eqref{eq:G} & \S\ref{sec:setting}\\
\(\alpha_m,\beta_m\); \(a_m\) & exterior coefficients (of \(G\), \(g\)); interior coefficients of \(f\) & \S\ref{sec:setting}; \S\S\ref{app:Wproof}, \ref{sec:entropy}\\
\(a\), \(b\), \(\Phi\), \(\Psi\) & deterministic circle maps (a WP map \(a\) is never a coefficient) & throughout\\
\(\alpha\) (no index) & uniform rotation angle in the definition of \(P\) & \S\ref{sec:setting}, \S\ref{app:Wproof}\\
\(s_n\); \(\mathfrak t_n\) & exterior and interior stress coefficients \eqref{eq:stress} & \S\ref{sec:finiteinput}\\
\(\cL\); \(\cm=c=c(\kappa)\) & Liouville charge; matter (central) charge, \(\cm=26-\cL\le1\) & \S\ref{sec:setting}; \S\ref{app:Wproof}, App.~\ref{app:normalization}\\
\(c(\gamma)\), \(\eta(\gamma)\) & log conformal radius \(\log f_\gamma'(0)\); shape \(e^{-c}\gamma\) & \S\S\ref{sec:setting}, \ref{app:Wproof}\\
\(C_n\); \(C\) & \(2n+\frac{\cL}{12}(n^3-n)\); \(C=C_2\) in \S\S\ref{sec:setting}--\ref{sec:finiteinput} and App.~\ref{app:mobius}, otherwise a generic constant (also \(C'\)) & \\
\(K\); \(K_t(\delta)\) & \(\log|g'(\infty)/f'(0)|\); compact sets \(\widehat{\mathbb C}\setminus\phi_t(D_\delta)\); a plain compact set is also \(K\) in \eqref{eq:LambdaNest} & \S\ref{sec:entropy}; App.~\ref{app:caiInput}\\
\(S_1\), \(I^L=S_1/\pi\), \(H_\kappa\) & universal Liouville action, Loewner energy, \(\frac{\cL}{24}I^L+2K\) & \S\ref{sec:entropy}\\
\(H(\cdot,\cdot)\); \(H\) & Hellinger distance; a bounded test function on loops & \S\ref{sec:setting}; \S\ref{app:Wproof}\\
\(\Lambda^*\), \(\mathcal B\) & Field--Lawler normalized loop mass; loop mass in a domain & Apps.~\ref{app:normalization}, \ref{app:caiInput}\\
\(\Lambda_pf\); \(\Lambda_\Psi\), \(\Lambda^{\rm e}\); \(\mathfrak l\) & Ward-flow variation \eqref{eq:Wp-LPhi}; radial extension and its exterior part; \([h]\mapsto[\log h']\) & \S\ref{app:Wproof}; App.~\ref{sec:entropyDet}; \S\ref{sec:markedcomparison}\\
\(\mathcal C\), \(\mathcal W\) & cylinder core; weldings of removable curves & \S\ref{sec:finiteinput}\\
\(D_v\), \(\mathfrak d_j\), \(L_j\), \(\mathbf L_j\) & right-flow derivative; its complex modes; coefficient derivations and their generators & \S\S\ref{sec:finiteinput}, \ref{app:Wproof}, App.~\ref{app:coeffResp}\\
\(\nu=\nu_\kappa\), \(\mu_\kappa\), \(\mu^0\), \(\mu^1_z\), \(\mu^{\rm loop}\) & SLE loop measure; shape law; Zhan's unrooted/rooted loop measures; Brownian loop measure & \S\ref{sec:setting}, App.~\ref{app:normalization}\\
\(\varpi\); \(\mu_F\) & Beltrami coefficients (variable); coefficient of a map \(F\) & \S\ref{app:Wproof}, Apps.\\
\(\mathsf S\), \(\|\cdot\|_F\), \(I\) & score of a finite Fourier field, Fisher norm \eqref{eq:Fishernorm}, Fisher metric (Theorem~\ref{thm:globalKahler}) & \S\ref{sec:global}\\
\(\mathcal M\), \(\pi\), \(g_{\rm VK}\), \(g_{\rm WP}\), \(g_{c,h}\) & marked WP space; forgetful map; Velling--Kirillov, WP and Kirillov--Yuriev metrics & \S\ref{sec:global}\\
\(m_t\); \(\iota\), \(\mathsf c\) & \((z+\bar t)/(1+tz)\); \(z\mapsto1/\bar z\), \(z\mapsto\bar z\) & \S\ref{sec:entropy}; \S\ref{sec:setting}\\
\bottomrule
\end{tabular}
\caption{Main notation.}\label{tab:notation}
\end{table}
\section{The welding experiment and its parameters}
\label{sec:setting}\label{sec:background}

Put \(\D=\{|z|<1\}\) and
\(\Ds=\{|z|>1\}\cup\{\infty\}\). The observation space
\(X=\operatorname{Homeo}_+(S^1)\) is the Polish group of
orientation-preserving circle homeomorphisms, with uniform convergence
of maps and inverses and its Borel sigma field. The maps
\(L_a(x)=a\circ x\), \(\mathsf R_a(x)=x\circ a\), and \(\mathsf J(x)=x^{-1}\)
are Borel bijections. Plain \(R_\alpha\) or \(R_y\) always denotes the
rotation \(z\mapsto e^{i\alpha}z\) of the circle.

We first define the reference probability law \(P\). Fix \(0<\kappa\le4\). Let \(\nu_\kappa\) be the SLE loop
measure restricted to simple loops separating \(0\) and \(\infty\),
with the normalization of \cite[Section~1.2]{BJloop}.
For such a loop \(\gamma\), let \(f_\gamma:\D\to\Omega_\gamma\)
be the interior Riemann map with \(f_\gamma(0)=0\) and
\(f_\gamma'(0)>0\), and put
\[
 c(\gamma)=\log f_\gamma'(0),\qquad
 \eta(\gamma)=e^{-c(\gamma)}\gamma.
\]
The shape \(\eta\) has conformal radius one. If \(\mathcal J_1\)
is the corresponding Borel space of shapes, the scale disintegration is
\[
 (c,\eta)_*\nu_\kappa=\pi^{-1/2}\,dc\otimes\mu_\kappa,
 \qquad \mu_\kappa(\mathcal J_1)=1.
\]
Equivalently, for Borel \(B\subset\mathcal J_1\),
\[
 \mu_\kappa(B)=\sqrt\pi\,
 \nu_\kappa\{\gamma:0\le c(\gamma)<1,\ \eta(\gamma)\in B\}.
\]
The loop measure and its disintegration are external inputs; no
uniform probability on the scale axis \(\mathbb R\) is used.

For \(\eta\sim\mu_\kappa\), choose the interior and exterior maps
\(f,g^+\), with
\[
 f(0)=0,\quad f'(0)=1,\qquad
 g^+(\infty)=\infty,\quad (g^+)'(\infty)>0.
\]
Independently choose \(\alpha\) uniform modulo \(2\pi\), set
\(g=g^+\circ R_\alpha\), and define
\[
 \sigma=g^{-1}\circ f:S^1\to S^1,\qquad
 P=P_\kappa=\operatorname{Law}(\sigma).
\]
This is \cite[Definition 1.1]{BJweld}. The range \(\kappa\le4\) gives simple
loops. Use \(\cL,C_n\) from the introduction, with matter charge
\[
 \cm=26-\cL=\frac{(6-\kappa)(3\kappa-8)}{2\kappa}=13-\frac{24}\kappa-\frac{3\kappa}2,
 \qquad C=C_2 .
\]
Since \(\frac{24}\kappa+\frac{3\kappa}2\ge12\), with equality only at \(\kappa=4\),
we have \(\cm\le1\) and \(\cL\ge25\) for \(0<\kappa\le4\).

We now list the properties of \(P\) that the paper uses, so that later sections can refer to this list. Besides the integration-by-parts formula
\eqref{eq:W} (Theorem~\ref{thm:Wproof}), the paper uses only the following
properties of \(P\): (a) the definition \cite[Definition~1.1]{BJweld} with the
scale disintegration above; (b) conformal removability of SLE\(_\kappa\) loops,
which rests on H\"older-domain estimates \cite{RohdeSchramm,JonesSmirnov} for
\(\kappa<4\) and on \cite{KMS} at \(\kappa=4\); it makes the normalized
uniformizers functions of the welding, so that the coefficient functionals
have Borel versions that are actual coefficients on the full-measure set
\(\mathcal W\) of Section~\ref{sec:finiteinput}; (c) invariance under left and
right rotations, and \(\mathsf J_*P=P\) (Lemma~\ref{lem:symmetries}, from
\cite{ZhanLoop}). Later sections refer to this paragraph rather than
repeating it.

Next we fix the coordinates on weldings in which the scores are computed, and
the elementary Fourier directions. Write the normalized exterior map as
\[
 \frac{g(z)}b=z+\beta_0+\beta_1z^{-1}+\cdots,
 \qquad b=g'(\infty)\ne0.
\]
Its reciprocal coordinate is the normalized univalent map
\begin{equation}\label{eq:G}
 G(z)=\frac b{g(1/z)}=z+\alpha_1z^2+\alpha_2z^3+\cdots\quad(z\in\D).
\end{equation}
Thus \(\alpha_1=-\beta_0\), \(\alpha_2=\beta_0^2-\beta_1\). The
coordinates are invariant under common complex scaling of \(f,g\)
and have sample-independent bounds by the schlicht coefficient theorem.

For \(t\in\mathbb C\), the second-mode boundary field is
\begin{equation}\label{eq:flow}
 u_t(z)=-tz^3+\bar t z^{-1},\qquad |z|=1.
\end{equation}
Its angular velocity is \(-2\Rea t\sin2\theta-2\Ima t\cos2\theta\).
Write \(\Phi_t\) for its time-one flow and \(P_t=(L_{\Phi_t})_*P\).
Only \(s\mapsto\Phi_{s\tau}\), for fixed \(\tau\), is a
one-parameter group. The second-mode score polynomial is
\begin{equation}\label{eq:score}
 S=s_2=C\alpha_2+(3-C)\alpha_1^2=3\beta_0^2-C\beta_1.
\end{equation}
The finite-flow theorem identifies it as a likelihood derivative.

We also fix the statistical conventions and recall the Weil--Petersson
class. For probabilities \(Q,R\), we use
\[
 H^2(Q,R)=\int(\sqrt{dQ/d\varpi}-\sqrt{dR/d\varpi})^2d\varpi,
 \qquad \varpi=Q+R.
\]
For a dominated family \(Q_t\), DQM with score
\(\ell\in L^2_0(Q_0;\mathbb R)\) means
\(\sqrt{dQ_t/dQ_0}=1+t\ell/2+o_{L^2}(|t|)\).
Its Fisher form is the score covariance. The metric induced by the
half-density embedding is one quarter of that form.
For a complex coordinate \(t=x+iy\) we write its linear score
as \(ts+\bar t\bar s\); the real score vector is therefore
\((2\Rea s,-2\Ima s)\). If \(\E s^2=0\), its real covariance
is \(2\E|s|^2I_2\).

The WP class consists of the quasisymmetric homeomorphisms admitting
a quasiconformal disc extension with square-integrable Beltrami
coefficient for the hyperbolic area. Equivalently \cite{Shen} (see \cite[Theorem~1.1]{ShenTang}, \cite[Theorem~1.1]{WHS}), they
are absolutely continuous and their angular logarithmic derivatives
belong to \(H^{1/2}(S^1;\mathbb R)\). The marked WP space and its standard
complex structure are recalled in Section~\ref{sec:markedcomparison}.

\begin{lemma}[Symmetries of \(P\)]\label{lem:symmetries}
Let \(0<\kappa\le4\) and \(\mathsf J(h)=h^{-1}\).
\begin{enumerate}
\item \(P\) is invariant under left and right rotations, and the shape
law \(\mu_\kappa\) under \(\eta\mapsto e^{i\beta}\eta\).
\item \(\mathsf J_*P=P\).
\end{enumerate}
\end{lemma}
\begin{proof}
(i) \(P\) is invariant under left rotations by the independent uniform
rotation \(\alpha\). Zhan's loop measure is invariant under M\"obius maps
\cite[Theorem~4.2(ii)]{ZhanLoop}, in particular under rotations, hence so is
\(\mu_\kappa\); for \(\eta'=e^{i\beta}\eta\) the normalized maps are
\(e^{i\beta}f(e^{-i\beta}\cdot)\) and \(e^{i\beta}g^+(e^{-i\beta}\cdot)\), so the welding of
\(\eta'\) is \(R_\beta\circ\sigma\circ R_{-\beta}\), and the uniform left rotation
absorbs \(R_\beta\). This gives right-rotation invariance of \(P\).

(ii) For a Jordan curve \(\gamma\) separating \(0\) from \(\infty\), let
\(f_\gamma,g^+_\gamma\) be its normalized interior and exterior maps
(\(f_\gamma(0)=0\), \(f_\gamma'(0)>0\), \(g^+_\gamma(\infty)=\infty\), \((g^+_\gamma)'(\infty)>0\)), and put
\(w(\gamma)=(g^+_\gamma)^{-1}\circ f_\gamma|_{S^1}\). Then \(w(\lambda\gamma)=w(\gamma)\) for
\(\lambda>0\) and \(w(e^{i\beta}\gamma)=R_\beta\circ w(\gamma)\circ R_{-\beta}\). The reflection \(\iota\)
exchanges \(0\) with \(\infty\) and \(\D\) with \(\Ds\). It maps \(\gamma\) to a Jordan curve
separating \(0\) from \(\infty\), exchanging interior and exterior. The maps
\(\iota\circ g^+_\gamma\circ\iota\) and \(\iota\circ f_\gamma\circ\iota\) are conformal onto the
interior and the exterior of \(\iota(\gamma)\), with derivatives \(1/(g^+_\gamma)'(\infty)>0\) at \(0\)
and \(1/f'_\gamma(0)>0\) at \(\infty\). Hence they are its normalized maps. Since \(\iota=\mathrm{id}\)
on \(S^1\), \(w(\iota(\gamma))=w(\gamma)^{-1}\). In the notation of
Section~\ref{sec:setting}, \(c(\iota(\gamma))=-c(\gamma)-K(\eta(\gamma))\), and since
\(\iota(\lambda z)=\lambda^{-1}\iota(z)\),
\(\eta(\iota(\gamma))=\iota^\sharp(\eta(\gamma))\) with \(\iota^\sharp(\eta)=e^{K(\eta)}\iota(\eta)\).
Here \(K(\eta)=\log(g^+_\eta)'(\infty)\) is finite and Borel by the Carath\'eodory kernel
theorem, and \(\iota\) is a homeomorphism for the spherical Hausdorff metric.

Next, \(\iota_*\nu_\kappa=\nu_\kappa\). Write \(\iota=\mathsf c\circ(z\mapsto1/z)\) with
\(\mathsf c(z)=\bar z\). By \cite[Theorem~4.2(ii)]{ZhanLoop} it suffices to show that
\(\mathsf c_*\mu^0=\mu^0\) for Zhan's unrooted loop measure \(\mu^0\). The map \(\mathsf c\)
reverses the winding of oriented loops. Let \(\mathcal R\) denote reversal of
orientation and put \(\hat{\mathsf c}=\mathcal R\circ\mathsf c=\mathsf c\circ\mathcal R\). The
rooted measures \(\mu^1_z\) and \(\mu^0\) are reversible \cite[Theorems~4.1(i)
and~4.2]{ZhanLoop}, so \(\hat{\mathsf c}_*\mu^1_{\bar z}=\mathsf c_*\mu^1_{\bar z}\) and
\(\hat{\mathsf c}_*\mu^0=\mathsf c_*\mu^0\); we apply the uniqueness statement to this
measure. The map \(\mathsf c\) sends chordal SLE\(_\kappa\) in \(D\) from \(a\) to \(b\) to
chordal SLE\(_\kappa\) in \(\mathsf c(D)\) from \(\bar a\) to \(\bar b\), with the same
direction of parametrization: in \(\mathbb H\), \(z\mapsto-\bar z\) turns the Loewner
chain driven by \(W\) into the one driven by \(-W\), and \(-\sqrt\kappa B\) has the law
of \(\sqrt\kappa B\). Being a Euclidean isometry, \(\mathsf c\) preserves diameters and
Minkowski content, \(M_{\mathsf c(\gamma)}=\mathsf c_*M_\gamma\). So for each
\(z\in\mathbb C\), \(\mathsf c_*\mu^1_{\bar z}\) is supported on loops rooted at \(z\),
satisfies the conformal Markov property of \cite[Theorem~4.1(ii)]{ZhanLoop}, and gives
the same mass \(C_1r^{d-2}\) to loops of diameter \(>r\)
\cite[Theorem~4.1(vi)]{ZhanLoop}. Uniqueness \cite[Theorem~4.1(vii)]{ZhanLoop} gives
\(\mathsf c_*\mu^1_{\bar z}=\mu^1_z\) as measures on oriented rooted loops. Applying
\(\mathsf c\) to \(\mu^0(d\gamma)\otimes M_\gamma(dz)=\mu^1_z(d\gamma)\otimes m_2(dz)\)
\cite[(4.6)]{ZhanLoop} gives the same identity for \(\mathsf c_*\mu^0\). Weighting by
\(\mathrm{Cont}(\gamma)^{-1}=M_\gamma(\mathbb C)^{-1}\) and taking first marginals gives
\(\mathsf c_*\mu^0=\mu^0\) on oriented loops, hence on unoriented traces.

Hence, for Borel \(B\subset\mathcal J_1\), the scale disintegration gives
\[
\begin{aligned}
 \pi^{-\frac12}\mu_\kappa(B)&=\nu_\kappa\{0\le c\circ\iota<1,\ \eta\circ\iota\in B\}\\
 &=\pi^{-\frac12}\!\int\!\mathbf 1_B(\iota^\sharp\eta)\,
   \bigl|\{c:0\le-c-K(\eta)<1\}\bigr|\,d\mu_\kappa(\eta)
 =\pi^{-\frac12}\mu_\kappa\bigl((\iota^\sharp)^{-1}B\bigr).
\end{aligned}
\]
So \(\iota^\sharp\) preserves \(\mu_\kappa\), and \(w\circ\iota^\sharp=w^{-1}\). Finally let
\(\sigma=R_{-\alpha}\circ w(\eta)\), with \(\eta\sim\mu_\kappa\) and \(\alpha\) independent and
uniform. Then
\(\sigma^{-1}=w(\iota^\sharp\eta)\circ R_\alpha=R_\alpha\circ w(e^{-i\alpha}\iota^\sharp\eta)\).
Since \(\mu_\kappa\) is rotation invariant by (i), the pair
\((e^{-i\alpha}\iota^\sharp\eta,\alpha)\) has the law of \((\eta,\alpha)\). Hence
\(\sigma^{-1}\) has the law of \(R_\alpha\circ w(\eta)\), which is the law of
\(R_{-\alpha}\circ w(\eta)=\sigma\). That is, \(\mathsf J_*P=P\).
\end{proof}

\phantomsection\label{sec:conventionsLit}%
The section ends with the conventions and constants of the literature.
Several statements that we quote can be read in more than one way, or
use conventions that differ from ours by a constant. We record here, in
one place, the reading we use and why it is the compatible one. None of
our proofs depends on any of these readings: in each case the statement
we need is proved in this paper, from published results and, at
\(\kappa=4\), the removability of SLE\(_4\) loops \cite{KMS}.

\begin{enumerate}
\item[(C1)] \emph{\cite[Theorem~1.2]{BJweld}, multiplier.} We read the
adjoint formula without complex conjugation of the multiplier; in our
notation it is \(\int\mathfrak d_nF\,dP=i\int F\,\bar{\mathfrak t}_n\,dP\)
(the complex form \eqref{eq:Wp-complex} of \eqref{eq:W}). This is the
reading that their proof gives, and the one compatible with
\(\E|\mathfrak t_n|^2=C_n>0\): the other sign would give \(-C_n\)
(Remark~\ref{rem:Wp-checks}(b)). We prove \eqref{eq:W} in
Theorem~\ref{thm:Wproof} without using \cite[Theorem~1.2]{BJweld}.
\item[(C2)] \emph{\cite[Corollary~1.3]{BJweld}, orientation.} With the
sign fixed by \cite[Theorem~1.2, (1.13)]{BJweld} and
\cite[Lemma~5.2]{BJweld}, the integrated identity reads
\(\int F(\Phi\circ h)e^{-\Omega(\Phi\circ h,h)}\,dP(h)=\int F\,dP\). We
therefore read the weight \(e^{-\Omega(\Phi\circ h,h)}\) as the density at
\(h\) of the law of \(\Phi^{-1}\circ h\); equivalently
\(d(L_\Phi)_*P/dP(x)=e^{\Omega(x,\Phi^{-1}\circ x)}\), formally
\(e^{H_\kappa(x)-H_\kappa(\Phi^{-1}\circ x)}\). This is the orientation, and the
sign, of Theorem~\ref{thm:finite}, Lemma~\ref{lem:Bfin} and
Corollary~\ref{cor:entMobius}, which are proved from \eqref{eq:W}.
\item[(C3)]\label{rem:GQWext} \emph{\cite[Corollary~3.9]{GQW}, the set
\(K_t\) and the Jordan case.} For an index \(n=-p<0\) the generating field
has a pole at \(0\), so the set written \(\varphi_{n,t}(B(0,\delta))\) there is
not an image of \(B(0,\delta)\). We use the set
\(K_t(\delta)=\widehat{\mathbb C}\setminus\phi_t(D_\delta)\) of
Appendix~\ref{app:caiInput}, which is the set
\(\widehat{\mathbb C}\setminus\omega^{\varepsilon\nu}(D)\) of
\cite[Proposition~3.1]{GQW}. The corollary is stated for all Jordan
curves separating \(0\) and \(\infty\); its proof passes from
Weil--Petersson quasicircles to Jordan curves by an approximation
\(\gamma_r=f(rS^1)\) whose continuity step is not written out. We prove
the formula in the form used here, first for analytic curves from the
Loewner energy (Lemma~\ref{lem:GQWanalytic}) and then for Jordan curves
with the continuity step (Lemma~\ref{lem:GQWext}); the values agree with
\cite[Corollary~3.9]{GQW}, which is itself not used.
\item[(C4)] \emph{\cite[Theorem~5.1]{ZhanLoop}, the constant.} The
conformal restriction covariance is displayed there with
\(e^{c\Lambda^*}\). For the Field--Lawler \(\Lambda^*\)
\cite[Theorem~1.1]{FieldLawler} we use the exponent \(\frac c2\Lambda^*\),
which is the one used in \cite[(2.3)]{CarfagniniWang},
\cite[Definition~2.1]{GQW}, \cite[Definition~1.8]{BJloop} and
\cite[Section~4.2]{FS}, and the one compatible with the constant of the
bubble decomposition in Zhan's proof \cite[(3.35)--(3.36)]{ZhanLoop}.
Appendix~\ref{app:normalization} verifies this constant
(Lemma~\ref{lem:normalization}, Step~2); equivalently, the constants of
\cite[Theorem~5.1]{ZhanLoop} as displayed are those of the loop measure
\(\frac12\mu^{\rm loop}\).
\item[(C5)] \emph{Takhtajan--Teo, Ricci form.} \cite{TTI} write
\(\mathrm{Ric}^{\rm TT}_{\rm WP}=\frac i2\sum R_{k\bar l}\,d\varepsilon_k\wedge d\bar\varepsilon_l\).
With the convention \(\mathrm{Ric}(J\cdot,\cdot)=i\sum R_{k\bar l}\,d\varepsilon_k\wedge d\bar\varepsilon_l\)
used for Kähler forms in Section~\ref{sec:KE}, the Ricci form of
\(g_{\rm WP}\) is \(2\mathrm{Ric}^{\rm TT}_{\rm WP}=-\omega_{26,0}\); the
literal value \(\mathrm{Ric}^{\rm TT}_{\rm WP}=-\frac12\omega_{26,0}\) differs
only by the factor \(\frac i2\) in the definition
(Lemma~\ref{lem:conventions}(iii)).
\end{enumerate}
\section{Finite-flow likelihoods}
\label{sec:finiteinput}\label{sec:finite}

The inputs are the integration-by-parts formula of
Theorem~\ref{thm:Wproof} for all modes \(n\ge1\), i.e.\ \eqref{eq:W}
(cf.\ \cite[Theorem~1.2]{BJweld}), and the deterministic
coefficient responses of Lemma~\ref{lem:coeffResp} and
Corollaries~\ref{cor:coeffBounded}--\ref{cor:coeffTm}, proved in
Appendix~\ref{app:coeffResp}; the properties of \(P\) used are those listed
in Section~\ref{sec:setting}, for \(0<\kappa\le4\). A determining-core argument
combines first and higher modes without complexifying sample orbits.

For the normalized welding maps in \eqref{eq:G}, define
\begin{equation}\label{eq:stress}
 T_g(z)=\frac{\cL}{12}\left\{\frac{g'''(z)}{g'(z)}
                 -\frac32\left(\frac{g''(z)}{g'(z)}\right)^2\right\}
       +\left(\frac{g'(z)}{g(z)}\right)^2-z^{-2}
       =\sum_{n\geq1}s_nz^{-n-2}.
\end{equation}
Write \(T_F\) for the same stress expression with \(g\) replaced
by \(F\). The Schwarzian chain rule and M\"obius invariance give
\(T_G(z)=z^{-4}T_g(1/z)=\sum_{n\geq1}s_nz^{n-2}\).
Each \(s_n\) is a polynomial in finitely many bounded schlicht
coefficients. In particular,
\[
 s_1=2\alpha_1=-2\beta_0,\qquad
 s_2=C_2\alpha_2+(3-C_2)\alpha_1^2=3\beta_0^2-C_2\beta_1,
 \qquad C_n=2n+\frac{\cL}{12}(n^3-n).
\]

\begin{lemma}[The law and the two composition actions]
\label{lem:CaiDictionary}
For the Borel involution \(\mathsf J(h)=h^{-1}\),
\begin{equation}\label{eq:inverseconjugacy}
 \mathsf J\circ L_a=\mathsf R_{a^{-1}}\circ \mathsf J,
 \qquad \mathsf R_a(h)=h\circ a.
\end{equation}
For the interior stress coefficient
\(\mathfrak t_n(h)=[z^{n-2}]T_{f_h}(z)\),
\begin{equation}\label{eq:stressdictionary}
 f_{\mathsf J(h)}(z)=\overline{G_h(\bar z)},\qquad
 \mathfrak t_n(\mathsf J(h))=\overline{s_n(h)}.
\end{equation}
Consequently, if \((\mathsf R_{a^{-1}})_*P\) has density
\(Z_{a^{-1}}\) relative to \(P\), then
\begin{equation}\label{eq:RNconjugacy}
 \frac{d(L_a)_*P}{dP}(x)=Z_{a^{-1}}(\mathsf J(x)).
\end{equation}
\end{lemma}
\begin{proof}
Equation~\eqref{eq:inverseconjugacy} is the identity
\((a\circ h)^{-1}=h^{-1}\circ a^{-1}\).

Let \(\iota(z)=1/\bar z\), which fixes the circle pointwise. The
pair \((\iota\circ g\circ\iota,\iota\circ f\circ\iota)\)
welds \(h^{-1}\). Common postmultiplication by
\(\bar b=\overline{g'(\infty)}\) normalizes its interior
derivative to one, proving \eqref{eq:stressdictionary}.
Changing variables against bounded Borel tests using
\eqref{eq:inverseconjugacy} and \(\mathsf J_*P=P\) proves
\eqref{eq:RNconjugacy}.

The capacity is inversion invariant: for the pair above,
\(\widetilde f'(0)=1/\overline{g'(\infty)}\) and
\(\widetilde g'(\infty)=1/\overline{f'(0)}\), so \(K\circ \mathsf J=K\).
\end{proof}

Coefficient identities are used on the full-measure set \(\mathcal W\)
of weldings of conformally removable curves. It is invariant under
quasisymmetric right composition (Lemma~\ref{lem:entQ}; cf.\
\cite[Lemma~4.1]{FS}), and under
inversion, since \(h^{-1}\) welds the reflected curve
\(\iota(\gamma)\), which is removable together with \(\gamma\). By
\eqref{eq:inverseconjugacy} it is therefore also invariant under
quasisymmetric left composition. We choose bounded Borel versions of the
coefficients elsewhere.

For smooth real \(a\) on \(\mathbb T\) and \(b\) on \(S^1\), put
\[
 \Gamma_{a,b}(h)=\frac1{2\pi}\int_0^{2\pi}
                  a(\theta)b(h(e^{i\theta}))\,d\theta.
\]
Let \(\mathcal C\) consist of
\(\phi(\Gamma_{a_1,b_1},\ldots,\Gamma_{a_m,b_m})\), where
\(\phi\in C_b^\infty(\mathbb R^m;\mathbb R)\), for all finite
lists. This bounded Borel algebra contains constants and each
\(\Gamma_{a,b}\), by a smooth cutoff equal to the identity on its range.

It determines finite signed Borel measures: the weights
\(1,\cos(n\theta),\sin(n\theta)\) and the two functions
\(b(z)=\Rea z,\Ima z\) recover all Fourier coefficients of
the continuous function \(h\). These countably many continuous
observables give an injective observation map from the Polish space
\(X\) to \(\mathbb R^{\mathbb N}\). Lusin--Souslin gives a
Borel inverse on its image, so they generate the Borel sigma field;
the bounded-algebra monotone class theorem proves the assertion.

For a smooth angular field \(v\), write \(\psi_t^v\) for its
circle flow and \(T_t h=h\circ\psi_t^v\). If
\(\widetilde\psi_t^v\) is an angular lift, a change of variables
gives
\begin{equation}\label{eq:corestability}
 \Gamma_{a,b}(T_t h)=\Gamma_{a_t,b}(h),\qquad
 a_t=(a\circ\widetilde\psi_{-t}^v)
                      (\widetilde\psi_{-t}^v)'.
\end{equation}
Thus \(\mathcal C\) is stable under every smooth right flow,
with \(D_v\Gamma_{a,b}=-\Gamma_{(av)',b}\). Smooth dependence of \(a_t\),
the bound \(|\Gamma_{a,b}(h)|\leq\|a\|_{L^1(d\theta/(2\pi))}
\|b\|_\infty\), and the cylinder chain rule give \(C^1\)
dependence in the uniform norm.

\begin{lemma}[Integration on a determining flow core]
\label{lem:measureODE}
Let \(P\) be a probability on a standard Borel space \(Y\) and
\(T_t\) a jointly Borel group of Borel bijections. Let a unital
bounded real algebra \(\mathcal A\) determine finite signed Borel
measures, be stable under \(T_t\), and have uniform-norm \(C^1\)
orbits \(t\mapsto F\circ T_t\). Write \(D\) for their derivative
at zero. If a bounded Borel \(\rho\) satisfies
\[
 \int DF\,dP=\int F\rho\,dP\qquad(F\in\mathcal A),
\]
then the actual image measures have density
\begin{equation}\label{eq:measureODEdensity}
 \frac{d(T_t)_*P}{dP}(y)
   =\exp\left\{\int_0^t\rho(T_{-s}y)\,ds\right\}.
\end{equation}
In particular the image measures and \(P\) are equivalent.
\end{lemma}
\begin{proof}
Take \(t\ge0\) (replace \(T_t,D,\rho\) by \(T_{-t},-D,-\rho\)
otherwise), put \(Q_t=(T_t)_*P\), \(c_t=\rho\circ T_{-t}\). For
\(F\in\mathcal A\), \(\frac d{dt}\int F\,dQ_t=\int D(F\circ T_t)\,dP
=\int Fc_t\,dQ_t\); integrating and using measure determination,
\(Q_t-P=\int_0^tc_sQ_s\,ds\) as measures, so \(t\mapsto Q_t\) is
Lipschitz in total variation. The bounded density
\(r_t=\exp\int_0^tc_s\,ds\) gives a second solution \(r_tP\) of this
Volterra equation, and \(\Delta_t=Q_t-r_tP\) satisfies
\(\|\Delta_t\|_{\rm var}\le\|\rho\|_\infty\int_0^t\|\Delta_s\|_{\rm var}ds\);
here \(c_s(y)\) is jointly Borel and bounded, and
\(|\partial_tr_t|\le\|\rho\|_\infty e^{\|\rho\|_\infty t}\), so both
solutions are Lipschitz in total variation and \(s\mapsto\|\Delta_s\|_{\rm var}\)
is continuous. Gronwall gives \(\Delta_t=0\), proving the density formula
and its positive bounds.
\end{proof}

For a real trigonometric field write
\[
 v(\theta)=c_0+\sum_{n=1}^N
            \{x_n\cos(n\theta)+y_n\sin(n\theta)\},\qquad
 \rho_v=\sum_{n=1}^N
       \{x_n\Ima\mathfrak t_n+y_n\Rea\mathfrak t_n\}.
\]
Choose bounded Borel coefficient versions on \(X\). By
Theorem~\ref{thm:Wproof}, which covers every mode \(n\ge1\), including
\(n=1\),
\[
 \int D_{\cos n\theta}F\,dP=\int F\,\Ima\mathfrak t_n\,dP,\qquad
 \int D_{\sin n\theta}F\,dP=\int F\,\Rea\mathfrak t_n\,dP
 \qquad(F\in\mathcal C,\ n\ge1);
\]
this is \eqref{eq:W}. Rotation has zero divergence:
\(\int D_1F\,dP=0\) by right-rotation invariance of \(P\), see
\eqref{eq:Wp-rot}. Since \(D_vF\) is real-linear in \(v\) on \(\mathcal C\) by
\eqref{eq:corestability}, \(\int D_vF\,dP=\int F\rho_v\,dP\) holds for every
real trigonometric field, including mixed modes and the first mode, so
Lemma~\ref{lem:measureODE} applies. For \(n=1\) this agrees with the
exact M\"obius likelihood \eqref{eq:mobiuscheck} of Appendix~\ref{app:mobius}
(Remark~\ref{rem:Wp-checks}(a)).

The coefficient responses of Lemma~\ref{lem:coeffResp} hold for every
real trigonometric field, mode one included, at every point of
\(\mathcal W\): derivatives of a fixed coefficient are polynomials
in finitely many schlicht coefficients, real-linear in the field. By
Corollary~\ref{cor:coeffBounded}, \(t\mapsto\rho_w(h\circ\psi^v_t)\) is
\(C^\infty\) on every orbit in \(\mathcal W\), with derivatives
\((D_v^j\rho_w)(h\circ\psi^v_t)\), and \(D_v\rho_w\) is bilinear in
\((v,w)\) and bounded, uniformly on each fixed finite Fourier unit sphere.
The equivalence just proved transfers these pointwise estimates to
\(L^\infty(P)\) versions along the flows, giving
\[
 \left\|\log\frac{d(T_t)_*P}{dP}-t\rho_v\right\|_\infty
 \leq\frac{t^2}{2}\|D_v\rho_v\|_\infty.
\]

\begin{theorem}[Finite Fourier likelihoods]\label{thm:finite}
Fix \(N<\infty\) and \(0<\kappa\leq4\). For
\(\mathbf t=(t_1,\ldots,t_N)\in\mathbb C^N\), let
\(\Phi_{\mathbf t}\) be the time-one flow of the circle field
\[
 \sum_{n=1}^N(-t_nz^{n+1}+\bar t_nz^{1-n}),
 \qquad
 b_{\mathbf t}=\sum_{n=1}^N(t_ns_n+\bar t_n\bar s_n).
\]
Then the actual left image law has density
\begin{equation}\label{eq:finiteRN}
 r_{\mathbf t}(h)
 :=\frac{dP_{\Phi_{\mathbf t}}}{dP}(h)
 =\exp\left\{\int_{-1}^0
          b_{\mathbf t}(\Phi_{u\mathbf t}\circ h)\,du\right\}.
\end{equation}
Jointly as \(\mathbf t\to0\),
\begin{equation}\label{eq:finiteDQM}
 \sqrt{r_{\mathbf t}}
 =1+\frac12b_{\mathbf t}
        +O_{L^\infty(P)}(|\mathbf t|^2).
\end{equation}
The score covariance is
\begin{equation}\label{eq:allcov}
 \E s_n=\E[s_ns_m]=0,\qquad
 \E[s_n\bar s_m]=\delta_{nm}C_n,\qquad
 C_n=2n+\frac{\cL}{12}(n^3-n),\quad m,n\geq1.
\end{equation}
Consequently the real Fisher matrix has blocks \(2C_nI_2\)
in the coordinates \((\Rea t_n,\Ima t_n)\). One may also
add a constant angular field \(c_0\), equivalently the circle field
\(ic_0z\), to the generator. The density formulas still hold
with that flow and the same \(b_{\mathbf t}\). In the joint
limit \((c_0,\mathbf t)\to0\), the DQM remainder is
\(O_{L^\infty(P)}((|c_0|+|\mathbf t|)^2)\); the constant
direction has zero score.
\end{theorem}
\begin{proof}
The angular field associated with the displayed circle field is
\[
 v_{\mathbf t}(\theta)
 =c_0-2\sum_{n=1}^N\{\Ima t_n\cos(n\theta)
                         +\Rea t_n\sin(n\theta)\}.
\]
The dictionary gives \(-\rho_{v_{\mathbf t}}\circ \mathsf J=b_{\mathbf t}\).
Conjugating Lemma~\ref{lem:measureODE} by \eqref{eq:RNconjugacy}
turns the inverse right flow into the left flow with opposite time,
proving \eqref{eq:finiteRN} in the image variable. For \(c_0=0\),
the bilinear coefficient-response bounds of Corollary~\ref{cor:coeffBounded} give
\[
 \|b_{\mathbf t}\|_\infty\leq M_N|\mathbf t|,\qquad
 \left\|\log r_{\mathbf t}-b_{\mathbf t}\right\|_\infty
       \leq L_N|\mathbf t|^2.
\]
The second bound becomes \(L_N(|c_0|+|\mathbf t|)^2\) when
a constant field is included. Exponential Taylor expansion proves
\eqref{eq:finiteDQM}, without commuting the individual mode flows.

For \(m,n\geq2\), \eqref{eq:allcov} follows from \eqref{eq:W} and
\eqref{eq:stressdictionary}, \(\mathsf J_*P=P\) and
Corollary~\ref{cor:coeffTm}, as derived in
Appendix~\ref{app:covariance}. For \(m=n=1\) the same derivation applies, Theorem~\ref{thm:Wproof}
covering \(n=1\), and gives \(\E|s_1|^2=C_1=2\), as does
Appendix~\ref{app:mobius}. The identity
\(G_{R_\alpha\circ h}(z)=e^{-i\alpha}G_h(e^{i\alpha}z)\)
gives \(s_n\) weight \(n\). Rotation invariance therefore
centers the scores and removes unconjugated moments and mixed
Hermitian moments with mode one. This proves \eqref{eq:allcov}.
The real score \((2\Rea s_n,-2\Ima s_n)\) has covariance
\(2C_nI_2\).
\end{proof}

\begin{corollary}[Second Fourier mode]\label{cor:secondMode}
Let \(\Phi_t\) be the second-mode flow in \eqref{eq:flow},
put \(P_t=(L_{\Phi_t})_*P\), and let \(S\) be defined by
\eqref{eq:score}. The actual image density is
\begin{equation}\label{eq:RN}
 \frac{dP_t}{dP}(h)
 =\exp\left\{\int_{-1}^{0}
       \bigl[tS(\Phi_{ut}\circ h)
            +\bar t\overline{S(\Phi_{ut}\circ h)}\bigr]\,du\right\}.
\end{equation}
Jointly as \(t\to0\),
\begin{equation}\label{eq:DQM}
 \sqrt{\frac{dP_t}{dP}}
   =1+\frac12(tS+\bar t\bar S)+O_{L^\infty(P)}(|t|^2),
\end{equation}
and
\begin{equation}\label{eq:variance}
 \E S=\E S^2=0,\qquad \E|S|^2=C.
\end{equation}
The Fisher matrix for \((\Rea t,\Ima t)\) is \(2CI_2\),
and \(H^2(P_t,P)=C|t|^2/2+O(|t|^3)\).
\end{corollary}
\begin{proof}
Specialize Theorem~\ref{thm:finite} to \(t_1=0,t_2=t\)
and square its half-density expansion in \(L^2(P)\).
\end{proof}
\section{The integration-by-parts formula}\label{app:Wproof}

We prove \eqref{eq:W} from published results on the SLE loop measure,
the definition of \(P\), Lemma~\ref{lem:GQWext} and classical complex
analysis, together with, at \(\kappa=4\), the removability of SLE\(_4\) loops
\cite{KMS} (an unrefereed preprint), which enters through property (b) of
Section~\ref{sec:setting}. Neither \cite[Lemma~5.2]{BJweld}, \cite[(5.3)]{BJweld} nor frame measures are used. The
mode \(n=1\) is included.

\begin{theorem}\label{thm:Wproof}
Let \(0<\kappa\le4\). For every \(F\in\mathcal C\) and every \(n\ge1\),
\[
 \int D_{\cos n\theta}F\,dP=\int F\,\Ima\mathfrak t_n\,dP,\qquad
 \int D_{\sin n\theta}F\,dP=\int F\,\Rea\mathfrak t_n\,dP .
\]
In particular \eqref{eq:W} holds.
\end{theorem}

The inputs of the proof are the following. Besides classical complex analysis
(\cite{AhlforsLectures,LehtoVirtanen,Duren,PommerenkeBB}, de Branges'
bound \(|a_m|\le m+1\) \cite{deBranges}), the proof uses:
\begin{enumerate}
\item[(P1)] the M\"obius invariance of Zhan's unrooted SLE\(_\kappa\)
loop measure \(\mu^0\) \cite[Theorem~4.2(ii)]{ZhanLoop};
\item[(P2)] the conformal restriction covariance of \(\mu^0\) with exponent
\(\frac c2\Lambda^*\), \(c=c(\kappa)=\cm\), and its
consequence \cite[Lemma~2.2]{GQW}, in the normalization
\cite[Definition~2.1, (1)]{GQW}; for the constant see
Section~\ref{sec:conventionsLit}, item (C4), and Lemma~\ref{lem:normalization};
\item[(P3)] Lemma~\ref{lem:GQWext}, with Lemma~\ref{lem:GQWanalytic} for
analytic curves (from \cite[Lemma~2.1, Corollary~2.6, Theorems~2.2--2.3]{CarfagniniWang},
\cite[Theorem~4.1]{WangAIF}, \cite[Theorem~1.4]{WangInv}, \cite[Lemma~3.7]{GQW} and
\cite{FieldLawler,LawlerWerner});
\item[(P4)] the properties (a)--(b) of \(P\) listed in Section~\ref{sec:setting}
(see also \cite[(1.3)]{BJweld}, \cite[Remark~2.14]{GQW}).
\end{enumerate}
Since \(\mu^0\) and the measure \(\nu_\kappa\) of Section~\ref{sec:setting}
differ by a constant factor on loops separating \(0\) and \(\infty\), (P1)
and (P2) hold for \(\nu_\kappa\), which we write \(\nu\).

We begin with notation for loops and their weldings. Let \(\mathcal J\) be the set of Jordan
curves \(\gamma\subset\mathbb C\) with \(0\in\operatorname{int}\gamma\).
For \(\gamma\in\mathcal J\) let \(f_\gamma:\D\to\operatorname{int}\gamma\) and
\(g_\gamma:\Ds\to\operatorname{ext}\gamma\) be the conformal maps with
\(f_\gamma(0)=0\), \(f_\gamma'(0)>0\), \(g_\gamma(\infty)=\infty\),
\(g_\gamma'(\infty)>0\), and put
\[
 c(\gamma)=\log f_\gamma'(0),\qquad \eta(\gamma)=e^{-c(\gamma)}\gamma,\qquad
 W(\gamma,\alpha)=R_{-\alpha}\circ g_\gamma^{-1}\circ f_\gamma|_{S^1}
 \quad(\alpha\in\mathbb R/2\pi\mathbb Z).
\]
Then \(W(e^c\gamma,\alpha)=W(\gamma,\alpha)\), and by Section~\ref{sec:setting}
\(P\) is the law of \(W(\eta,\alpha)\) under \(\mu_\kappa\otimes d\alpha/2\pi\).
We write \(\Xi=\mathcal J\times\mathbb R/2\pi\mathbb Z\) with the
measure \(\nu\otimes d\alpha/2\pi\). For \(a<b\) put
\(S_{a,b}=\{\gamma\in\mathcal J:a\le c(\gamma)\le b\}\); then
\(\nu(S_{a,b})=\pi^{-1/2}(b-a)\), and by Koebe's one-quarter theorem
\(B(0,e^{c(\gamma)}/4)\subset\operatorname{int}\gamma\).

For the shape \(\eta\) we write \(f=f_\eta=z(1+\sum_{m\ge1}a_mz^m)\); as in
Appendix~\ref{app:caiInput}, \(a_m\) are interior coefficients (the
exterior coefficients of Section~\ref{sec:setting} are written
\(\alpha_m,\beta_m\)), \(c\)
denotes the log conformal radius and \(\cm\) the central charge, \(\nu\)
the loop measure and \(\varpi\) Beltrami coefficients, and \(P_p\) is never
the law \(P\). A
\emph{coefficient polynomial} is an element of \(\mathbb C[a_m,\bar a_m:m\ge1]\),
regarded as a function of \(\gamma\in\mathcal J\) through \(\eta(\gamma)\).
By removability, \(a_m(W(\eta,\alpha))=a_m(\eta)\) for \(\mu_\kappa\)-almost
every \(\eta\), where on the left \(a_m\) is the Borel coefficient version
on \(X\); so a coefficient polynomial is also a bounded Borel function on
\(X\), \(P\)-almost surely equal to the loop-side one. All coefficient
polynomials are bounded, by de Branges' theorem.

Next we introduce the Laurent data of the interior map, which appear in the
change of basis below. Here \([z^k]\) denotes a Laurent coefficient at
\(0\), and \(\Pi_{\ge m}\) the projection
\(\sum_kc_kz^k\mapsto\sum_{k\ge m}c_kz^k\). For \(p\ge1\) and \(0\le j\le p\) put
\begin{equation}\label{eq:Wp-E}
 E_{p,j}=[z^{1-j}]\frac{f^{1-p}}{f'},\qquad e_p=E_{p,0},\qquad
 M_{n,p}=[z^{n-2}]\bigl(f^{p-2}f'^2\bigr)\quad(n\ge1).
\end{equation}
These are holomorphic polynomials in \(a_1,a_2,\ldots\), and
\(E_{p,p}=M_{p,p}=1\), \(E_{p,j}=0\) for \(j>p\), \(M_{n,p}=0\) for \(p>n\).
Under \(a_m\mapsto e^{im\beta}a_m\) they have weights
\(E_{p,j}\mapsto e^{i(p-j)\beta}E_{p,j}\), \(M_{n,p}\mapsto e^{i(n-p)\beta}M_{n,p}\).
With \(\psi=f^{-1}\), \(P_p=[w^{p-2}]\mathcal S\psi\) for \(p\ge2\) as in
Appendix~\ref{app:caiInput}, and we set \(P_1=0\).

Finally we complexify the right-flow derivatives. For a real \(C^1\) field \(\dot\theta\) on
\(\mathbb T\) and \(F\in\mathcal C\), \(D_{\dot\theta}F\) is given by
\eqref{eq:corestability} and the chain rule: for \(F=\phi(\Gamma_{a_1,b_1},\dots)\),
\(D_{\dot\theta}F=-\sum_r\partial_r\phi\,\Gamma_{(a_r\dot\theta)',b_r}\). It is
real-linear in \(\dot\theta\), and we extend it complex-linearly. Put
\(\mathfrak d_j=D_{e^{ij\theta}}\), \(j\in\mathbb Z\). Then
\(\mathfrak d_0=D_1\), \(\mathfrak d_{\pm n}=D_{\cos n\theta}\pm iD_{\sin n\theta}\),
\(D_{\dot\theta}=\sum_j\dot\theta_j\mathfrak d_j\) if
\(\dot\theta=\sum_j\dot\theta_je^{ij\theta}\), and
\(\mathfrak d_{-n}F=\overline{\mathfrak d_nF}\) for real \(F\). Hence, for real
\(F\),
\begin{equation}\label{eq:Wp-complex}
 \int\mathfrak d_{-n}F\,dP=-i\int F\,\mathfrak t_n\,dP
\end{equation}
is equivalent, by taking real and imaginary parts, to the pair
\(\int D_{\cos n\theta}F\,dP=\int F\Ima\mathfrak t_n\,dP\),
\(\int D_{\sin n\theta}F\,dP=\int F\Rea\mathfrak t_n\,dP\).
We prove \eqref{eq:Wp-complex}.

\begin{lemma}[Cylinder derivatives along \(C^1\) families]\label{lem:Wp-cyl}
Let \(h\in X\) and let \(k_s\), \(|s|<s_0\), be orientation-preserving
\(C^1\) diffeomorphisms of \(S^1\) with \(k_s(e^{i\theta})=e^{i\vartheta_s(\theta)}\),
\(\vartheta_s(\theta)=\theta+s\dot\theta(\theta)+o(s)\) in \(C^1(\mathbb T)\). Then
\(F(h\circ k_s)=F(h)+sD_{\dot\theta}F(h)+o(s)\) for \(F\in\mathcal C\).
\end{lemma}
\begin{proof}
As in \eqref{eq:corestability}, \(\Gamma_{a,b}(h\circ k_s)=\Gamma_{a_s,b}(h)\) with
\(a_s=(a\circ\vartheta_s^{-1})(\vartheta_s^{-1})'\), and
\(\vartheta_s^{-1}=\mathrm{id}-s\dot\theta+o(s)\) in \(C^1\). Hence
\(a_s=a-s(a\dot\theta)'+o(s)\) uniformly, and
\(|\Gamma_{a',b}|\le\|a'\|_\infty\|b\|_\infty\). The chain rule finishes the proof.
\end{proof}

The rotation \(h\mapsto h\circ R_\beta\) replaces \(f_h\) by
\(e^{-i\beta}f_h(e^{i\beta}z)\), i.e.\ \(a_m\mapsto e^{im\beta}a_m\). Hence a
holomorphic coefficient polynomial \(\Phi\) of weight \(k\) satisfies
\(\mathfrak d_0\Phi=ik\Phi\). Right-rotation invariance of \(P\)
(Section~\ref{sec:setting}; it follows from (P1)) gives, for every bounded
Borel \(H\) with \(\beta\mapsto H(h\circ R_\beta)\) differentiable at \(0\) with
derivative \(D_1H\) bounded uniformly in \(h\) (and along each orbit, so
that by the mean value theorem the difference quotients are uniformly
bounded and dominated convergence applies),
\begin{equation}\label{eq:Wp-rot}
 \int D_1H\,dP=0 .
\end{equation}

The first step of the proof is a family of Ward identities for the loop
measure, one for each holomorphic field \(-\lambda\xi^{1-p}\).
For \(p\ge1\) and \(\lambda\in\mathbb C\) let \(\varphi^{p,\lambda}_s\) be the flow of
the holomorphic field \(u_{p,\lambda}(\xi)=-\lambda\xi^{1-p}\), \(\dot\xi=u(\xi)\).
For \(p=1\) it is the translation \(\xi\mapsto\xi-s\lambda\). For \(p\ge2\) and
\(\lambda\in\{1,i\}\) it is the flow \(\phi_s\) of Appendix~\ref{app:caiInput}.
For \(r>0\), \(\varphi^{p,r\lambda}_s=\varphi^{p,\lambda}_{rs}\), and for
\(c\in\mathbb R\),
\begin{equation}\label{eq:Wp-scale}
 \varphi^{p,\lambda}_s(e^c\xi)=e^c\varphi^{p,\lambda e^{-pc}}_s(\xi).
\end{equation}
For a function \(H\) on \(\Xi\) write
\(X^\lambda_pH(\gamma,\alpha)=\frac{d}{ds}\big|_{s=0}H(\varphi^{p,\lambda}_s\gamma,\alpha)\)
when the derivative exists.

\begin{proposition}[Ward identity]\label{prop:Wp-loopWard}
Let \(p\ge1\), \(\lambda\in\{1,i\}\) and \(a<b\). Let \(H\) be a bounded Borel
function on \(\Xi\) vanishing unless \(\gamma\in S_{a,b}\). Assume that
\(X^\lambda_pH(\gamma,\alpha)\) exists for all \(\gamma\in S_{a-1,b+1}\) and all
\(\alpha\), and is bounded there. Then
\begin{equation}\label{eq:Wp-loopWard}
 \iint X^\lambda_pH\,d\nu\,\frac{d\alpha}{2\pi}
 =-\frac{\cm}{6}\iint H\,\Rea\bigl(\lambda P_p(\gamma)\bigr)\,d\nu\,\frac{d\alpha}{2\pi}.
\end{equation}
\end{proposition}

\begin{proof}
Put \(\delta=e^{a-1}/8\), so that every \(\gamma\in S_{a-1,b+1}\) lies in
\(\mathcal J_\delta\) (Appendix~\ref{app:caiInput}). Along a trajectory,
\(|\frac d{ds}\varphi_s(\xi)|\le|\xi|^{1-p}\) and, by Proposition~\ref{prop:Wp-var}
below, \(\frac d{ds}c(\varphi_s\gamma)=-\Rea(\lambda e^{-pc}e_p(\eta))\) with
\(|e_p|\) bounded by de Branges' theorem. Hence there is \(\varepsilon'>0\) such
that, for \(|t|\le\varepsilon'\) and \(\gamma\in S_{a,b}\), the loops
\(\varphi_{\pm s}\gamma\), \(|s|\le|t|\), are defined and lie in
\(S_{a-1,b+1}\); for \(p\ge2\) we also take \(\varepsilon'\le\varepsilon\) of
Lemma~\ref{lem:GQWext}. Let \(|t|\le\varepsilon'\). If \(H(\gamma,\alpha)\ne0\) or
\(H(\varphi_t\gamma,\alpha)\ne0\), then \(\gamma\) or \(\varphi_t\gamma\) lies in \(S_{a,b}\), so
\(\varphi_s\gamma\in S_{a-1,b+1}\) for all \(s\) between \(0\) and \(t\). By the flow rule
\(\frac d{ds}H(\varphi_s\gamma,\alpha)=X^\lambda_pH(\varphi_s\gamma,\alpha)\) there, and the
mean value theorem gives, for all \(\gamma\) and \(\alpha\),
\begin{equation}\label{eq:Wp-Lip}
 |H(\varphi_t\gamma,\alpha)-H(\gamma,\alpha)|\le C'|t|\,\mathbf 1_{S_{a-1,b+1}}(\gamma)
 \qquad(|t|\le\varepsilon'),
\end{equation}
with \(C'=\sup|X^\lambda_pH|\) over \(S_{a-1,b+1}\).

\emph{The case \(p=1\).} Let \(|t|<\min(\varepsilon',e^{a}/8)\), write \(\mathcal J(z)\) for the
Jordan curves surrounding \(z\), and take \(\nu\) to be \(\mu^0\) restricted to
\(\mathcal J(0)=\mathcal J\) (constant factors are irrelevant). If
\(H(\gamma-t\lambda,\alpha)\ne0\), then \(\gamma-t\lambda\in S_{a,b}\) surrounds
\(B(0,e^a/4)\ni-t\lambda\), so \(\gamma\in\mathcal J(0)\cap\mathcal J(t\lambda)\). Hence
\[
 \int_{\mathcal J(0)}H(\gamma-t\lambda,\alpha)\,d\mu^0(\gamma)
 =\int_{\mathcal J(t\lambda)}H(\gamma-t\lambda,\alpha)\,d\mu^0(\gamma)
 =\int_{\mathcal J(0)}H(\gamma',\alpha)\,d\mu^0(\gamma'),
\]
because the translation \(\gamma\mapsto\gamma-t\lambda\) maps \(\mathcal J(t\lambda)\) onto
\(\mathcal J(0)\) and preserves \(\mu^0\) (P1). Differentiate at \(t=0\), using
\eqref{eq:Wp-Lip}, \(\nu(S_{a-1,b+1})<\infty\) and dominated convergence. This is
\eqref{eq:Wp-loopWard} with \(P_1=0\).

\emph{The case \(p\ge2\).} Fix \(|t|\le\varepsilon'\) and \(\alpha\). The map
\(\varphi_t\) is conformal from \(D=D_\delta\) onto \(D'=\varphi_t(D_\delta)\), with
complement \(K_t=K_t(\delta)\) (Appendix~\ref{app:caiInput}). By
Lemma~\ref{lem:normalization}(ii), for every bounded Borel \(\Psi\) on loops in
\(D'\) vanishing outside a set of finite measure,
\[
 \int_{\gamma\subset D}\Psi(\varphi_t\gamma)\,d\nu(\gamma)
 =\int_{\gamma'\subset D'}\Psi(\gamma')
 \exp\Bigl\{\frac{\cm}2\bigl[\Lambda^*(\gamma',K_t)-\Lambda^*(\varphi_t^{-1}\gamma',K_0)\bigr]\Bigr\}
 d\nu(\gamma').
\]
(This is an identity of measures on all loops; a loop \(\gamma\subset D\)
separates \(0\) from \(\infty\) iff it surrounds \(K_0\), iff \(\varphi_t\gamma\) surrounds
the connected set \(K_t\ni0\), so it restricts to \(\mathcal J\).) Apply it to
\(\Psi(\gamma')=H(\gamma',\alpha)\exp\{-\frac{\cm}2[\cdots]\}\). Every loop in
\(S_{a,b}\) lies in \(D'\), since its \(\varphi_{-t}\)-image lies in
\(S_{a-1,b+1}\subset\mathcal J_\delta\). With \(A_t\) as in
Appendix~\ref{app:caiInput} we obtain
\[
 \int H(\gamma',\alpha)\,d\nu(\gamma')
 =\int_{\mathcal J_\delta}H(\varphi_t\gamma,\alpha)\,e^{-\frac{\cm}2A_t(\gamma)}\,d\nu(\gamma).
\]
By Lemma~\ref{lem:GQWext}, \(|A_t(\gamma)|\le\frac{|t|}3\sup_{\mathcal J_{\delta/2}}|P_p|\) on
\(\mathcal J_\delta\) and \(\frac d{dt}A_t(\gamma)|_{t=0}=-\frac13\Rea(\lambda P_p(\gamma))\).
Together with \eqref{eq:Wp-Lip}, the difference quotients of the integrand on
the right are bounded by a constant and vanish off \(S_{a-1,b+1}\); they
converge to \(X^\lambda_pH+\frac{\cm}6H\Rea(\lambda P_p)\). Dominated convergence and
the \(t\)-independence of the left side give
\(\int(X^\lambda_pH+\frac{\cm}6H\Rea(\lambda P_p))\,d\nu=0\). Integrate in \(\alpha\).
\end{proof}

The second step computes how the welding and the shape of a curve move
along the Ward flows.
The next two lemmas are deterministic and hold for every Jordan curve.

\begin{lemma}[First variation of the Beltrami equation]\label{lem:Wp-AB}
Let \(R>0\), \(0\le k<1\), and let \(\varpi_s\), \(|s|\le s_0\), be measurable on
\(\mathbb C\), supported in \(B(0,R)\), with \(\|\varpi_s\|_\infty\le k\) and
\(\|\varpi_s-s\dot\varpi\|_\infty\le M_0s^2\). Let \(F_s\) be the principal solution of
\(\partial_{\bar z}F=\varpi_s\partial_zF\) (\(F_s(z)=z+O(1/z)\) at \(\infty\)) and
\(\mathsf C[g](z)=-\frac1\pi\int g(\zeta)(\zeta-z)^{-1}dA(\zeta)\). Then
\(\sup_{z\in\mathbb C}|F_s(z)-z-s\mathsf C[\dot\varpi](z)|=O(s^2)\).
\end{lemma}
\begin{proof}
By the construction of the principal solution
\cite[Chapter~V]{AhlforsLectures}, see also \cite{AhlforsBers,LehtoVirtanen}, \(F_s=z+\mathsf C[\omega_s]\), where
\(\omega_s\in L^p\) solves \(\omega_s=\varpi_s+\varpi_sT\omega_s\), \(T\) is the Beurling
transform, and \(p>2\) is chosen with \(k\|T\|_{L^p\to L^p}<1\) (possible since
\(\|T\|_{L^p\to L^p}\to1\) as \(p\to2\)). Then
\(\|\omega_s\|_p\le C\|\varpi_s\|_p\) and \(\|\omega_s-\varpi_s\|_p\le C\|\varpi_s\|_\infty\|\varpi_s\|_p=O(s^2)\),
because the supports lie in \(B(0,R)\). By H\"older's inequality,
\(|\mathsf C[g](z)|\le C_{p,R}\|g\|_p\) for \(g\) supported in \(B(0,R)\), uniformly in
\(z\in\mathbb C\). Hence \(F_s-z-\mathsf C[\varpi_s]=\mathsf C[\omega_s-\varpi_s]=O(s^2)\) and
\(\mathsf C[\varpi_s]-s\mathsf C[\dot\varpi]=\mathsf C[\varpi_s-s\dot\varpi]=O(s^2)\), uniformly.
\end{proof}

\begin{proposition}[Welding and shape variation]\label{prop:Wp-var}
Let \(\gamma\in\mathcal J\), \(f=f_\gamma\), \(p\ge1\), \(\lambda\in\mathbb C\),
\(\varphi_s=\varphi^{p,\lambda}_s\) and \(U=(u_{p,\lambda}\circ f)/f'=-\lambda f^{1-p}/f'\),
a meromorphic function on \(\D\) with its only pole at \(0\). For small \(|s|\)
the loop \(\varphi_s\gamma\) lies in \(\mathcal J\), and:
\begin{enumerate}
\item \(W(\varphi_s\gamma,\alpha)=W(\gamma,\alpha)\circ k_s\) for all \(\alpha\), where \(k_s\)
is a real-analytic circle diffeomorphism and
\(k_s(e^{i\theta})=e^{i(\theta+s\dot\theta(\theta)+o(s))}\) in \(C^1(\mathbb T)\), with
the trigonometric polynomial
\begin{equation}\label{eq:Wp-velocity}
 \dot\theta_{-j}=iU_{1-j}\ (j\ge1),\qquad \dot\theta_j=\overline{\dot\theta_{-j}},
 \qquad \dot\theta_0=-\Ima U_1,\qquad U_k=[z^k]U ;
\end{equation}
\item \(f_{\varphi_s\gamma}=f+s\,\delta f+o(s)\) locally uniformly on \(\D\), where
\(\delta f=f'Y\) and
\begin{equation}\label{eq:Wp-Y}
 Y=(\Rea U_1)\,z+\sum_{m\ge2}\bigl(U_m+\overline{U_{2-m}}\bigr)z^m ;
\end{equation}
in particular \(\frac d{ds}c(\varphi_s\gamma)|_{s=0}=\Rea U_1\).
\end{enumerate}
\end{proposition}

\begin{proof}
Choose \(\rho>0\) with \(\overline{B(0,4\rho)}\subset\operatorname{int}\gamma\), and a smooth
radial \(\chi\) with \(\chi=0\) on \(B(0,\rho)\) and \(\chi=1\) off \(B(0,2\rho)\). Since
\(u=u_{p,\lambda}\) is holomorphic on \(\widehat{\mathbb C}\setminus\{0\}\) and vanishes at
\(\infty\) (to order \(p+1\ge2\)), \(\chi u\) is a smooth vector field on
\(\widehat{\mathbb C}\). Let \(Q_s\) be its flow. For small \(|s|\), \(Q_s=\mathrm{id}\) on
\(B(0,\rho/2)\), \(Q_s=\varphi_s\) on \(\{|\xi|\ge3\rho\}\cup\{\infty\}\), which contains
\(\gamma\cup\operatorname{ext}\gamma\), and the Beltrami coefficient of \(Q_s\) is smooth
in \((s,\xi)\), supported in \(\{\rho/2\le|\xi|\le3\rho\}\), and equals
\(s\,\partial_{\bar\xi}(\chi u)+O(s^2)\) uniformly. Hence \(\varphi_s\gamma=Q_s\gamma\in\mathcal J\),
\(g_{\varphi_s\gamma}=\varphi_s\circ g_\gamma\) (the normalization at \(\infty\) is
preserved because \(\varphi_s(\xi)=\xi+O(1)\) for \(p=1\) and
\(\varphi_s(\xi)=\xi+O(\xi^{1-p})\) for \(p\ge2\)), and \(Q_s\circ f:\D\to\operatorname{int}\varphi_s\gamma\)
is a quasiconformal homeomorphism fixing \(0\), with Beltrami coefficient
\(\varpi_s=(\mu_{Q_s}\circ f)\,\overline{f'}/f'\), smooth, supported in a compact
\(L\subset\D\setminus\{0\}\), and \(\varpi_s=s\dot\varpi+O(s^2)\) uniformly.

Let \(\varpi^\sharp_s=\varpi_s\) on \(\D\) and
\(\varpi^\sharp_s(z)=\overline{\varpi_s(1/\bar z)}\,z^2/\bar z^2\) on \(\Ds\); it is supported in
the compact \(L\cup\iota(L)\subset\mathbb C\setminus\{0\}\), \(\iota(z)=1/\bar z\). Let \(E_s\) be
the solution with coefficient \(\varpi^\sharp_s\) fixing \(0,1,\infty\). Then
\(\iota\circ E_s\circ\iota\) has the same coefficient and the same three fixed points, so
equals \(E_s\) by uniqueness of normalized solutions \cite[Chapter~V]{AhlforsLectures}; thus \(E_s\) preserves
\(S^1\) and \(\D\). By Lemma~\ref{lem:Wp-AB} (applied to the principal solution
and the affine normalization \(E_s=(F_s-F_s(0))/(F_s(1)-F_s(0))\)),
\(E_s=\mathrm{id}+s\dot E+O(s^2)\) uniformly on compact sets. Let
\(\theta_s=\arg E_s'(0)\) (\(E_s\) is conformal near \(0\)) and
\(\widetilde E_s=e^{-i\theta_s}E_s\). The maps \(\widetilde E_s\) are conformal on a closed
annulus \(A=\{r_*\le|z|\le1/r_*\}\) disjoint from \(L\cup\iota(L)\), with
\(f(\{r_*\le|z|<1\})\subset\{|\xi|>3\rho\}\); Cauchy estimates give
\(\widetilde E_s^{-1}=\mathrm{id}+sw+O(s^2)\) in \(C^1\) on a slightly smaller annulus, where
\(w=-\frac d{ds}\widetilde E_s|_{s=0}\) is holomorphic on \(A^\circ\) and, since
\(\widetilde E_s\) commutes with \(\iota\), satisfies
\begin{equation}\label{eq:Wp-reflect}
 w(z)=-z^2\,\overline{w(1/\bar z)} .
\end{equation}
In particular \(w(z)/(iz)\) is real on \(S^1\).

The map \(f_s:=Q_s\circ f\circ\widetilde E_s^{-1}\) is conformal on \(\D\), because
\(Q_s\circ f\) and \(\widetilde E_s\) have the same Beltrami coefficient there
\cite[Chapter~I]{AhlforsLectures}, \cite{LehtoVirtanen}. It maps \(\D\) onto
\(\operatorname{int}\varphi_s\gamma\), fixes \(0\), and
\(f_s'(0)=f'(0)/\widetilde E_s'(0)>0\). Hence \(f_s=f_{\varphi_s\gamma}\). On \(S^1\),
\(f_s=\varphi_s\circ f\circ\widetilde E_s^{-1}\), so
\(W(\varphi_s\gamma,\alpha)=R_{-\alpha}\circ g^{-1}\circ\varphi_s^{-1}\circ\varphi_s\circ f\circ\widetilde E_s^{-1}
=W(\gamma,\alpha)\circ k_s\) with \(k_s=\widetilde E_s^{-1}|_{S^1}\). Writing
\(w(e^{i\theta})=ie^{i\theta}\dot\theta(\theta)\) gives the \(C^1\) expansion in (i).

From \(f_s=Q_s\circ f\circ\widetilde E_s^{-1}\) and the expansions of \(Q_s\) and
\(\widetilde E_s^{-1}\), the difference quotients \((f_s-f)/s\) converge locally
uniformly on \(\D\); their limit \(\delta f\) is holomorphic, \(\delta f(0)=0\), and
\(\delta f'(0)\in\mathbb R\). Put \(Y=\delta f/f'\), holomorphic on \(\D\) with
\(Y(0)=0\), \(Y'(0)=\delta f'(0)/f'(0)\in\mathbb R\). On \(\{r_*<|z|<1\}\) we have
\(f_s\circ\widetilde E_s=\varphi_s\circ f\); differentiating, \(\delta f=u\circ f+f'w\), i.e.
\(w=Y-U\) there. So \(w\) extends meromorphically to \(\D\) with its only pole at
\(0\), and by \eqref{eq:Wp-reflect} to \(\Ds\) with its only pole at \(\infty\):
\(w=\sum_kw_kz^k\) is a Laurent polynomial. Comparing coefficients in
\(w=Y-U\) on \(0<|z|<1\): \(w_k=-U_k\) for \(k\le0\); by \eqref{eq:Wp-reflect},
\(w_m=-\overline{w_{2-m}}=\overline{U_{2-m}}\) for \(m\ge2\) and \(w_1\in i\mathbb R\); since
\(Y_1=U_1+w_1\) is real, \(w_1=-i\Ima U_1\) and \(Y_1=\Rea U_1\); and
\(Y_m=U_m+w_m\) for \(m\ge2\). This is \eqref{eq:Wp-Y}, and
\(\dot\theta=w/(iz)\) gives \eqref{eq:Wp-velocity}. Finally
\(\frac d{ds}\log f_s'(0)=Y'(0)=\Rea U_1\).
\end{proof}

The third step transfers the Ward identity from the loop measure to the
welding law.
Fix \(p\ge1\) and \(\lambda\in\mathbb C\). At a shape \(\eta\), with \(f=f_\eta\), we have
\(U=\lambda V_p\), \(V_p=-f^{1-p}/f'\), and \([z^{1-j}]V_p=-E_{p,j}\). By
Proposition~\ref{prop:Wp-var}, the welding moves by the trigonometric polynomial
\(w^\lambda_p=w^\lambda_p(\eta)\) with
\begin{equation}\label{eq:Wp-wfield}
 \dot\theta_{-j}=-i\lambda E_{p,j}\ (1\le j\le p),\qquad
 \dot\theta_0=\Ima(\lambda e_p),
\end{equation}
the scale by \(\Rea U_1=-\Rea(\lambda e_p)\), and the normalized shape map
\(f_{\varphi_s\eta}/f'_{\varphi_s\eta}(0)\) by
\(\delta^\lambda_pf=f'Y-(\Rea U_1)f\). For \(F\in\mathcal C\) and a coefficient
polynomial \(\Phi\) define
\[
 \mathfrak D^\lambda_p(F\Phi)=(D_{w^\lambda_p}F)\,\Phi+F\,d\Phi[\delta^\lambda_pf],
\]
regarded as a function of \(h=W(\eta,\alpha)\) through \(f_h=f_\eta\) on
\(\mathcal W\). It is the derivative of \(s\mapsto F(W(\varphi_s\eta,\alpha))\Phi(\eta(\varphi_s\eta))\) at
\(s=0\) (Lemma~\ref{lem:Wp-cyl} and Proposition~\ref{prop:Wp-var}). All three
ingredients are real-linear in \(\lambda\), and their coefficients are
polynomials in \(a_m,\bar a_m\), hence bounded.

\begin{proposition}[Ward identity for the welding law]\label{prop:Wp-weldWard}
For \(p\ge1\), \(\lambda\in\mathbb C\), \(F\in\mathcal C\) and every coefficient
polynomial \(\Phi\),
\begin{equation}\label{eq:Wp-weldWard}
 \int\mathfrak D^\lambda_p(F\Phi)\,dP=\int F\Phi\,\Rea(\lambda Q_p)\,dP,\qquad
 Q_p=p\,e_p-\frac{\cm}6P_p .
\end{equation}
\end{proposition}
\begin{proof}
By real-linearity it suffices to take \(\lambda\in\{1,i\}\) and \(\Phi\) real. Fix
\(\chi\in C_c^\infty((0,1))\), \(\chi\ge0\), \(\chi\ne0\), and put
\(H(\gamma,\alpha)=F(W(\gamma,\alpha))\,\Phi(\eta(\gamma))\,\chi(c(\gamma))\). By
\eqref{eq:Wp-scale}, the scale invariance of \(W\) and of \(\eta\), and
Proposition~\ref{prop:Wp-var} applied at \(\eta\) with \(\lambda e^{-pc}\), for
\(\gamma=e^c\eta\in\mathcal J\),
\[
 X^\lambda_pH(\gamma,\alpha)=e^{-pc}\bigl\{\chi(c)\,\mathfrak D^\lambda_p(F\Phi)(\eta,\alpha)
 -\chi'(c)\Rea(\lambda e_p(\eta))\,(F\Phi)(\eta,\alpha)\bigr\},
\]
which is bounded on \(S_{-1,2}\times\mathbb R/2\pi\mathbb Z\). Since
\(f_{e^c\eta}=e^cf_\eta\), the inverse maps satisfy
\(\psi_{e^c\eta}(w)=\psi_\eta(e^{-c}w)\), so by the scaling of the Schwarzian
\(\mathcal S\psi_{e^c\eta}(w)=e^{-2c}\mathcal S\psi_\eta(e^{-c}w)\) and
\(P_p(e^c\eta)=e^{-pc}P_p(\eta)\). Proposition~\ref{prop:Wp-loopWard} with \(a=0\),
\(b=1\), the disintegration \(\nu=\pi^{-1/2}dc\otimes\mu_\kappa\) and the definition of
\(P\) give
\begin{multline*}
 \int e^{-pc}\chi(c)\,dc\int\mathfrak D^\lambda_p(F\Phi)\,dP
 -\int e^{-pc}\chi'(c)\,dc\int F\Phi\,\Rea(\lambda e_p)\,dP\\
 =-\frac{\cm}6\int e^{-pc}\chi(c)\,dc\int F\Phi\,\Rea(\lambda P_p)\,dP.
\end{multline*}
Since \(\int e^{-pc}\chi'(c)\,dc=p\int e^{-pc}\chi(c)\,dc>0\), dividing gives
\eqref{eq:Wp-weldWard}.
\end{proof}

It is convenient to complexify these identities. Put
\(\mathcal L_{-p}=\frac12(\mathfrak D^1_p-i\mathfrak D^i_p)\). Since
\(\Rea(Q)-i\Rea(iQ)=Q\), \eqref{eq:Wp-weldWard} gives, for complex coefficient
polynomials \(\Phi\),
\begin{equation}\label{eq:Wp-complexWard}
 \int\mathcal L_{-p}(F\Phi)\,dP=\frac12\int F\Phi\,Q_p\,dP,
 \qquad\mathcal L_{-p}(F\Phi)=(\mathcal L_{-p}F)\Phi+F\,\mathcal L_{-p}\Phi .
\end{equation}
For a quantity \(\lambda\alpha+\bar\lambda\beta\) the operation
\(\frac12(\cdot|_{\lambda=1}-i\cdot|_{\lambda=i})\) returns \(\alpha\). Applied to
\eqref{eq:Wp-wfield} and to \(\delta^\lambda_pf\), this gives, for \(F\in\mathcal C\) and
holomorphic \(\Phi\in\mathbb C[a_1,a_2,\ldots]\),
\begin{align}
 \mathcal L_{-p}F&=-i\sum_{j=1}^pE_{p,j}\,\mathfrak d_{-j}F-\frac i2e_p\,\mathfrak d_0F,
 \label{eq:Wp-LF}\\
 \mathcal L_{-p}\Phi&=d\Phi[\Lambda_pf],\qquad
 \Lambda_pf=-f'\,\Pi_{\ge2}\Bigl[\frac{f^{1-p}}{f'}\Bigr]-\frac{e_p}2\,(zf'-f).
 \label{eq:Wp-LPhi}
\end{align}
(In \eqref{eq:Wp-wfield} the negative modes are complex-linear in \(\lambda\), the
positive ones antilinear, and \(\dot\theta_0=\Ima(\lambda e_p)\) contributes
\(-\frac i2e_p\). In \(\delta^\lambda_pf=f'Y-(\Rea U_1)f\) the terms
\(\overline{U_{2-m}}\) are antilinear and \(\Rea U_1=-\Rea(\lambda e_p)\) contributes
\(-\frac12e_p\).) Note \(\Lambda_pf=O(z^2)\), so \(\Lambda_pf\) is a variation within
normalized maps.

The fourth step is algebraic. The Ward fields \(w^\lambda_p\) are not the
Fourier fields \(e^{\pm in\theta}\); the next lemmas give the triangular change
of basis between them and the identity that produces the ghost charge
\(26\).

\begin{lemma}\label{lem:Wp-basis}
For \(n\ge1\): \(\sum_{p=j}^nM_{n,p}E_{p,j}=\delta_{nj}\) \((1\le j\le n)\), and
\(\sum_{p=1}^nM_{n,p}P_p=-[z^{n-2}]\mathcal Sf\).
\end{lemma}
\begin{proof}
Substituting \(z=\psi(w)\) in the residue,
\[
 E_{p,j}=\operatorname{Res}_{z=0}\frac{z^{j-2}f^{1-p}}{f'}\,dz
 =\operatorname{Res}_{w=0}\psi^{j-2}w^{1-p}\psi'^2\,dw=[w^{p-2}](\psi^{j-2}\psi'^2).
\]
Hence \(\sum_pE_{p,j}f^{p-2}f'^2=(\psi^{j-2}\psi'^2)\circ f\cdot f'^2=z^{j-2}\) near
\(0\), and taking \([z^{n-2}]\) gives \(\sum_pM_{n,p}E_{p,j}=\delta_{nj}\), the
terms with \(p>n\) being zero. Likewise
\(\sum_pM_{n,p}P_p=[z^{n-2}]\bigl((\mathcal S\psi\circ f)f'^2\bigr)=-[z^{n-2}]\mathcal Sf\)
by the chain rule \(\mathcal S(\psi\circ f)=(\mathcal S\psi\circ f)f'^2+\mathcal Sf=0\);
the series \(\mathcal S\psi(w)=\sum_{p\ge2}P_pw^{p-2}\) converges near \(0\).
\end{proof}

\begin{lemma}[Ghost identity]\label{lem:Wp-ghost}
For every normalized univalent \(f\) and \(n\ge1\),
\[
 \sum_{p=1}^n\Bigl\{\Lambda_pM_{n,p}-\frac{n+p}2\,e_pM_{n,p}\Bigr\}
 =[z^{n-2}]\Bigl\{\frac{26}{12}\,\mathcal Sf+\frac{f'^2}{f^2}-\frac1{z^2}\Bigr\},
\]
where \(\Lambda_pM_{n,p}=\frac d{d\epsilon}|_{0}M_{n,p}(f+\epsilon\Lambda_pf)\) is the derivative of
the polynomial \(M_{n,p}\) in the direction \(\Lambda_pf\).
\end{lemma}

\begin{proof}
For a holomorphic field \(X\) and a quadratic differential \(\phi\,dz^2\) put
\(\mathcal L_X\phi=X\phi'+2X'\phi\). Let \(\phi_p=f^{p-2}f'^2\). A direct computation
gives \(\frac d{d\epsilon}(f+\epsilon\delta)^{p-2}(f'+\epsilon\delta')^2|_0=\mathcal L_{\delta/f'}\phi_p\),
\(\mathcal L_{f/f'}\phi_p=p\,\phi_p\) and \([z^{n-2}]\mathcal L_z\phi=n[z^{n-2}]\phi\). Hence the
left side equals \([z^{n-2}]\sum_p\mathcal L_{Y_p}\phi_p\) with
\[
 Y_p=\frac{\Lambda_pf}{f'}-\frac{e_p}2\Bigl(z+\frac f{f'}\Bigr)
    =-\Pi_{\ge2}\Bigl[\frac{f^{1-p}}{f'}\Bigr]-e_pz=-\Pi_{\ge1}\Bigl[\frac{f^{1-p}}{f'}\Bigr],
\]
by \eqref{eq:Wp-LPhi}. Since \(Y_p=O(z)\) and \(\phi_p=O(z^{p-2})\), the terms with
\(p>n\) contribute nothing to \([z^{n-2}]\), and we may sum over all \(p\ge1\).

Fix \(0<r<1\) and \(r'>0\) so small that \(|f(x)|<\min_{|t|=r}|f(t)|\) for
\(|x|\le r'\). For \(0<|z|<r'\), Cauchy's formula gives
\(Y_p(z)=-\frac1{2\pi i}\oint_{|t|=r}\frac{f(t)^{1-p}}{f'(t)}k(t,z)\,dt\) with
\(k(t,z)=\frac1{t-z}-\frac1t\). The geometric series
\(\sum_{p\ge1}f(t)^{1-p}f(x)^{p-2}=\frac1{f(t)-f(x)}+\frac1{f(x)}\) converges uniformly for
\(|t|=r\), \(|x|\le r'\). Hence
\[
 B(z):=\sum_{p\ge1}\mathcal L_{Y_p}\phi_p(z)
 =\frac1{2\pi i}\oint_{|t|=r}\Bigl\{k(t,z)\,\partial_x\mathcal H(t,x)\big|_{x=z}
   +2\,\partial_zk(t,z)\,\mathcal H(t,z)\Bigr\}dt,
\]
\[
 \mathcal H(t,x)=-\frac{f'(x)^2}{f'(t)}\Bigl(\frac1{f(t)-f(x)}+\frac1{f(x)}\Bigr),
\]
and \([z^{n-2}]\sum_p\mathcal L_{Y_p}\phi_p=[z^{n-2}]B\) (uniform convergence on
\(|z|=r'/2\)). In \(|t|<r\) the integrand is meromorphic with poles at most at
\(t=0\) and \(t=z\) (univalence). At \(t=0\), \(\mathcal H(\cdot,x)\) is holomorphic
with \(\mathcal H(0,x)=0\), so \(\partial_x\mathcal H(0,x)=0\) and the residue vanishes. At
\(t=z\) write \(\mathcal H=\mathcal K+\mathcal R\), with
\(\mathcal R(t,x)=-f'(x)^2/(f'(t)f(x))\) regular and
\[
 \mathcal K(t,x)=-\frac{f'(x)^2}{f'(t)(f(t)-f(x))}
 =-\frac1{t-x}+\mathsf A(x)+\mathsf B(x)(t-x)+O((t-x)^2),
\]
\(\mathsf A=\frac32\frac{f''}{f'}\), \(\mathsf B=\frac23\frac{f'''}{f'}-\frac74\Bigl(\frac{f''}{f'}\Bigr)^2\);
this follows from
\(f'(t)(f(t)-f(x))=f'^2s\{1+\frac32\frac{f''}{f'}s+(\frac23\frac{f'''}{f'}+\frac12(\frac{f''}{f'})^2)s^2+O(s^3)\}\),
\(s=t-x\), all derivatives at \(x\). Thus
\(\partial_x\mathcal K(t,x)|_{x=z}=-(t-z)^{-2}+\mathsf A'(z)-\mathsf B(z)+O(t-z)\), and
\[
 \operatorname{Res}_{t=z}\{k\,\partial_x\mathcal K\}=\mathsf A'-\mathsf B-\frac1{z^2},\quad
 \operatorname{Res}_{t=z}\{2\partial_zk\,\mathcal K\}=2\mathsf B,
\]
\[
 \operatorname{Res}_{t=z}\{k\,\partial_x\mathcal R\}=\partial_x\mathcal R(z,z)
   =-\frac{2f''}f+\frac{f'^2}{f^2},\quad
 \operatorname{Res}_{t=z}\{2\partial_zk\,\mathcal R\}=2\partial_t\mathcal R(z,z)=\frac{2f''}f,
\]
using \(\partial_zk=(t-z)^{-2}\) and \(\operatorname{Res}_{t=z}t^{-1}(t-z)^{-2}=-z^{-2}\).
Summing,
\(B=\mathsf A'+\mathsf B+\frac{f'^2}{f^2}-\frac1{z^2}
=\frac{13}6\bigl(\frac{f'''}{f'}-\frac32(\frac{f''}{f'})^2\bigr)+\frac{f'^2}{f^2}-\frac1{z^2}\).
\end{proof}

The diagonal value \(\mathsf A'+\mathsf B=\frac{13}6\mathcal Sf\) is the analogue of
\cite[Lemma~5.5]{BJweld}; the factor \(26/12=2\cdot\frac{13}{12}\) is the ghost
charge. Lemmas~\ref{lem:Wp-basis} and~\ref{lem:Wp-ghost} and the assembly in
the proof of Theorem~\ref{thm:Wproof} below were also checked by computer, in exact
arithmetic for \(n\le6\) and in 40-digit floating point for \(n\le9\)
(ancillary scripts; see the README).

\phantomsection\label{sec:Wp-assembly}%
We can now assemble the proof of Theorem~\ref{thm:Wproof}.
Fix \(n\ge1\) and a real \(F\in\mathcal C\). By \eqref{eq:Wp-LF} and
Lemma~\ref{lem:Wp-basis}, pointwise on \(\mathcal W\),
\[
 \mathfrak d_{-n}F=\sum_{p=1}^nM_{n,p}\sum_{j=1}^pE_{p,j}\mathfrak d_{-j}F
 =\sum_{p=1}^nM_{n,p}\Bigl\{i\,\mathcal L_{-p}F-\frac{e_p}2\,\mathfrak d_0F\Bigr\}.
\]
By \eqref{eq:Wp-complexWard} with \(\Phi=M_{n,p}\), and by \eqref{eq:Wp-rot} with
\(H=Fe_pM_{n,p}\) and \(\mathfrak d_0(e_pM_{n,p})=in\,e_pM_{n,p}\) (weights
\(p\) and \(n-p\)),
\[
 \int M_{n,p}\,\mathcal L_{-p}F\,dP=\int F\Bigl\{\frac12M_{n,p}Q_p-\mathcal L_{-p}M_{n,p}\Bigr\}dP,\qquad
 -\frac12\int e_pM_{n,p}\,\mathfrak d_0F\,dP=\frac{in}2\int Fe_pM_{n,p}\,dP .
\]
With \eqref{eq:Wp-LPhi}, \(\mathcal L_{-p}M_{n,p}=\Lambda_pM_{n,p}\). Therefore
\[
 \int\mathfrak d_{-n}F\,dP=i\int F\sum_{p=1}^n\Bigl\{\frac p2e_pM_{n,p}
 -\frac{\cm}{12}M_{n,p}P_p-\Lambda_pM_{n,p}+\frac n2e_pM_{n,p}\Bigr\}dP .
\]
By Lemmas~\ref{lem:Wp-basis} and~\ref{lem:Wp-ghost} the sum equals
\[
 \frac{\cm}{12}[z^{n-2}]\mathcal Sf-[z^{n-2}]\Bigl\{\frac{26}{12}\mathcal Sf+\frac{f'^2}{f^2}-\frac1{z^2}\Bigr\}
 =-[z^{n-2}]\Bigl\{\frac{\cL}{12}\mathcal Sf+\frac{f'^2}{f^2}-\frac1{z^2}\Bigr\}=-\mathfrak t_n,
\]
since \(26-\cm=\cL\). Hence \(\int\mathfrak d_{-n}F\,dP=-i\int F\mathfrak t_n\,dP\), and
\eqref{eq:Wp-complex} proves the theorem. \qed

\begin{remark}[Checks]\label{rem:Wp-checks}
(a) For \(n=1\) only the translation \(p=1\) enters; \(M_{1,1}=1\), \(e_1=-2a_1\),
and the theorem gives \(\int D_{\cos\theta}F=\int F\Ima\mathfrak t_1\),
\(\int D_{\sin\theta}F=\int F\Rea\mathfrak t_1\) with \(\mathfrak t_1=2a_1\). This is
the first-mode identity used in Section~\ref{sec:finiteinput}; its
agreement with the exact M\"obius likelihood \eqref{eq:mobiuscheck}, which
does not pass through the Ward identities, checks the orientation
conventions in \eqref{eq:Wp-velocity}. (b) The sign: by the score covariance (Appendix~\ref{app:covariance}), \(\mathfrak d_n\mathfrak t_n=iC_n\) and
\(\int\mathfrak d_nF=i\int F\bar{\mathfrak t}_n\) give \(\E|\mathfrak t_n|^2=C_n>0\); the
opposite sign would give \(-C_n\). (c) The complex form
\(\int\mathfrak d_nF=i\int F\bar{\mathfrak t}_n\) of \eqref{eq:Wp-complex}; for its relation to
\cite[Theorem~1.2]{BJweld} see Section~\ref{sec:conventionsLit}, (C1).
(d) The Ward-flow derivation \eqref{eq:Wp-LPhi} is the response of
Lemma~\ref{lem:coeffResp} at the velocity \eqref{eq:Wp-wfield}:
\(\Lambda_pf=\sum_{j=1}^pE_{p,j}\mathbf L_{-j}f+\frac{e_p}2\mathbf L_0f\), i.e.\
\(\mathcal L_{-p}\) acts on holomorphic \(\Phi\) by the formula \eqref{eq:Wp-LF} with
\(\mathfrak d_j=iL_j\). Indeed \(\sum_{j=0}^pE_{p,j}w^{1-j}=(1-\Pi_{\ge2})[f^{1-p}/f']\). In the
integral defining \(\mathbf L_{-j}f(z)\), the term \(f^{1-p}/f'\) contributes
\(\oint\zeta^{-1-p}(\zeta-f(z))^{-1}d\zeta=0\) (\(\zeta=f(w)\); residue at \(\infty\)), and
\(\Pi_{\ge2}[f^{1-p}/f']=O(w^2)\) only its residue at \(w=z\). So
\(\sum_{j=0}^pE_{p,j}\mathbf L_{-j}f=-f'\Pi_{\ge2}[f^{1-p}/f']\); now use \(E_{p,0}=e_p\),
\(\mathbf L_0f=zf'-f\).
\end{remark}

\begin{remark}[Status]\label{rem:Wp-status}
The proof is complete relative to (P1)--(P4) and classical complex analysis.
We record the points where it rests on statements we have not reproved.
(i) Proposition~\ref{prop:Wp-loopWard} for \(p\ge2\) uses
Lemma~\ref{lem:normalization} and Lemma~\ref{lem:GQWext}.
Lemma~\ref{lem:normalization} rests on the proof of \cite[Theorem~5.1]{ZhanLoop}
(Step~2 computes the constant of its bubble decomposition) and on the boundary
perturbation rule \cite[(1.4)]{FieldLawler}; see (C4) of
Section~\ref{sec:conventionsLit}. Lemma~\ref{lem:GQWext} rests on
Lemma~\ref{lem:GQWanalytic}, i.e.\ on
\cite[Corollary~2.6, Theorems~2.2--2.3]{CarfagniniWang},
\cite[Theorem~4.1]{WangAIF} and \cite[Theorem~1.4]{WangInv}. (ii) For \(p=1\) and for the rotation identity
\eqref{eq:Wp-rot} only the M\"obius invariance \cite[Theorem~4.2(ii)]{ZhanLoop}
is used. (iii) From \cite{BJweld} only the definition of \(P\) is used
(Section~\ref{sec:setting}). (iv) The classical inputs are the
measurable Riemann mapping theorem with its \(L^p\) construction, uniqueness
of normalized solutions, and the composition rule for Beltrami coefficients
\cite{AhlforsLectures,AhlforsBers,LehtoVirtanen}, together with Koebe's
theorem and de Branges' bound. (v) Measurability of \(H\) in
Proposition~\ref{prop:Wp-loopWard} and of the maps \(\gamma\mapsto(c(\gamma),\eta(\gamma),W(\gamma,\alpha))\)
is for the Hausdorff metric on Jordan curves: \(\gamma\mapsto f_\gamma\) is
continuous for locally uniform convergence (Carath\'eodory kernel theorem), and
the boundary values \(f_\gamma(e^{i\theta})=\lim_{r\uparrow1}f_\gamma(re^{i\theta})\)
\cite[Theorem~2.6]{PommerenkeBB} are pointwise limits of continuous functions
of \((\gamma,\theta)\); likewise for \(g_\gamma\). Hence \(W\) and all coefficient
functionals used are Borel.
\end{remark}
\section{Relative entropy equals the Liouville--capacity action}
\label{sec:entropy}

For \(\Phi\in X\) put \(\mathsf R_\Phi(h)=h\circ\Phi\) and
\(P^{\mathsf R}_\Phi=(\mathsf R_\Phi)_*P\). For a welding
\(y=g^{-1}\circ f|_{S^1}\) of a Jordan curve, with \(f(0)=0\),
\(g(\infty)=\infty\), let
\[
 H_\kappa=\frac{\cL}{24\pi}S_1+2K,\qquad K=\log\left|\frac{g'(\infty)}{f'(0)}\right|,\qquad
 S_1=\int_\D\Bigl|\frac{f''}{f'}\Bigr|^2\dd+\int_{\Ds}\Bigl|\frac{g''}{g'}\Bigr|^2\dd-4\pi K,
\]
the Liouville--capacity action, whose exponential \(e^{-H_\kappa}\) is the
formal density of \cite[(1.9)--(1.10), (1.14)]{BJweld}; \(S_1=\pi I^L\) is
the universal Liouville action of the welded curve \cite{TT,WangInv},
finite exactly on Weil--Petersson (WP) weldings, and the normalization is
that of Proposition~\ref{prop:potential}. Both terms are unchanged under
common affine post-composition of \((f,g)\). Hence \(H_\kappa(\mathrm{id})=0\),
and \(H_\kappa\) is invariant under left rotations (which rotate \(g\)) and
under right rotations (which replace \((f,g)\) by
\((e^{-i\theta}f(e^{i\theta}\cdot),e^{-i\theta}g)\)). Since the inverse welds the
reflected curve \(\iota(\gamma)\), \(\iota(z)=1/\bar z\), \(K\circ\mathsf J=K\)
(Lemma~\ref{lem:CaiDictionary}), and the loop Loewner energy is invariant
under M\"obius maps and complex conjugation \cite{WangInv} (equivalently,
\(S_1\) is symmetric \cite[Theorem~4.8]{TT}), \(H_\kappa(\Phi^{-1})=H_\kappa(\Phi)\).

The law \(P\) is invariant under left and right rotations
(Lemma~\ref{lem:symmetries}). Write
\(\|Q-R\|_{\rm TV}=\sup_B|Q(B)-R(B)|\). In this section and in
Appendix~\ref{app:caiInput}, \(a_m(h)\) denotes the coefficients of the
interior map \(f_h(z)=z(1+\sum_ma_mz^m)\); by \eqref{eq:stressdictionary}
they are the conjugates of the coefficients of \eqref{eq:G} for
\(\mathsf J(h)\).

\begin{theorem}[Entropy--action identity]\label{thm:entropyWP}
Let \(0<\kappa\le4\) and let \(\Phi\) be a Weil--Petersson
homeomorphism. Then
\[
 \operatorname{Ent}(P^{\mathsf R}_\Phi\,|\,P)
 =\operatorname{Ent}(P\,|\,P^{\mathsf R}_\Phi)
 =H_\kappa(\Phi)
 =\frac{\cL}{24}I^L(\Phi)+2K(\Phi).
\]
In particular \(P^{\mathsf R}_\Phi\sim P\), and for the left action
\(\operatorname{Ent}(P_a|P)=\operatorname{Ent}(P|P_a)=H_\kappa(a)\) for every
WP homeomorphism \(a\). At \(\kappa=4\) this proves quasi-invariance of
the SLE\(_4\) welding law under all WP maps.
\end{theorem}

The proof has two steps: an exact identity for a dense group of
finite-mode maps, obtained by a mean-value argument along radial
holomorphic discs, and a two-sided approximation. Its probabilistic input
is the integration-by-parts formula of Theorem~\ref{thm:Wproof} (all modes
\(n\ge1\)), turned into densities by Lemma~\ref{lem:measureODE}, together with
the symmetries of \(P\) (Lemma~\ref{lem:symmetries}). Everything else is deterministic and is
proved in Appendix~\ref{sec:entropyDet} from classical complex analysis and
the published results of Takhtajan--Teo on the WP topology. Nothing from the
unrefereed preprint discussed in Remark~\ref{rem:entCai} is used.

We begin with the finite-mode identity, for maps that are finite
compositions of Fourier flows.
A word is a finite sequence \(w=(\psi^{v_1}_{t_1},\ldots,\psi^{v_\ell}_{t_\ell})\) of flows
of real trigonometric polynomials (rotations included); put
\(\Phi_{w,i}=\psi^{v_1}_{t_1}\circ\cdots\circ\psi^{v_i}_{t_i}\), \(\Phi_w=\Phi_{w,\ell}\), and
\[
 B_w(h)=\sum_{i=1}^\ell\int_0^{t_i}\rho_{v_i}\bigl(h\circ\Phi_{w,i-1}\circ\psi^{v_i}_s\bigr)\,ds,
 \qquad h\in X,
\]
with the bounded Borel versions of Section~\ref{sec:finiteinput}. Every element
of \(\mathcal G_{\rm fin}\) is some \(\Phi_w\), and concatenation gives
\(B_{ww'}(h)=B_w(h)+B_{w'}(h\circ\Phi_w)\).

For a normalized univalent \(f\) on \(\D\) and \(\zeta\in\D\), put
\(F^f_\zeta(z)=\zeta^{-1}f(\zeta z)\) (\(F^f_0=\mathrm{id}\)), univalent on
\(\{|z|<1/|\zeta|\}\); its boundary is an analytic curve. For \(h\in\mathcal W\) let
\(x^h_\zeta\) be a welding with interior map \(F^{f_h}_\zeta\); it is determined up to a
left rotation, \(a_m(x^h_\zeta)=\zeta^ma_m(h)\), \(x^h_0=\mathrm{id}\), and \(x^h_\zeta\) is
real analytic, since \(F^{f_h}_\zeta\) and, by Schwarz reflection across the
analytic boundary curve, the exterior map extend conformally across \(S^1\). Since \(f_{h\circ R_\theta}(z)=e^{-i\theta}f_h(e^{i\theta}z)\), the welding \(x^h_{re^{i\theta}}\)
equals \(x^{h\circ R_\theta}_r\) modulo left rotations. The weldings \(x^h_r\),
\(0<r<1\), are the radial approximations of \(h\) (Remark~\ref{rem:entCai}).

\begin{lemma}\label{lem:Bfin}
Let \(0<\kappa\le4\) and let \(w\) be a word.
\begin{enumerate}
\item \(|B_w|\le\sum_i|t_i|\,\|\rho_{v_i}\|_\infty\), and
\(dP^{\mathsf R}_{\Phi_w}/dP(y)=\exp B_w(y\circ\Phi_w^{-1})\) for \(P\)-almost every \(y\).
\item On \(\mathcal W\), \(B_w\) is invariant under left rotations, and for every
\(h\in\mathcal W\), \(B_w(h)=\lim_{r\uparrow1}B_w(x^h_r)\). The function
\((f,\zeta)\mapsto B_w(x_\zeta)\), \(x_\zeta\) a welding with interior map \(F^f_\zeta\), is
continuous for locally uniform convergence of \(f\) and \(\zeta\in\D\).
\item For every real-analytic diffeomorphism \(\hat h\),
\(B_w(\hat h)=H_\kappa(\hat h\circ\Phi_w)-H_\kappa(\hat h)\). Hence
\(B_w(\mathrm{id})=H_\kappa(\Phi_w)\), and, by (ii), for \(h\in\mathcal W\),
\(B_w(h)=\lim_{r\uparrow1}[H_\kappa(x^h_r\circ\Phi_w)-H_\kappa(x^h_r)]\). In particular, on
\(\mathcal W\), \(B_w\) depends only on \(\Phi_w\) and satisfies
\(B_{\Phi\circ\Psi}(h)=B_\Phi(h)+B_\Psi(h\circ\Phi)\).
\item For every \(h\in\mathcal W\), \(\zeta\mapsto B_w(x^h_\zeta)\) is harmonic on \(\D\).
\end{enumerate}
\end{lemma}

\begin{proof}
(i) The bound follows from the definition. For a single letter, let
\(T_th=h\circ\psi^v_t\). The cylinder algebra \(\mathcal C\) of
Section~\ref{sec:finiteinput} determines finite signed measures, is stable
under \(T_t\) and has uniform-norm \(C^1\) orbits, with \(D_vF\) real-linear in
\(v\) by \eqref{eq:corestability}. Theorem~\ref{thm:Wproof} gives
\(\int D_{\cos n\theta}F\,dP=\int F\Ima\mathfrak t_n\,dP\) and
\(\int D_{\sin n\theta}F\,dP=\int F\Rea\mathfrak t_n\,dP\) for all \(n\ge1\), and
\(\int D_1F\,dP=0\) by right-rotation invariance \eqref{eq:Wp-rot}. Hence
\(\int D_vF\,dP=\int F\rho_v\,dP\) on \(\mathcal C\), and Lemma~\ref{lem:measureODE}
gives \(dP^{\mathsf R}_{\psi^v_t}/dP(y)=\exp\int_0^t\rho_v(y\circ\psi^v_{-s})\,ds
=\exp B_{\psi^v_t}(y\circ\psi^v_{-t})\). For a word \(w=w'\psi\) with last letter \(\psi\),
\(P^{\mathsf R}_{\Phi_w}=(\mathsf R_\psi)_*P^{\mathsf R}_{\Phi_{w'}}\), and the change of variables
\(\int F(u\circ\psi)Z(u)\,dP(u)=\int F(u)Z(u\circ\psi^{-1})\,dP^{\mathsf R}_\psi(u)\) gives
\(Z_{\Phi_w}(u)=Z_{\Phi_{w'}}(u\circ\psi^{-1})Z_\psi(u)=\exp\{B_{w'}(h)+B_\psi(h\circ\Phi_{w'})\}\),
\(h=u\circ\Phi_w^{-1}\). Induction on \(\ell\) proves the claim.

(ii) By Lemma~\ref{lem:entQ}, \(h\circ\Phi_{w,i-1}\circ\psi^{v_i}_s\in\mathcal W\), so on \(\mathcal W\) all
coefficients in \(B_w\) are actual coefficients; left rotations do not change
\(f_{h\circ\Psi}\). Fix \(i\) and \(s\), and \(\Psi=\Phi_{w,i-1}\circ\psi^{v_i}_s\). As \(r\uparrow1\),
\(F^{f_h}_r\to f_h\) uniformly on \(\overline\D\), because \(f_h\) is uniformly continuous on
\(\overline\D\); Lemma~\ref{lem:entC} gives
\(\rho_{v_i}(x^h_r\circ\Psi)\to\rho_{v_i}(h\circ\Psi)\), \(\rho_{v_i}\) being a polynomial in
finitely many coefficients. The same lemma, with \(\Psi_k=\Phi_{w,i-1}\circ\psi^{v_i}_{s_k}\),
shows that the integrands are continuous in \(s\); they are bounded, so
dominated convergence gives the limit. If \(f_k\to f\) locally uniformly and
\(\zeta_k\to\zeta\in\D\), then \(F^{f_k}_{\zeta_k}\to F^f_\zeta\) uniformly on \(\overline\D\), and
Lemma~\ref{lem:entC} with dominated convergence gives the continuity
statement.

(iii) Each \(\hat h\circ\Phi_{w,i-1}\circ\psi^{v_i}_s\) is real analytic. By
Lemma~\ref{lem:entV} applied at these points and the flow property,
\(s\mapsto H_\kappa(\hat h\circ\Phi_{w,i-1}\circ\psi^{v_i}_s)\) is differentiable everywhere with
derivative \(\rho_{v_i}(\hat h\circ\Phi_{w,i-1}\circ\psi^{v_i}_s)\), which is bounded. Such a
function is Lipschitz, so the fundamental theorem of calculus applies, and
summing over \(i\) telescopes. The welding \(x^h_r\) is real analytic, which
gives the limit formula; its right side does not involve the word.
Concatenation gives the cocycle.

(iv) Fix \(i\), \(s\) and \(\Psi\) as in (ii). The family \(F^{f_h}_\zeta\), \(\zeta\in\D\), satisfies
the hypothesis of Lemma~\ref{lem:entH} (near \(\zeta_0\) take \(1<\varrho<1/|\zeta_0|\)
and \(|\zeta|<1/\varrho\); \(F^{f_h}_\zeta(z)=z(1+\sum_ma_m(\zeta z)^m)\) is jointly holomorphic). So
each \(a_m(x^h_\zeta\circ\Psi)\) is holomorphic in \(\zeta\); each \(\mathfrak t_n\) is a holomorphic
polynomial in the \(a_m\), so \(\rho_{v_i}(x^h_\zeta\circ\Psi)\) is the real part of a
holomorphic function of \(\zeta\). By Lemma~\ref{lem:entC} (with
\(F_k=F^{f_h}_{\zeta_k}\) and \(\Psi_k=\Phi_{w,i-1}\circ\psi^{v_i}_{s_k}\)) these integrands are jointly
continuous in \((s,\zeta)\), and they are bounded; Fubini's theorem transfers the mean-value
property on circles to \(\zeta\mapsto B_w(x^h_\zeta)\), which is continuous, hence harmonic.
\end{proof}

\begin{corollary}[M\"obius letters]\label{cor:entMobius}
For every M\"obius automorphism \(m\) of \(\D\), for \(P\)-almost every \(y\) and \(x\),
\[
 \frac{dP^{\mathsf R}_m}{dP}(y)=e^{2K(y)-2K(y\circ m^{-1})},\qquad
 \frac{dP_m}{dP}(x)=e^{2K(x)-2K(m^{-1}\circ x)} .
\]
Equivalently, \(e^{2K}P\) is invariant under right and left M\"obius composition.
\end{corollary}
\begin{proof}
Every M\"obius automorphism is \(R_\alpha\circ m_t\), and
\(m_t(z)=(z+\bar t)/(1+tz)\) is the time-one map of a field
\(v=x_1\cos\theta+y_1\sin\theta\) (conjugate the one-parameter group \(m_{\tanh s}\), generated by
\(-2\sin\theta\), by a rotation); \(P\) and \(K\) are rotation invariant. For such \(v\),
\(\rho_v=x_1\Ima\mathfrak t_1+y_1\Rea\mathfrak t_1\) with \(\mathfrak t_1=k_1\), so by
Lemma~\ref{lem:entV}, valid on \(\mathcal W\) for \(K\), \(s\mapsto2K(h\circ\psi^v_s)\) has
derivative \(\rho_v(h\circ\psi^v_s)\) everywhere, for \(h\in\mathcal W\). This derivative is
bounded, so \(B_{\psi^v_t}(h)=2K(h\circ\psi^v_t)-2K(h)\), and Lemma~\ref{lem:Bfin}(i) gives the
first formula. The second follows from \eqref{eq:RNconjugacy}, \(K\circ\mathsf J=K\)
and \((y^{-1}\circ m)^{-1}=m^{-1}\circ y\).
\end{proof}

This is \cite[Lemma~5.2]{BJweld}, here a consequence of
Theorem~\ref{thm:Wproof}.

\begin{proposition}[Finite-mode identity]\label{prop:entropyFin}
For \(0<\kappa\le4\) and every word \(w\),
\(\operatorname{Ent}(P^{\mathsf R}_{\Phi_w}|P)=\E_P[B_w]=H_\kappa(\Phi_w)\).
\end{proposition}
\begin{proof}
By Lemma~\ref{lem:Bfin}(i), \(\operatorname{Ent}(P^{\mathsf R}_{\Phi_w}|P)
=\int B_w(y\circ\Phi_w^{-1})\,dP^{\mathsf R}_{\Phi_w}(y)=\E_P[B_w]\). Since \(P(\mathcal W)=1\) (Section~\ref{sec:finiteinput}),
Lemma~\ref{lem:Bfin}(ii) and bounded convergence give
\(\E_P[B_w]=\lim_{r\uparrow1}\E_P[B_w(x^h_r)]\). Here \(B_w(x^h_\zeta)=\beta(f_h,\zeta)\) depends
only on the interior map, so on the loop side
\(\E_P[B_w(x^h_\zeta)]=\int\beta(f_\eta,\zeta)\,d\mu_\kappa(\eta)\), where \(\eta\mapsto f_\eta\) is Borel
(Remark~\ref{rem:Wp-status}(v)) and \(\beta\) is jointly continuous by
Lemma~\ref{lem:Bfin}(ii). For fixed \(r<1\), right-rotation invariance of \(P\),
the identity \(x^{h\circ R_\theta}_r=x^h_{re^{i\theta}}\) modulo left rotations, and
left-rotation invariance of \(B_w\) give
\(\E_P[B_w(x^h_r)]=\E_P[B_w(x^h_{re^{i\theta}})]\) for every \(\theta\). Averaging in \(\theta\)
(Fubini, \(\beta\) being bounded and jointly measurable) and using the mean-value
property of Lemma~\ref{lem:Bfin}(iv),
\(\E_P[B_w(x^h_r)]=\E_P[B_w(x^h_0)]=B_w(\mathrm{id})=H_\kappa(\Phi_w)\) by
Lemma~\ref{lem:Bfin}(iii).
\end{proof}

We now prove Theorem~\ref{thm:entropyWP}, by approximation from the
finite-mode identity.
Both sides of the identity depend on \(\Phi\) only modulo left rotations:
\(H_\kappa\) is left-rotation invariant, and
\(P^{\mathsf R}_{R_\alpha\circ\Phi}=P^{\mathsf R}_\Phi\) by right-rotation invariance of \(P\). We
therefore work in \(\widetilde T_0(1)\), the WP part of the universal Teichm\"uller
curve \(S^1\backslash\mathrm{Homeo}_{\rm qs}(S^1)\), with representatives fixing
\(1\). It is a topological group, smooth maps are dense in it, and the
inclusion of \(\mathrm{Diff}_+^\infty(S^1)\) with its \(C^\infty\) topology is
continuous \cite[Appendix~A]{TTI}. By Lemma~\ref{lem:entD} and a diagonal
argument there are \(\Phi_j\in\mathcal G_{\rm fin}\) with \(\Phi_j\to\Phi\) in \(\widetilde T_0(1)\),
hence, after left rotations (which lie in \(\mathcal G_{\rm fin}\)), uniformly with
their inverses. \(K\) and \(S_1\) are continuous on \(\widetilde T_0(1)\)
\cite[Lemma~5.1]{TTI}, \cite[Theorems~3.10 and~4.8]{TT}.

\emph{Upper bound.} Since \(h\circ\Phi_j\to h\circ\Phi\) in \(X\) for every \(h\),
\(P^{\mathsf R}_{\Phi_j}\Rightarrow P^{\mathsf R}_\Phi\). Weak lower semicontinuity of relative
entropy and Proposition~\ref{prop:entropyFin} give
\(\operatorname{Ent}(P^{\mathsf R}_\Phi|P)\le\liminf_jH_\kappa(\Phi_j)=H_\kappa(\Phi)\), for every WP map \(\Phi\).

\emph{Lower bound.} Relative entropy is invariant under the Borel bijection
\(\mathsf R_{\Phi_j^{-1}}\), so with \(\Theta_j=\Phi\circ\Phi_j^{-1}\to\mathrm{id}\) in \(\widetilde T_0(1)\) the upper
bound gives \(\operatorname{Ent}(P^{\mathsf R}_\Phi|P^{\mathsf R}_{\Phi_j})=\operatorname{Ent}(P^{\mathsf R}_{\Theta_j}|P)
\le H_\kappa(\Theta_j)\to0\), and \(P^{\mathsf R}_{\Phi_j}\to P^{\mathsf R}_\Phi\) in total variation by
Pinsker's inequality. Fix \(\Psi\in\mathcal G_{\rm fin}\) and let
\(F=\log(dP^{\mathsf R}_\Psi/dP)\), bounded by Lemma~\ref{lem:Bfin}(i), with
\(\E_Pe^F=1\). Since all densities involved are bounded above and below,
\[
 \int F\,dP^{\mathsf R}_{\Phi_j}
 =\operatorname{Ent}(P^{\mathsf R}_{\Phi_j}|P)-\operatorname{Ent}(P^{\mathsf R}_{\Phi_j}|P^{\mathsf R}_\Psi)
 =H_\kappa(\Phi_j)-H_\kappa(\Phi_j\circ\Psi^{-1}),
\]
by Proposition~\ref{prop:entropyFin} and
\(\operatorname{Ent}(P^{\mathsf R}_{\Phi_j}|P^{\mathsf R}_\Psi)=\operatorname{Ent}(P^{\mathsf R}_{\Phi_j\circ\Psi^{-1}}|P)\) (apply
\(\mathsf R_{\Psi^{-1}}\)). Letting \(j\to\infty\), using total-variation convergence on the
left and continuity of \(H_\kappa\) and of right translation by \(\Psi^{-1}\) on the
right, \(\int F\,dP^{\mathsf R}_\Phi=H_\kappa(\Phi)-H_\kappa(\Phi\circ\Psi^{-1})\). The Donsker--Varadhan
inequality gives \(\operatorname{Ent}(P^{\mathsf R}_\Phi|P)\ge\int F\,dP^{\mathsf R}_\Phi-\log\E_Pe^F
=H_\kappa(\Phi)-H_\kappa(\Phi\circ\Psi^{-1})\). Take \(\Psi=\Phi_k\) and let \(k\to\infty\).

The second equality follows by applying the first to \(\Phi^{-1}\) and
using \(H_\kappa(\Phi^{-1})=H_\kappa(\Phi)\), since
\(\operatorname{Ent}(P|P^{\mathsf R}_\Phi)=\operatorname{Ent}(P^{\mathsf R}_{\Phi^{-1}}|P)\). Finite entropy in both
directions gives \(P^{\mathsf R}_\Phi\sim P\). For the left action,
\(\mathsf J_*P=P\) and \(\mathsf J\circ L_a=\mathsf R_{a^{-1}}\circ\mathsf J\) give
\((L_a)_*P=\mathsf J_*P^{\mathsf R}_{a^{-1}}\). \qed

\begin{corollary}[Consequences]\label{cor:entropyConsequences}
Let \(0<\kappa\le4\).
\begin{enumerate}
\item For WP maps \(a,b\),
\(\operatorname{Ent}(P_a|P_b)=H_\kappa(b^{-1}\circ a)\), and
\(\|P_a-P_b\|_{\rm TV}\le\sqrt{H_\kappa(b^{-1}\circ a)/2}\).
\item On the M\"obius family, \(\operatorname{Ent}(P_{m_t}|P)=-2\log(1-|t|^2)\)
for every \(\kappa\).
\item \(\kappa\operatorname{Ent}(P_a|P)\to I^L(a)\) as \(\kappa\to0\).
\item All laws \(P_a\), \(a\) WP, are mutually equivalent; in particular so
are the laws \(P_{h^{-1}}\) of the marked WP family of Section~\ref{sec:global},
and their mutual relative entropies are finite.
\end{enumerate}
\end{corollary}
\begin{proof}
(i) Pushing forward by \(L_{b^{-1}}\) preserves entropy; then Pinsker.
(ii) \(S_1(m_t)=0\) and \(K(m_t)=-\log(1-|t|^2)\), since the normalized pair
is \((A\circ m_t,A)\) with \(A(w)=(w-\bar t)/(1-|t|^2)\), where
\(m_t(z)=(z+\bar t)/(1+tz)\). (iii)
\(\kappa\cL/24=1+\frac{13}{24}\kappa+\frac{1}{16}\kappa^2\to1\) (Proposition~\ref{prop:kappa0}). (iv) By (i), both relative entropies of
\(P_a\) and \(P_b\) are finite, since the WP maps form a group.
\end{proof}

Part (ii) agrees with an independent fourth-order computation: with
\(X=2\Rea\log G'(t)\), \eqref{eq:mobiuscheck} and \(\E e^X=(1-|t|^2)^{-2}\) give
\(\operatorname{Ent}(P_{m_t}|P)=2\log(1-|t|^2)+(1-|t|^2)^2\E[Xe^X]\);
Proposition~\ref{prop:moments} gives
\((1-|t|^2)^2\E[Xe^X]=4|t|^2+2|t|^4+O(|t|^6)\), the dependence on \(C_2\)
cancelling, so the entropy is \(2|t|^2+|t|^4+O(|t|^6)\), as it should be. The
second-order expansion of
\(H_\kappa\) along mode \(n\) is \(C_n|t_n|^2\)
(Proposition~\ref{prop:potential}), matching \(\frac12I\). The identity is
the rigorous form of what the formal density \(e^{-H_\kappa}\) of
\cite[(1.10), (1.14)]{BJweld} predicts; the new step is the mean-value
argument. It is an entropy analogue of the Cameron--Martin formula
\(\operatorname{Ent}(\mu(\cdot-h)|\mu)=\frac12|h|^2_{\rm CM}\). Part (iii) is
consistent with the Loewner-energy large deviations as \(\kappa\to0\)
\cite{WangJEMS,PeltolaWang}: it has the same leading coefficient \(1/\kappa\)
as the exact chordal Cameron--Martin identity recalled in
Remark~\ref{rem:kappa0}. This is a consistency check of the
normalization in Lemma~\ref{lem:normalization}, not a calibration.

Theorem~\ref{thm:entropyWP} has further consequences, which use only
general properties of the relative entropy. They give the action \(H_\kappa\),
and the Loewner energy itself, a meaning at fixed \(\kappa\), in addition to
the large-deviation meaning of the Loewner energy as \(\kappa\to0\).

\begin{corollary}[Further consequences]\label{cor:entropyFurther}
Let \(0<\kappa\le4\), let \(a\) be a WP homeomorphism, and let \(\sigma\sim P\).
\begin{enumerate}
\item \emph{Symmetries.} \(P_a=P\) if and only if \(a\) is a rotation. Hence
\(a\mapsto P_a\) is injective on \(\operatorname{Rot}(S^1)\backslash\operatorname{WP}(S^1)\).
\item \emph{Testing.} Let \(\beta_n(\varepsilon)\) be the smallest probability under
\(P_a^{\otimes n}\) of a Borel set \(A\subset X^n\) with \(P^{\otimes n}(A)\ge1-\varepsilon\).
Then for every \(0<\varepsilon<1\),
\[
 \lim_{n\to\infty}-\frac1n\log\beta_n(\varepsilon)=H_\kappa(a)
 =\frac{\cL}{24}I^L(a)+2K(a).
\]
\item \emph{Recovery of the Loewner energy.} Writing \(\operatorname{Ent}_\kappa\) for
the entropy at parameter \(\kappa\), for \(\kappa\ne\kappa'\) in \((0,4]\),
\[
 I^L(a)=24\,\frac{\operatorname{Ent}_\kappa(P_a|P)-\operatorname{Ent}_{\kappa'}(P_a|P)}
 {\cL(\kappa)-\cL(\kappa')},\qquad
 2K(a)=\frac{\cL(\kappa)\operatorname{Ent}_{\kappa'}(P_a|P)-\cL(\kappa')\operatorname{Ent}_\kappa(P_a|P)}
 {\cL(\kappa)-\cL(\kappa')} .
\]
\item \emph{Monotonicity in \(\kappa\).} \(\kappa\mapsto\operatorname{Ent}_\kappa(P_a|P)\) is
non-increasing on \((0,4]\), strictly decreasing unless \(a\) is a M\"obius
map, and \(\operatorname{Ent}_\kappa(P_a|P)\ge\frac{25}{24}I^L(a)+2K(a)\), with equality at
\(\kappa=4\).
\item \emph{Variational formula.} For every bounded Borel \(F:X\to\mathbb R\),
\[
 \E\bigl[F(a\circ\sigma)\bigr]-\log\E\bigl[e^{F(\sigma)}\bigr]\le H_\kappa(a),
\]
and the supremum of the left side over such \(F\) is \(H_\kappa(a)\).
\item \emph{Partial observations.} For every Borel map \(T\) from \(X\) to a
measurable space, \(\operatorname{Ent}(T_*P_a|T_*P)\le H_\kappa(a)\); in particular this
bounds the relative entropy of the laws of
\((a(\sigma(z_1)),\ldots,a(\sigma(z_m)))\) and \((\sigma(z_1),\ldots,\sigma(z_m))\) for
any points \(z_1,\ldots,z_m\in S^1\).
\item \emph{A mean identity for the capacity.} For the disc automorphisms
\(m_t\), with \(K(\sigma)=\log|g'(\infty)/f'(0)|\) for the SLE loop welded by
\(\sigma\),
\[
 \E\bigl[K(m_t^{-1}\circ\sigma)-K(\sigma)\bigr]=-\log(1-|t|^2)
\]
for every \(\kappa\in(0,4]\).
\end{enumerate}
\end{corollary}
\begin{proof}
(i) Rotations preserve \(P\) (Lemma~\ref{lem:symmetries}). Conversely, if
\(P_a=P\) then \(H_\kappa(a)=0\). Both terms of \(H_\kappa\) are nonnegative:
\(I^L\ge0\), and for a pair \((f,g)\) of a Jordan curve \(\gamma\),
\(\pi|f'(0)|^2\le\operatorname{area}(\operatorname{int}\gamma)\le\pi|g'(\infty)|^2\), the first
inequality by the area formula
\(\operatorname{area}=\pi\sum_{n\ge1}n|[z^n]f|^2\) and the second by the area
theorem \cite[Theorem~2.1]{Duren}, with equality in both only if \(f\) is
linear. Hence \(K(a)=0\), \(\gamma\) is a circle centred at \(0\), \(f\) and \(g\) are
linear with \(|f'(0)|=|g'(\infty)|\), and \(a=g^{-1}\circ f\) is a rotation. If
\(P_a=P_b\) then \(P_{b^{-1}\circ a}=P\).
(ii) This is the Chernoff--Stein lemma \cite[Theorem~11.8.3]{CoverThomas}
for the null hypothesis \(P\) and the alternative \(P_a\), with exponent
\(\operatorname{Ent}(P|P_a)=H_\kappa(a)\); its proof uses only the law of large numbers
for \(\log(dP/dP_a)\in L^1(P)\) and applies on any measurable space.
(iii) \(H_\kappa\) is affine in \(\cL\), and \(\cL(\kappa)=13+24/\kappa+3\kappa/2\) is
strictly decreasing on \((0,4]\), so \(\cL(\kappa)\ne\cL(\kappa')\); solve the two
linear equations. (iv) follows from the same monotonicity, \(\cL(4)=25\) and
\(I^L\ge0\), with \(I^L(a)=0\) exactly for M\"obius maps. (v) is the
Donsker--Varadhan formula for \(\operatorname{Ent}(P_a|P)\), and (vi) the
data-processing inequality. (vii) By Corollary~\ref{cor:entMobius},
\(\log(dP/dP_{m_t})(x)=2K(m_t^{-1}\circ x)-2K(x)\), which is \(P\)-integrable
since \(\operatorname{Ent}(P|P_{m_t})<\infty\); its \(P\)-mean is
\(\operatorname{Ent}(P|P_{m_t})=-2\log(1-|t|^2)\) by
Corollary~\ref{cor:entropyConsequences}(ii).
\end{proof}

\begin{remark}[Meaning of the consequences]\label{rem:entropyMeaning}
Part (ii) says that \(H_\kappa(a)\) is the optimal exponential rate at which
\(n\) independent weldings distinguish \(P\) from \(P_a\); by the symmetry
\(\operatorname{Ent}(P_a|P)=\operatorname{Ent}(P|P_a)\) the rate is the same when the roles
of the hypotheses are exchanged. This symmetry of the relative entropy is
not a general property of statistical models; it holds for the
translates of a Gaussian measure, and here it reflects
\(H_\kappa(a^{-1})=H_\kappa(a)\). Part (iii) defines the Loewner energy and the
capacity term from the welding laws at two values of \(\kappa\), without a
limit \(\kappa\to0\). By part (iv) the SLE\(_4\) welding is the least sensitive
to deformations. Part (v) gives lower bounds for
\(\frac{\cL}{24}I^L+2K\) by expectations of SLE observables. Since
\(V\circ\mathsf R_b-V\) is pluriharmonic (Section~\ref{sec:global}), the defect
\(\operatorname{Ent}(E_{h\circ b}|P)-\operatorname{Ent}(E_h|P)\) is pluriharmonic in \(h\). The
convergence \(\kappa\operatorname{Ent}_\kappa(P_a|P)\to I^L(a)\) does not by itself give
a large deviation principle for the welding law. Conditionally on
\cite{FanCM}, finiteness of \(\operatorname{Ent}(P_a|P)\) characterizes the WP class
(Proposition~\ref{prop:FanComparison}(iv)), and part (i) holds for every
\(a\in\operatorname{Homeo}_+(S^1)\).
\end{remark}

\begin{remark}[Sources]\label{rem:entCai}
The radial approximations \(x^h_r\), and the description of the
log-likelihood of a finite-mode right composition as the radial limit of
action increments (Lemma~\ref{lem:Bfin}(iii)), appear in Cai's preprint
\cite[(1.32)--(1.33)]{CaiResponse}, together with versions of the
deterministic Lemmas~\ref{lem:entC} and~\ref{lem:entV} and of the densities in
Lemma~\ref{lem:Bfin}(i). Our proofs do not use \cite{CaiResponse}: the
densities come from Theorem~\ref{thm:Wproof} and Lemma~\ref{lem:measureODE};
Lemmas~\ref{lem:entR}--\ref{lem:entD} are proved here from classical
quasiconformal theory; and the M\"obius letters and
\cite[Lemma~5.2]{BJweld} are not needed, because Theorem~\ref{thm:Wproof}
covers the mode \(n=1\); Corollary~\ref{cor:entMobius} recovers
\cite[Lemma~5.2]{BJweld}. The properties of \(P\) used are those of
Section~\ref{sec:setting}.
\end{remark}

As our results are also closely related to the concurrent and independent work of S.~Fan \cite{FanCM}, it is worthwhile to describe the connections precisely. All parts of the next proposition (\ref{prop:FanComparison}) rely on the preprint \cite{FanCM}. We emphasize that we do not use \cite{FanCM} or Proposition~{\ref{prop:FanComparison}} anywhere else in this work.

\begin{proposition}[Comparison with \cite{FanCM}]\label{prop:FanComparison}
Let \(0<\kappa\le4\). For a Weil--Petersson map \(\phi\) let \(\mathcal A(\phi,\cdot)\) be the
relative action of \cite[Definition~6.4]{FanCM} (denoted \(I(\phi,\cdot)\) there), a
real random variable defined modulo \(P\)-null sets, with
\(\mathcal A(\phi,\cdot)\in L^p(P)\) for all \(p<\infty\) \cite[Theorem~1.9]{FanCM}. For
deterministic Weil--Petersson \(h\), \cite[Proposition~6.17]{FanCM} gives it a
pointwise value, \(\mathcal A(\phi,h)=H_\kappa(h\circ\phi^{-1})-H_\kappa(h)\)
\cite[(1.5)]{FanCM}; this applies in particular at \(h=\mathrm{id}\), although the
Weil--Petersson maps form a \(P\)-null set.
\begin{enumerate}
\item The law \(P\) is the law \(P^{\rm weld}\) of \cite[Definition~1.1, (6.3)]{FanCM},
and \(P^{\mathsf R}_\phi=(R_\phi)_*P^{\rm weld}\) in the notation of \cite[(6.7)]{FanCM}.
\item For every Weil--Petersson \(\phi\), \(dP^{\mathsf R}_\phi/dP=e^{-\mathcal A(\phi,\cdot)}\) \(P\)-a.s.
\cite[Theorem~1.9]{FanCM}. For every word \(w\),
\(B_w=\mathcal A(\Phi_w^{-1},\cdot)\) \(P\)-almost surely.
\item For every Weil--Petersson \(\phi\),
\[
 \operatorname{Ent}(P\,|\,P^{\mathsf R}_\phi)=\E_P[\mathcal A(\phi,\cdot)],\qquad
 \operatorname{Ent}(P^{\mathsf R}_\phi\,|\,P)=\E_P[\mathcal A(\phi^{-1},\cdot)].
\]
Consequently, given \cite[Theorem~1.9]{FanCM}, Theorem~\ref{thm:entropyWP} is
equivalent to the identity
\begin{equation}\label{eq:FanMean}
 \E_P[\mathcal A(\phi,\cdot)]=\mathcal A(\phi,\mathrm{id})=H_\kappa(\phi)
 \qquad\text{for every Weil--Petersson }\phi,
\end{equation}
with \(\mathcal A(\phi,\mathrm{id})\) in the pointwise sense above, and
Theorem~\ref{thm:entropyWP} proves \eqref{eq:FanMean}.
\item \emph{(Conditional on \cite[Theorem~1.4(ii)]{FanCM}.)} If
\(a\in\operatorname{Homeo}_+(S^1)\) is not Weil--Petersson, then
\(\operatorname{Ent}(P_a|P)=\operatorname{Ent}(P|P_a)=\infty\). Hence, with \(H_\kappa:=+\infty\) off the
Weil--Petersson class,
\(\operatorname{Ent}(P_a|P)=\operatorname{Ent}(P|P_a)=H_\kappa(a)\) for every \(a\in\operatorname{Homeo}_+(S^1)\).
\end{enumerate}
\end{proposition}

\begin{proof}
For this proof alone, we assume the results of \cite{FanCM}.

(i) By \cite[(6.2)--(6.3)]{FanCM}, \(P^{\rm weld}\) is the law of
\(e^{-i\alpha}\circ g^{-1}\circ f\), where the loop is drawn from the shape measure of
the SLE\(_\kappa\) loop measure \cite[Section~4.2, Definition~4.7]{FS} and \(\alpha\) is
independent and uniform. The SLE\(_\kappa\) loop measure is unique up to a
constant \cite[Theorem~1.9 and Section~1.1]{BJloop} (as recalled in
\cite[Section~4.2]{FS}), and the shape measure \(\mu_\kappa\), a probability measure, does
not depend on that constant; so this is
\cite[Definition~1.1]{BJweld}, that is, \(P\). In \cite[Definition~1.1]{FanCM}
the exterior map \(g_0\) is normalized by
\(f(1)=g_0(1)\) instead of \((g^+)'(\infty)>0\); then \(g_0=g^+\circ R_\beta\) with
\(\beta\) a function of the loop, and \(\mathrm{Rot}_Z\circ R_{-\beta}\) is again a uniform
rotation independent of the loop, so the two laws coincide. The maps
\(R_\phi(h)=h\circ\phi\) and the push-forward \(T_*Q=Q\circ T^{-1}\) are those of
\cite[(1.10), (6.7)]{FanCM}.

(ii) The first statement is \cite[Theorem~1.9]{FanCM}. By
Lemma~\ref{lem:Bfin}(i) and uniqueness of Radon--Nikodym derivatives,
\(B_w(y\circ\Phi_w^{-1})=-\mathcal A(\Phi_w,y)\) for \(P\)-a.e.\ \(y\). Since \(P^{\mathsf R}_{\Phi_w}\sim P\),
we may substitute \(y=h\circ\Phi_w\), and get \(B_w(h)=-\mathcal A(\Phi_w,h\circ\Phi_w)\) for
\(P\)-a.e.\ \(h\). The cocycle identity \cite[(6.66)]{FanCM} with \(\chi=\Phi_w^{-1}\),
together with \(\mathcal A(\mathrm{id},\cdot)=0\) (the density of the identity is \(1\)),
gives \(-\mathcal A(\Phi_w,h\circ\Phi_w)=\mathcal A(\Phi_w^{-1},h)\) \(P\)-a.s.

(iii) Since \(\mathcal A(\phi,\cdot)\in L^1(P)\),
\(\operatorname{Ent}(P|P^{\mathsf R}_\phi)=\int\log(dP/dP^{\mathsf R}_\phi)\,dP=\E_P[\mathcal A(\phi,\cdot)]\).
Relative entropy is invariant under the Borel bijection \(\mathsf R_{\phi^{-1}}\),
so \(\operatorname{Ent}(P^{\mathsf R}_\phi|P)=\operatorname{Ent}(P|P^{\mathsf R}_{\phi^{-1}})
=\E_P[\mathcal A(\phi^{-1},\cdot)]\). The pointwise value at \(h=\mathrm{id}\), together with
\(H_\kappa(\mathrm{id})=0\) and \(H_\kappa(\phi^{-1})=H_\kappa(\phi)\) (Section~\ref{sec:entropy}), gives
\(\mathcal A(\phi,\mathrm{id})=H_\kappa(\phi)\). Thus the identity
\(\operatorname{Ent}(P|P^{\mathsf R}_\phi)=H_\kappa(\phi)\) for all \(\phi\) is \eqref{eq:FanMean}, and it
implies \(\operatorname{Ent}(P^{\mathsf R}_\phi|P)=H_\kappa(\phi^{-1})=H_\kappa(\phi)\); the left action
follows as in the proof of Theorem~\ref{thm:entropyWP}.

(iv) We have \(P_a=(T_{a^{-1},\mathrm{id}})_*P\) in the notation of
\cite[(1.10)]{FanCM}, and \(a^{-1}\) is not Weil--Petersson, since the
Weil--Petersson maps form a group. By \cite[Theorem~1.4(ii)]{FanCM}, which
includes \(\kappa=4\) \cite[Theorem~1.8]{FanCM}, \(P_a\perp P\), so neither
measure is absolutely continuous with respect to the other, and both
entropies are infinite. The last statement then follows from
Theorem~\ref{thm:entropyWP}.
\end{proof}

\begin{remark}[What the comparison tests]\label{rem:FanTests}
Once both Theorem~\ref{thm:entropyWP} (through Lemma~\ref{lem:Bfin}(i)) and
\cite[Theorem~1.9]{FanCM} hold, part~(ii) is automatic, since two
Radon--Nikodym derivatives of the same pair of laws agree almost surely.
Proposition~\ref{prop:FanComparison}
verifies that the sign and orientation conventions of
Lemma~\ref{lem:Bfin} and Corollary~\ref{cor:entMobius} agree with those of
\cite{FanCM}. It is not an independent derivation: for finite-mode maps,
Fan's density is also obtained from finite Fourier flows
\cite[Lemma~6.9]{FanCM}, as the exponential of a time integral of the same
kind as \(B_w\), and both proofs rest on the same integration-by-parts input
(\cite[Theorem~4.1]{BJloop} there, \eqref{eq:W} here); Grunsky operators enter
only for general \(\phi\). Identity~\eqref{eq:FanMean} is not a consequence of
\cite{FanCM}: the pointwise formula \(\mathcal A(\phi,h)=H_\kappa(h\circ\phi^{-1})-H_\kappa(h)\) holds only
for Weil--Petersson \(h\), a \(P\)-null set, and the expectation of the renormalized
increment is the mean-value argument of Proposition~\ref{prop:entropyFin}.


We also note that this comparison does not independently confirm the value of \(\cL\): for us as for \cite{FanCM}, the value of \(\cL\) comes from the normalization in \cite{BJloop}.
The law \(P\) is defined from the
loop measure of \cite{BJloop}, whose restriction covariance carries the
exponent \(\frac{\cm}2\Lambda^*\) by \cite[Definition~1.8]{BJloop}; it is unique up
to a constant \cite[Theorem~1.9]{BJloop}, and Zhan's construction provides it
\cite[Section~1.1]{BJloop}. In \cite{FanCM} the constant enters only through
\cite[(6.4)]{FanCM}\(=\)\cite[Theorem~4.1]{BJloop}, with coefficient \(\cm/12\), for
this measure. Our proof of \eqref{eq:W} in Section~\ref{app:Wproof} starts instead
from Zhan's measure and the restriction covariance as displayed in
\cite[Theorem~5.1]{ZhanLoop}; the only open logical point is therefore that
Zhan's measure satisfies the covariance with exponent \(\frac{\cm}2\Lambda^*\), i.e.\ is
the measure of \cite{BJloop}, and this is the content of
Lemma~\ref{lem:normalization} (Remark~\ref{rem:Wp-normalization}(b)). For
\(\kappa<4\) the Gaussian multiplicative chaos comparison
\cite[Lemma~1.2]{FanCM} concerns only the measure class.
\end{remark}
\section{The Fisher metric of the welding family}
\label{sec:global}

Let \(0<\kappa\le4\). In this section, the Fisher information of the
finite-mode submodels of the welding family (a statistical quantity
computed from the finite-flow likelihoods of Section~\ref{sec:finiteinput})
is shown to extend uniquely to a continuous right-invariant metric on
the marked WP space. This metric is the right-invariant Kähler metric
\(2g_{\rm VK}+\frac{\cL}6\pi^*g_{\rm WP}\) of Takhtajan--Teo, which at the
identity is the Kirillov--Yuriev metric \(g_{\cL,1}\); by
Theorem~\ref{thm:entropyWP} its Kähler potential is the relative entropy,
and it is the Hessian of the relative entropy in every direction.

\phantomsection\label{sec:markedcomparison}%
We first recall the marked WP space, its complex structure and its chart. Let
\(\mathcal M=\operatorname{Rot}(S^1)\backslash\operatorname{WP}(S^1)\), rotations acting
by postcomposition; every class has a unique representative fixing \(1\),
which we always use. Takhtajan--Teo identify the universal Teichm\"uller
curve \(\mathcal T(1)\) with \(S^1\backslash\operatorname{Homeo}_{\rm qs}(S^1)\), the coset
space being identified with the stabilizer of \(1\)
\cite[Section~2.2.1, Remark~2.6]{TTI}; the projection
\(\pi:\mathcal T(1)\to T(1)=\operatorname{M\ddot{o}b}(S^1)\backslash\operatorname{Homeo}_{\rm qs}(S^1)\)
forgets the marking, with fibres \(S^1\backslash\operatorname{M\ddot{o}b}(S^1)\simeq\Ds\).
Thus \(\mathcal M\) is the Weil--Petersson Teichm\"uller curve
\(\mathcal T_0(1)=\pi^{-1}(T_0(1))\) of \cite[Section~2]{TT}, fibred over the
Weil--Petersson Teichm\"uller space \(T_0(1)\) with fibre \(\Ds\); it is a
topological group under composition \cite[Appendix~A]{TTI}. For
\(h\in\mathcal M\) write
\[
 h=g_h^{-1}\circ f_h,\qquad f_h(0)=0,\quad f_h'(0)=1,
 \quad g_h(\infty)=\infty
\]
(properties CW2\('\)--CW3\('\) of \cite[Section~2.2.1]{TTI}), and
\(q(h)=\log f_h'\). The complex Hilbert manifold structure of
\(\mathcal T_0(1)\) is modelled on \(A_2(\D)\oplus\mathbb C\) through the Bers
embedding \(h\mapsto(\mathcal Sf_h,\frac12f_h''(0))\), and by
\cite[Appendix~A, Theorem~A.3]{TT} the pre-Bers embedding
\(h\mapsto f_h''/f_h'\) into
\(A^1_2(\D)=\{\psi\text{ holomorphic on }\D:\int_\D|\psi|^2d^2z<\infty\}\)
induces the same structure (a result parallel to \cite[Appendix]{Teo}). Since \(q\mapsto q'\) maps
\[
 AD_0(\mathbb D)=\left\{q\text{ holomorphic}:q(0)=0,
       \ \|q\|^2=\frac1\pi\int_{\mathbb D}|q'|^2dA<\infty\right\},
 \qquad \Bigl\|\sum_{n\ge1}q_nz^n\Bigr\|^2=\sum_{n\ge1}n|q_n|^2,
\]
isometrically, up to the factor \(\sqrt\pi\), onto \(A^1_2(\D)\), the map \(q\)
is a holomorphic chart of \(\mathcal M\) onto an open subset of \(AD_0(\D)\),
with \(q(\mathrm{id})=0\). This is the \emph{standard} (Takhtajan--Teo) complex
structure of \(\mathcal M\). The projection \(\pi\) is holomorphic, being the
first Bers coordinate followed by the Bers embedding of \(T_0(1)\),
which is biholomorphic onto its image \cite[Theorem~3.13]{TTI}. The
inclusion \(\iota:\mathcal M\hookrightarrow\mathcal T(1)\) into the Banach curve is
holomorphic: in Bers coordinates (the Bers fibre-space model of
\cite[Section~2.2.1]{TTI}, \cite[Appendix~A]{TT}) it is the restriction of
the bounded inclusion \(A_2(\D)\oplus\mathbb C\hookrightarrow A_\infty(\D)\oplus\mathbb C\)
\cite[Lemma~3.1]{TTI}. Statements proved in \cite{TTI} on the Banach curve
\(\mathcal T(1)\) therefore restrict to \(\mathcal M\).

Next we describe the tangent space at the identity in this chart. By \cite[TS1, Section~2.2.3]{TTI},
if \(\gamma_s=g_s^{-1}\circ f_s\) is a smooth curve tangent at \(s=0\) to the real
field \(u=\sum_{n\ne0}c_ne^{in\theta}\frac d{d\theta}\), then
\(\dot f=u_+=i\sum_{n\ge1}c_nz^{n+1}\), so that in the chart \(q\) its velocity is
\(u_+'=i\sum_{n\ge1}(n+1)c_nz^n\). The properties TS1--TS2 are formulated in
\cite{TTI} for the Banach (Zygmund) structure; in the Hilbert chart \(q\)
the vector \(u_+'\) has squared norm \(\sum_{n\ge1}n(n+1)^2|c_n|^2\), which is
finite exactly when \(u\in H^{3/2}\). We therefore identify
\(T_{\mathrm{id}}\mathcal M\) with \(H^{3/2}(S^1;\mathbb R)/\mathbb R\) through the
topological isomorphism
\begin{equation}\label{eq:tangentId}
 Dq_{\mathrm{id}}(v)=v_+',\qquad v_+=i\sum_{n\ge1}c_nz^{n+1}
 \quad\text{for}\quad v=c_0+\sum_{n\ne0}c_ne^{in\theta},
\end{equation}
the constants being the forgotten rotation direction. On real
trigonometric polynomials it is the velocity map of TS1, and these form a
dense subspace. We write \(c_n=it_n\), i.e.\
\(v=c_0+\sum_{n\ge1}(it_ne^{in\theta}-i\bar t_ne^{-in\theta})\). The complex
structure \(J\) of the chart is multiplication of \((t_n)_{n\ge1}\) by \(i\),
that is, the conjugation operator \(Ju=i\sum\operatorname{sgn}(n)c_ne^{in\theta}\) of
\cite[Section~2.2.3]{TTI}, and \(v^{1,0}=(v-iJv)/2=\sum_{n\ge1}c_ne^{in\theta}\).
For velocities of flows we use Shen--Tang \cite[Theorem~1.2]{ShenTang}: for a
continuous \(H^{3/2}\) vector field \(v\) on \(S^1\) the flow \(s\mapsto\psi^v_s\) lies in
the Weil--Petersson class and is continuously differentiable for the
Takhtajan--Teo Hilbert manifold structure, with
\(\frac d{ds}\psi^v_s=v\circ\psi^v_s\) (see also \cite[Theorem~7.3]{Shen} for the
continuous differentiability of \(s\mapsto\log(\psi^v_s)'\) in \(H^{1/2}\)). On the
rotation quotient \(\mathcal M\) this is immediate from the chart
\(\mathfrak l\) below: \([\psi^v_s]=\mathfrak l^{-1}([\log(\psi^v_s)'])\), where
\(\mathfrak l^{-1}\) is real-analytic \cite[Theorem~1.2]{WHS} and
\(s\mapsto\log(\psi^v_s)'\) is \(C^1\) into \(H^{1/2}\) \cite[Theorem~7.3]{Shen}
(compare \cite[proof of Theorem~1.3]{WHS}). We use this only for real
trigonometric polynomials \(v\): the velocity of \(s\mapsto[\psi^v_s]\) at
\(s=0\) is the class of \(v\).

The group structure enters through right translations. For \(b\in\mathcal M\) put
\(\mathsf R_b(k)=k\circ b\); with representatives fixing \(1\) this is well
defined on classes.

\begin{lemma}[Right translations]\label{lem:rightTranslations}
For every \(b\in\mathcal M\), \(\mathsf R_b\) is a biholomorphic automorphism of the
Hilbert manifold \(\mathcal M\), with inverse \(\mathsf R_{b^{-1}}\), and
\(\pi\circ\mathsf R_b=\mathsf R_{\pi(b)}\circ\pi\).
\end{lemma}
\begin{proof}
Since \(\mathcal M\) is a group \cite[Appendix~A]{TTI}, \(\mathsf R_b\) is a bijection of
\(\mathcal M\) with inverse \(\mathsf R_{b^{-1}}\); the equivariance of \(\pi\) is immediate
from the definitions. It suffices to show that both Bers coordinates of
\(\mathsf R_b(k)\) are holomorphic in \(k\). The Schwarzian is invariant under
M\"obius postcomposition, so the first coordinate is the Bers embedding of
\(\mathsf R_{\pi(b)}(\pi(k))\); right translations of the Hilbert manifold \(T(1)\)
are biholomorphic \cite[Corollary~3.12]{TTI} and preserve its component
\(T_0(1)\), \(\pi(b)\in T_0(1)\) \cite[Section~3.4]{TTI}, and the Bers embedding
is biholomorphic \cite[Theorem~3.13]{TTI}. The second coordinate
\(k\mapsto\frac12f''_{k\circ b}(0)\) is the composition of the holomorphic
inclusion \(\iota\), the right translation of the Banach curve \(\mathcal T(1)\),
which is biholomorphic \cite[Section~3, Theorem~2]{Bers73} (see \cite[Section~2.2.2]{TTI}), and
the second Bers coordinate of \(\mathcal T(1)\); hence it is holomorphic.
\end{proof}

Hence every tangent vector at \(h\) is \(X=D\mathsf R_h(v)\) for a unique
\(v\in T_{\mathrm{id}}\mathcal M\); we write \(v=X\circ h^{-1}:=D\mathsf R_{h^{-1}}(X)\), which
pointwise is the angular field \(X\circ h^{-1}\).

We now recall the two Kähler potentials. Let \(K=\log|g'(\infty)/f'(0)|\), which is
\(\log|g_h'(\infty)|\) on \(\mathcal M\), and let \(S_1\) be the universal Liouville
action on \(T_0(1)\) (Section~\ref{sec:entropy}, Proposition~\ref{prop:potential}).
The function \(K\) is real-analytic on \(\mathcal T(1)\) and \(K\circ\mathsf R_b-K\) is
pluriharmonic there \cite[Lemma~5.1; Lemma~5.4 and its proof]{TTI}; composing
with \(\iota\), the same holds on \(\mathcal M\). The function \(S_1\) is a Kähler
potential of the Weil--Petersson metric on \(T_0(1)\),
\(i\partial\bar\partial S_1=2\omega^{\rm TT}_{\rm WP}\)
\cite[Theorems~4.8 and~5.1, Corollary~5.2]{TT}; since the Weil--Petersson
metric is right-invariant \cite[Section~4.2]{TTI} and right translations are
biholomorphic, \(S_1\circ\mathsf R_{\pi(b)}-S_1\) is pluriharmonic. Moreover
\(S_1\) is real-analytic on \(T_0(1)\): by \cite[Theorem~4.8]{TT},
\(S_1=-12\pi S_2\), and \(S_2=\log\operatorname{Det}(I-Z\bar Z)\) is the pullback by
the restricted period map \(Z\) of the Kähler potential
\(\log\operatorname{Det}(I-Z\bar Z)\) of the Hilbert--Schmidt Siegel disc
\(\mathfrak D^0_\infty\), the period map being a holomorphic map of Hilbert
manifolds into \(\mathfrak D^0_\infty\) \cite[Theorem~3.10, Corollary~6.4, Remark~6.5]{TT}. Of the potential property we use only
that \(S_1\) is of class \(C^2\) on \(T_0(1)\), with this mixed Hessian. Consequently
\[
 V=2K+\frac{\cL}{24\pi}\,\pi^*S_1
\]
is a \(C^2\) function on \(\mathcal M\) and \(V\circ\mathsf R_b-V\) is pluriharmonic for
every \(b\in\mathcal M\).

A second, real coordinate on \(\mathcal M\) is the logarithmic derivative. The map
\(\mathfrak l([h])=[\log h']\in H^{1/2}(S^1;\mathbb R)/\mathbb R\) is well defined on \(\mathcal M\)
(postcomposition by rotations preserves the angular derivative, and the
total increment \(2\pi\) fixes the constant) and is a homeomorphism for the
Hilbert topology \cite[Theorem~8.1]{Shen}; Shen and Shen--Tang write the
quotient as \(\operatorname{WP}(S^1)/\operatorname{Rot}(S^1)\), and since \(\log h'\) is invariant only
under postcomposition by rotations, it is the quotient \(\mathcal M\). The transport
identity
\begin{equation}\label{eq:WHStransport}
 [\log(h\circ h_0^{-1})']
   =[(\log h'-\log h_0')\circ h_0^{-1}]
\end{equation}
shows that right translations are affine in \(\mathfrak l\), and bounded, since
composition with a fixed quasisymmetric homeomorphism is bounded on
homogeneous \(H^{1/2}\) \cite[Theorem~3.1]{NagSullivan}. Moreover \(\mathfrak l\) is a
real-analytic diffeomorphism of \(\mathcal M\), with its Takhtajan--Teo
structure, onto \(H^{1/2}(S^1;\mathbb R)/\mathbb R\) \cite[Theorem~1.2]{WHS},
which settles \cite[Conjecture~9.1]{Shen} (the analogue on the real line is
\cite[Theorem~2.3]{ShenTang}). We use \(\mathfrak l\) only for the flows above.

Finally we introduce the statistical objects. The welding family is
\(E_h=P_{h^{-1}}\), \(h\in\mathcal M\). It is well defined
on the rotation quotient, since \(P_{(R\circ h)^{-1}}=(L_{h^{-1}})_*(L_{R^{-1}})_*P=E_h\)
by left-rotation invariance of \(P\) (Lemma~\ref{lem:symmetries}). Inversion
only relabels the laws and is never differentiated.
For a real angular field
\(v=c_0+\sum_{n\ge1}(it_ne^{in\theta}-i\bar t_ne^{-in\theta})\) put
\begin{equation}\label{eq:Fishernorm}
 \|v\|_F^2=2\sum_{n\ge1}C_n|t_n|^2;
\end{equation}
modulo constants this norm is equivalent to the norm of \eqref{eq:tangentId},
hence to \(H^{3/2}\), and \(\|Jv\|_F=\|v\|_F\). For a finite
Fourier field, \(\mathsf S(v)=\sum_n(t_ns_n+\bar t_n\bar s_n)\) is the score
\(b_{\mathbf t}\) of Theorem~\ref{thm:finite}: the angular field of the
generator of \(\Phi_{\mathbf t}\) is \(v_{\mathbf t}=\sum_{n\le N}(it_ne^{in\theta}-i\bar t_ne^{-in\theta})\)
(proof of Theorem~\ref{thm:finite}), and \eqref{eq:allcov} gives
\(\E[\mathsf S(v)\mathsf S(w)]=2\Rea\sum_nC_nt_n\bar q_n\) for coefficients
\(t_n,q_n\) of \(v,w\), in particular \(\E\mathsf S(v)^2=\|v\|_F^2\).

\phantomsection\label{sec:markedKahler}%
We can now state the main result of this section, the identification of
the Fisher metric and of its Kähler potential.
For \(h\in\mathcal M\), the \emph{finite-mode tangent vectors} at \(h\) are the vectors
\(X_{\mathbf t}=D\mathsf R_h(v_{\mathbf t})\), \(N<\infty\), \(\mathbf t\in\mathbb C^N\); they form a
dense subspace \(\mathcal V_h\subset T_h\mathcal M\), since \(D\mathsf R_h\) is a topological
isomorphism and the finite Fourier fields are dense in \(T_{\mathrm{id}}\mathcal M\). For a
real potential \(A\) of class \(C^2\), put
\(\mathfrak g_A(X,Y)=2\Rea(\partial\bar\partial A)
(X^{1,0},\overline{Y^{1,0}})\), \(X^{1,0}=(X-iJX)/2\), where
\(\partial\bar\partial A\) is the Hermitian coefficient tensor of the
real form \(i\partial\bar\partial A\); thus \(i\partial\bar\partial A=\mathfrak g_A(J\cdot,\cdot)\). In a
holomorphic chart the real Hessian is
\begin{equation}\label{eq:realHessian}
 D^2A(X,X)=\mathfrak g_A(X,X)+2\Rea(\partial\partial A)(X^{1,0},X^{1,0}),
\end{equation}
and at a critical point of \(A\) both terms are independent of the
holomorphic chart. Let \(\pi^*\) denote pull-back by
\(\pi:\mathcal M\to T_0(1)\), \(E=\pi^*S_1/\pi\), \(g_{\rm VK}=\mathfrak g_K\) and
\(\pi^*g_{\rm WP}=\mathfrak g_E/4\). At the identity,
\(g_{\rm VK}(v,w)=2\Rea\sum_{n\ge1}n t_n\bar q_n\) and
\((\pi^*g_{\rm WP})(v,w)=\Rea\sum_{n\ge2}(n^3-n)t_n\bar q_n\) on finite Fourier
fields (Appendix~\ref{sec:round}), hence on all of \(T_{\mathrm{id}}\mathcal M\) by
continuity; the WP pullback vanishes on the forgotten directions. On
smooth marked weldings the potential \(V\) is the Liouville--capacity action
\(H_\kappa\), the formal log-density of \cite[(1.10), (1.14)]{BJweld}.

\begin{theorem}[Fisher metric of the welding family]
\label{thm:globalKahler}
Let \(0<\kappa\le4\), give \(\mathcal M\) its standard
Takhtajan--Teo complex Hilbert structure, and let
\(V=2K+\cL\pi^*S_1/(24\pi)\) be as defined above.
\begin{enumerate}
\item \emph{Entropy.} For all \(h,k\in\mathcal M\) the laws \(E_h\), \(E_k\)
are equivalent, and
\[
 \operatorname{Ent}(E_h\,|\,P)=\operatorname{Ent}(P\,|\,E_h)=V(h)=H_\kappa(h),
\]
\[
 \operatorname{Ent}(E_k\,|\,E_h)=\operatorname{Ent}(E_h\,|\,E_k)=H_\kappa(k\circ h^{-1})
 =V(\mathsf R_{h^{-1}}k).
\]
\item \emph{Finite-mode Fisher information.} Let \(h\in\mathcal M\),
\(N<\infty\), and let \(\Phi_{\mathbf t}\), \(b_{\mathbf t}\) and
\(v_{\mathbf t}\) be as in Theorem~\ref{thm:finite}, \(\mathbf t\in\mathbb C^N\).
The laws of the submodel \(\mathbf t\mapsto E_{\Phi_{\mathbf t}\circ h}\) are
equivalent to \(E_h\) with densities bounded above and below, and jointly as
\(\mathbf t\to0\),
\[
 \sqrt{\frac{dE_{\Phi_{\mathbf t}\circ h}}{dE_h}}
 =1-\frac12\,b_{\mathbf t}\circ L_h+O_{L^\infty(E_h)}(|\mathbf t|^2),
 \qquad
 \E_{E_h}\bigl[(b_{\mathbf t}\circ L_h)^2\bigr]=\|v_{\mathbf t}\|_F^2,
\]
where \(X_{\mathbf t}=D\mathsf R_h(v_{\mathbf t})\) is the velocity at \(s=0\) of
the path \(s\mapsto\Phi_{s\mathbf t}\circ h\). Moreover
\[
 \operatorname{Ent}(E_{\Phi_{\mathbf t}\circ h}\,|\,E_h)
 =H_\kappa(\Phi_{\mathbf t})
 =\tfrac12\|v_{\mathbf t}\|_F^2+O(|\mathbf t|^3).
\]
\item \emph{Extension and Kähler identification.} For every \(h\in\mathcal M\)
there is a unique continuous quadratic form \(I_h\) on \(T_h\mathcal M\) such that
\(I_h(X_{\mathbf t},X_{\mathbf t})\) is the Fisher information
\(\E_{E_h}[(b_{\mathbf t}\circ L_h)^2]\) of (ii) for every finite-mode tangent
vector \(X_{\mathbf t}\); it is \(I_h(X,X)=\|X\circ h^{-1}\|_F^2\). The metric
\(I=(I_h)_h\) is right-invariant and positive definite, \(I_h\) is equivalent to the
\(H^{3/2}\) norm on \(T_h\mathcal M\), and for real \(X,Y\in T_h\mathcal M\), with
\(v=X\circ h^{-1}\), \(w=Y\circ h^{-1}\) having coefficients \(t_n,q_n\),
\[
 I_h(X,Y)=\mathfrak g_V(X,Y)=2\Rea\bigl(\partial\bar\partial V\bigr)_h
                  (X^{1,0},\overline{Y^{1,0}})
 =2\Rea\sum_{n\ge1}C_nt_n\bar q_n .
\]
Thus \(I=2g_{\rm VK}+\frac{\cL}{6}\pi^*g_{\rm WP}\) is a Kähler metric whose Kähler
potential is the relative entropy \(h\mapsto\operatorname{Ent}(E_h|P)=V(h)\).
\item \emph{Entropy Hessian in all directions.} For every \(h\in\mathcal M\) the
function \(F_h(k)=\operatorname{Ent}(E_k|E_h)\) is \(C^2\) on \(\mathcal M\), with \(F_h(h)=0\) and
\(dF_h(h)=0\), and its real Hessian at \(h\) is \(D^2F_h(h)=I_h\). Equivalently,
for every \(X\in T_h\mathcal M\) and every \(C^2\) path \(s\mapsto k_s\) in \(\mathcal M\) with
\(k_0=h\) and \(\dot k_0=X\),
\[
 \operatorname{Ent}(E_{k_s}\,|\,E_h)=\operatorname{Ent}(E_h\,|\,E_{k_s})
 =\frac{s^2}2\,I_h(X,X)+o(s^2)\qquad(s\to0).
\]
\end{enumerate}
\end{theorem}

In (iv) the second derivative at \(s=0\) does not depend on \(\ddot k_0\), because
\(dF_h(h)=0\). Inversion \(h\mapsto h^{-1}\) is not asserted to be holomorphic or
\(C^1\); it enters only through the labels \(E_h=P_{h^{-1}}\).

\begin{proof}
(i) By Theorem~\ref{thm:entropyWP} with \(a=h^{-1}\), which is WP, and
\(H_\kappa(h^{-1})=H_\kappa(h)\) (Section~\ref{sec:entropy}),
\(\operatorname{Ent}(E_h|P)=\operatorname{Ent}(P|E_h)=H_\kappa(h)\).
Corollary~\ref{cor:entropyConsequences}(i) with \(a=h^{-1}\), \(b=k^{-1}\) gives
\(\operatorname{Ent}(E_h|E_k)=H_\kappa(k\circ h^{-1})\), and exchanging \(h\) and \(k\),
\(\operatorname{Ent}(E_k|E_h)=H_\kappa(h\circ k^{-1})=H_\kappa(k\circ h^{-1})\) by inversion
symmetry of \(H_\kappa\). Finiteness in both directions gives equivalence. On
\(\mathcal M\), \(H_\kappa=V\) by the definitions of \(K\) and \(S_1\)
(Section~\ref{sec:entropy}, Proposition~\ref{prop:potential}), and
\(k\circ h^{-1}=\mathsf R_{h^{-1}}(k)\); \(H_\kappa\) is invariant under rotations on both
sides, so the value does not depend on the representatives.

(ii) Since \(\Phi_{\mathbf t}^{-1}=\Phi_{-\mathbf t}\) (the generator is real-linear
in \(\mathbf t\)), \(E_{\Phi_{\mathbf t}\circ h}=P_{h^{-1}\circ\Phi_{-\mathbf t}}
=(L_{h^{-1}})_*P_{\Phi_{-\mathbf t}}\) and \(E_h=(L_{h^{-1}})_*P\). Pushing forward
by the fixed Borel bijection \(L_{h^{-1}}\), whose inverse is \(L_h\),
\(dE_{\Phi_{\mathbf t}\circ h}/dE_h=r_{-\mathbf t}\circ L_h\), with \(r\) from
\eqref{eq:finiteRN}. By the proof of Theorem~\ref{thm:finite},
\(\|b_{\mathbf t}\|_\infty\le M_N|\mathbf t|\) and
\(\|\log r_{\mathbf t}-b_{\mathbf t}\|_\infty\le L_N|\mathbf t|^2\), so the density is bounded
above and below. Since \(b_{-\mathbf t}=-b_{\mathbf t}\),
\eqref{eq:finiteDQM} gives the expansion, and \((L_h)_*E_h=P\) together
with \eqref{eq:allcov} gives
\(\E_{E_h}[(b_{\mathbf t}\circ L_h)^2]=\E b_{\mathbf t}^2=2\sum_nC_n|t_n|^2=\|v_{\mathbf t}\|_F^2\).
The path \(s\mapsto\Phi_{s\mathbf t}\) is the flow of the trigonometric field
\(v_{\mathbf t}\), so its velocity at \(s=0\) is \(v_{\mathbf t}\)
(recalled at the beginning of this section), and by Lemma~\ref{lem:rightTranslations}
that of \(s\mapsto\mathsf R_h(\Phi_{s\mathbf t})\) is \(X_{\mathbf t}\). By (i),
\(\operatorname{Ent}(E_{\Phi_{\mathbf t}\circ h}|E_h)=H_\kappa(\Phi_{\mathbf t})\).
For the expansion, put \(Q=E_h\) and
\(\ell=\log(dE_{\Phi_{\mathbf t}\circ h}/dE_h)=\log r_{-\mathbf t}\circ L_h\). The two bounds
give \(\|\ell+b_{\mathbf t}\circ L_h\|_{L^\infty(Q)}\le L_N|\mathbf t|^2\) and
\(\|\ell\|_{L^\infty(Q)}=O(|\mathbf t|)\). From \(\E_Q e^\ell=1\) and
\(\operatorname{Ent}(E_{\Phi_{\mathbf t}\circ h}|E_h)=\E_Q[\ell e^\ell]\), Taylor expansion of the
exponential gives
\(0=\E_Q[\ell+\frac12\ell^2]+O(|\mathbf t|^3)\) and
\(\E_Q[\ell e^\ell]=\E_Q[\ell+\ell^2]+O(|\mathbf t|^3)\), hence
\(\operatorname{Ent}(E_{\Phi_{\mathbf t}\circ h}|E_h)=\frac12\E_Q\ell^2+O(|\mathbf t|^3)
=\frac12\E_Q(b_{\mathbf t}\circ L_h)^2+O(|\mathbf t|^3)=\frac12\|v_{\mathbf t}\|_F^2+O(|\mathbf t|^3)\).

(iii) \emph{Uniqueness and existence.} By (ii) the finite-mode Fisher
informations define the quadratic form \(X_{\mathbf t}\mapsto\|v_{\mathbf t}\|_F^2\) on
\(\mathcal V_h\). The form \(X\mapsto\|D\mathsf R_{h^{-1}}X\|_F^2\) is continuous on \(T_h\mathcal M\)
(Lemma~\ref{lem:rightTranslations} and \eqref{eq:Fishernorm}) and extends it.
Two continuous quadratic forms that agree on the dense subspace \(\mathcal V_h\)
agree, since their polarizations are continuous bilinear forms. Positivity
and equivalence with \(H^{3/2}\) follow from \(C_n>0\), constants being the
forgotten rotation direction, and from \(D\mathsf R_h\) being a topological
isomorphism. \emph{Right invariance.} Since
\(\mathsf R_{(k\circ b)^{-1}}\circ\mathsf R_b=\mathsf R_{k^{-1}}\),
\(I_{k\circ b}(D\mathsf R_bX,D\mathsf R_bX)=\|D\mathsf R_{k^{-1}}X\|_F^2=I_k(X,X)\).
\emph{Kähler identification.} By (i), \(V\ge0=V(\mathrm{id})\), so
\(dV(\mathrm{id})=0\), and \(D^2V(\mathrm{id})\) is a continuous symmetric bilinear form on
\(T_{\mathrm{id}}\mathcal M\). The path \(f_{\mathbf t}\) in the proof of
Proposition~\ref{prop:potential} is holomorphic in the chart \(q\), with first
jet \(Dq_{\mathrm{id}}(v_{\mathbf t})\) by \eqref{eq:tangentId}, and only the first jet
matters at a critical point; so Proposition~\ref{prop:potential} gives
\(D^2V(\mathrm{id})(v,v)=\|v\|_F^2\) for finite Fourier fields \(v\), hence, by
continuity and density, for all \(v\in T_{\mathrm{id}}\mathcal M\). Since \(\|Jv\|_F=\|v\|_F\) and
\((Jv)^{1,0}=iv^{1,0}\), \eqref{eq:realHessian} gives
\(\mathfrak g_V(v,v)=\frac12[D^2V(\mathrm{id})(v,v)+D^2V(\mathrm{id})(Jv,Jv)]=\|v\|_F^2\) at the
identity. As \(V\circ\mathsf R_b-V\) is pluriharmonic and \(\mathsf R_b\) is biholomorphic,
\(\mathsf R_b^*\mathfrak g_V=\mathfrak g_{V\circ\mathsf R_b}=\mathfrak g_V\); with \(b=h^{-1}\),
\(\mathfrak g_V(X,X)=\mathfrak g_V(\mathrm{id})(X\circ h^{-1},X\circ h^{-1})=I_h(X,X)\). Polarization
gives the Fourier formula. Since \(A\mapsto\mathfrak g_A\) is linear, the definitions of \(g_{\rm VK}\) and
\(\pi^*g_{\rm WP}\) give \(\mathfrak g_V=2g_{\rm VK}+\frac{\cL}6\pi^*g_{\rm WP}\). The form
\(\omega_I=I(J\cdot,\cdot)=i\partial\bar\partial V\) is closed, so \(I\) is Kähler.

(iv) By (i), \(F_h=V\circ\mathsf R_{h^{-1}}\), which is \(C^2\) by
Lemma~\ref{lem:rightTranslations}; \(F_h\ge0=F_h(h)\), so \(dF_h(h)=0\). The chain rule at
the critical point gives
\(D^2F_h(h)(X,X)=D^2V(\mathrm{id})(D\mathsf R_{h^{-1}}X,D\mathsf R_{h^{-1}}X)=\|X\circ h^{-1}\|_F^2=I_h(X,X)\),
by the computation of \(D^2V(\mathrm{id})\) in (iii). For a \(C^2\) path,
\(\frac{d^2}{ds^2}F_h(k_s)|_{s=0}=D^2F_h(h)(X,X)+dF_h(h)(\ddot k_0)=I_h(X,X)\), and
Taylor's formula gives the expansion; \(\operatorname{Ent}(E_h|E_{k_s})=F_h(k_s)\) by (i).
\end{proof}

We now separate what is known in Theorem~\ref{thm:globalKahler} from what
is new. The metric
\(2g_{\rm VK}+\frac{\cL}6\pi^*g_{\rm WP}\), its right invariance and the potentials
\(K\) and \(S_1\) are due to Takhtajan--Teo \cite[Theorem~5.3, Lemma~5.4]{TTI},
\cite[Theorem~5.1, Corollary~5.2]{TT}, and its value \(g_{\cL,1}\) at the
identity is the Kirillov--Yuriev normalization \cite[Remark~5.6]{TTI}
(Lemma~\ref{lem:conventions}). The Fisher form of the finite-mode
submodels at the identity is in \cite[(1.24), (1.36)]{CaiResponse}, where
(1.36) expresses the covariance as the averaged second derivative of the
negative log density along finite Fourier flows. Cai also proves Kähler
compatibility at a general Weil--Petersson base point
\cite[Section~2.3, Propositions~2.10--2.11]{CaiResponse}; this is a
tangent-level statement in period-mapping coordinates, in which the
central charge enters only through the covariance at the identity, and
it involves neither a Kähler potential nor the relative entropy. These are
infinitesimal statements, about the covariance of the score; the
relative entropy is a global functional of the laws, which the metric
determines only to second order at each point (see the discussion after
Theorem~\ref{thm:Akahler}). Since the
metric of part (ii) is right invariant, its value at a general base point is
determined by its value at the identity; the content of (ii) is that the
finite-mode Fisher informations of the laws \(E_h\), computed directly at
\(h\), are those of this right-invariant metric. That the identity value is a
Kirillov--Yuriev form reflects the Virasoro structure of the score
covariance: \(\theta_n=2hn+\frac c{12}(n^3-n)\) is the Shapovalov norm of
\(L_{-n}v\) for a highest-weight vector \(v\) of weight \(h\), cf.\
\cite{AiraultMalliavin,GQW,BJloop}. The statements proved here are the identification of the potential with the
relative entropy (by Theorem~\ref{thm:entropyWP}), the finite-mode Fisher
informations through every base point, their unique continuous
extension, and the entropy Hessian in all directions, part (iv), which
gives the metric a statistical meaning beyond finite-mode submodels.

\begin{remark}[Fisher information along general WP paths]
\label{rem:criticalPaths}
Theorem~\ref{thm:globalKahler}(ii) concerns finite-mode flows through an
arbitrary base point, and (iv) the relative entropy along all \(C^2\) paths.
That \(I_h(X,X)\) is the Fisher information of
\(t\mapsto E_{h_t}\) along every standard \(C^1\) path with \(\dot h_0=X\), i.e.\
differentiability in quadratic mean of the welding family at the critical
\(H^{3/2}\) regularity, and whether the Hellinger length of the path of
laws driven by an integrable \(H^{3/2}\) control \(v\) equals
\(\frac12\int\|v_t\|_F\,dt\), are not addressed here. Part (iv) gives the
corresponding second-order statement for the relative entropy along all
\(C^2\) paths.
\end{remark}

\phantomsection\label{sec:KE}%
The rest of this section decomposes the Fisher metric according to
central charges and compares it with the Kirillov--Yuriev metrics. For
this comparison the sign conventions for the complex structure matter.
We use the complex
structure \(J\) of Takhtajan--Teo (recalled at the beginning of this section), which is
the negative of the one of Kirillov and Nag--Verjovsky \cite[Remark~2.8]{TTI};
Takhtajan--Teo also use right cosets where Kirillov uses left cosets
\cite[Remark~2.7]{TTI}, and accordingly we write
\(S^1\backslash\operatorname{Diff}_+(S^1)\) for Kirillov's homogeneous space. The Kähler form of
a metric \(g\) is \(\omega=g(J\cdot,\cdot)\), and Ricci forms are
\(\mathrm{Ric}(J\cdot,\cdot)\), both with this \(J\). Riemannian metrics and Ricci
tensors do not depend on the sign of \(J\); Kähler and Ricci forms change
sign with it, so identities between them, such as \(\mathrm{Ric}^{\rm reg}=-\omega_{26,1}\)
below, hold in either convention, while each form taken alone is the
negative of Kirillov's.

Let \(g_{c,h}\) be the Kirillov--Yuriev metrics on
\(S^1\backslash\operatorname{Diff}_+(S^1)\) \cite{KirillovYuriev}, with Kähler forms
\(\omega_{c,h}=g_{c,h}(J\cdot,\cdot)\), normalized as in
\cite[Proposition 3.5]{GordinaLescot} by
\(g_{c,h}(u_n,u_n)=\theta_n/2\), \(\theta_n=2hn+\frac c{12}(n^3-n)\),
where \(u_n=\cos n\theta\,\partial_\theta\). Takhtajan--Teo use
\(\langle\mu,\nu\rangle_{\rm WP}=\iint_{\D^*}\mu\bar\nu\,\rho\,d^2z\),
\(\rho=4(1-|z|^2)^{-2}\), \(g_{k\bar l}=\langle\mu_k,\mu_l\rangle_{\rm WP}\)
in Bers coordinates,
\(\omega^{\rm TT}_{\rm WP}=\frac i2\sum g_{k\bar l}\,d\varepsilon_k\wedge
d\bar\varepsilon_l\) (Riemannian metric \(g^{\rm TT}_{\rm WP}\)) and
\(\mathrm{Ric}^{\rm TT}_{\rm WP}=\frac i2\sum R_{k\bar l}\,d\varepsilon_k
\wedge d\bar\varepsilon_l\), with \(R_{k\bar l}=\sum_nR_{k\bar nn\bar l}\)
summed over an orthonormal basis of \(H^{-1,1}(\D^*)\).

\begin{lemma}[Convention dictionary]\label{lem:conventions}
At the identity, for \(n\ge1\) (the WP terms vanish at \(n=1\)); all
tensors below are diagonal in \(n\) by rotation invariance, and right
invariance refers to right composition on \(\operatorname{Rot}(S^1)\backslash\operatorname{WP}(S^1)\)
as above:
\begin{enumerate}
\item \(g_{\rm VK}(u_n,u_n)=n/2\), \(\pi^*g_{\rm WP}(u_n,u_n)=(n^3-n)/4\),
\(g_{c,h}=\mathfrak g_{2hK+cS_1/(24\pi)}\), and
\(I(u_n,u_n)=C_n/2\). Hence \(I\) is the right-invariant extension of
\(g_{\cL,1}\), and at the identity it coincides with \(g_{\cL,1}\).
\item \(g^{\rm TT}_{\rm WP}=2\pi g_{\rm WP}=\frac12\mathfrak g_{S_1}\)
(\(i\partial\bar\partial S_1=2\omega^{\rm TT}_{\rm WP}\)), with value
\(\frac\pi2(n^3-n)\) on \((u_n,u_n)\).
\item \(R_{k\bar l}=-\frac{13}{12\pi}g_{k\bar l}\)
\cite[Theorem~7.11, Corollary~7.12]{TTI} (so \(T(1)\), \(T_0(1)\) are
Kähler--Einstein); the Riemannian Ricci tensor is
\(\mathrm{Ric}(g_{\rm WP})=-\frac{13}{6\pi}g^{\rm TT}_{\rm WP}\), with
Ricci form \(2\mathrm{Ric}^{\rm TT}_{\rm WP}=-\omega_{26,0}\) and
\(\pi^*\mathrm{Ric}(g_{\rm WP})(u_n,u_n)=-\frac{13}{12}(n^3-n)
=-g_{26,0}(u_n,u_n)\).
\item The regularized Ricci tensor of every \(g_{c,h}\) on
\(S^1\backslash\operatorname{Diff}_+(S^1)\) is \(-g_{26,1}\) (value
\(-(13n^3-n)/12\) on \((u_n,u_n)\); Ricci form \(-\omega_{26,1}\))
\cite{BowickRajeev}, \cite[Theorem 4.11]{GordinaLescot}.
\end{enumerate}
\end{lemma}
\begin{proof}
(i) The coefficients of \(u_n\) are \(t_n=-i/2\); the first two values
follow from the formulas before Theorem~\ref{thm:globalKahler},
\(g_{c,h}\) is matched with \(2hg_{\rm VK}+\frac c6\pi^*g_{\rm WP}\)
mode by mode, and \(I(u_n,u_n)=C_n/2\) is \eqref{eq:allcov}. Since
\(\theta_n=C_n\) for \((c,h)=(\cL,1)\), (i) follows from
Theorem~\ref{thm:globalKahler} and right invariance.
(ii) By \cite[TS1--TS2, Section~2.2.3]{TTI}, \(u=\sum c_ke^{ik\theta}\partial_\theta\)
corresponds to
\(\mu(z)=-\frac12(1-|z|^2)^2\varphi(1/\bar z)\bar z^{-4}\),
\(\varphi=i\sum_{k\ge2}(k^3-k)c_kz^{k-2}\). With \(w=1/\bar z\) and
\(\iint_\D|w|^{2k-4}(1-|w|^2)^2d^2w=2\pi/(k^3-k)\),
\(\langle\mu,\mu\rangle_{\rm WP}=2\pi\sum_{k\ge2}(k^3-k)|c_k|^2\), which
is \(\frac\pi2(n^3-n)\) for \(u_n\) (\(c_{\pm n}=\frac12\)). Equivalently
\(\partial_{\varepsilon_\mu}\partial_{\bar\varepsilon_\nu}S_1
=\langle\mu,\nu\rangle_{\rm WP}\) \cite[Theorem~5.1]{TT}.
(iii) The Ricci form is
\(\mathrm{Ric}(J\cdot,\cdot)=-i\partial\bar\partial\log\det(g_{k\bar l})
=i\sum R_{k\bar l}\,d\varepsilon_k\wedge d\bar\varepsilon_l\), and
\(R_{k\bar l}=\mathrm{Ric}(\partial_k,\bar\partial_l)\), whereas
\(g^{\rm TT}_{\rm WP}(\partial_k,\bar\partial_l)=\frac12g_{k\bar l}\).
Hence \(\mathrm{Ric}(g_{\rm WP})=-\frac{13}{6\pi}g^{\rm TT}_{\rm WP}\),
and by (ii) its Ricci form is
\(-\frac{13}{6\pi}\omega^{\rm TT}_{\rm WP}
=-i\partial\bar\partial\bigl(\frac{26}{24\pi}S_1\bigr)=-\omega_{26,0}\).
For the factor \(\frac i2\) in the definition of
\(\mathrm{Ric}^{\rm TT}_{\rm WP}\) see Section~\ref{sec:conventionsLit}, (C5).
The value on \((u_n,u_n)\) follows from (ii).
(iv) \cite[Theorem~4.11]{GordinaLescot} gives
\(\mathrm{Ric}(L_n,L_{-n})=-(13n^3-n)/6\); in the normalization of (i)
this is \(-\theta_n/2\) for \((c,h)=(26,1)\). Its WP part is (iii); the
statement on \(\operatorname{M\ddot{o}b}(S^1)\backslash
\operatorname{Diff}_+(S^1)\) is due to \cite{KirillovYuriev,BowickRajeev}
(see \cite[Remark~7.13]{TTI}).
\end{proof}

\begin{corollary}[Decomposition of the Fisher metric]\label{cor:KE}
For \(0<\kappa\le4\), \(I=g_{26,1}-c_{\rm m}\,g_{1,0}\), \(c_{\rm m}=26-\cL\),
where \(g_{1,0}=\mathfrak g_{E/24}\) has value \((n^3-n)/24\) on \((u_n,u_n)\); and
\(I\) is proportional to \(g_{26,1}\) if and only if \(\kappa=8/3\), when \(I=g_{26,1}\).
\end{corollary}
\begin{proof}
By Lemma~\ref{lem:conventions}(i), compare values on \((u_n,u_n)\):
\(C_n=\bigl(2n+\frac{26}{12}(n^3-n)\bigr)-\frac{c_{\rm m}}{12}(n^3-n)\), and
\(\cL=26\) only for \(\kappa=8/3\) in \((0,4]\).
\end{proof}

\begin{remark}[Relation to Kähler--Einstein metrics]\label{rem:regRicci}
By Theorems~\ref{thm:globalKahler} and~\ref{thm:entropyWP}, the Kähler
potential of \(I\) is the relative entropy
\(V=H_\kappa=\operatorname{Ent}(P_{h^{-1}}|P)\), globally and without
regularization. This gives a probabilistic reading of known
Kähler--Einstein statements; it does not prove a new one, and the
curvature of \(I\) on the marked WP space is not computed here.

(a) \emph{The base.} On \(T_0(1)\) the Weil--Petersson metric is
Kähler--Einstein (Lemma~\ref{lem:conventions}(iii)), with Ricci form
\(-\omega_{26,0}\). Ricci forms are unchanged by scaling, so the horizontal
part \(g_{\cL,0}=\frac{\cL}6g_{\rm WP}\) of the Fisher metric is Kähler--Einstein
for every \(\kappa\), with Ricci form \(-\frac{26}{\cL}\omega_{\cL,0}\). Its Einstein
constant \(-26/\cL\) equals \(-1\) at \(\kappa=8/3\), \(-26/25\) at \(\kappa=4\), and
tends to \(0\) as \(\kappa\to0\).

(b) \emph{The value \(\kappa=8/3\).} Write \(\mathrm{Ric}^{\rm reg}=-\omega_{26,1}\) for
the regularized Ricci \emph{form} of \cite{BowickRajeev,GordinaLescot}
(Lemma~\ref{lem:conventions}(iv)), transported by right invariance; it is
the same for all \(g_{c,h}\), and it is not the Ricci tensor
\(\mathrm{Ric}(g_{\rm WP})\) of (a), which accounts only for the \(g_{26,0}\) part.
In this regularized sense \(g_{c,h}\) is Einstein exactly when \(c=26h\). For
the Fisher metric the weight is \(h=1\) for every \(\kappa\), because the
capacity term \(2K\) of the entropy does not depend on \(\kappa\); so the
Einstein condition is \(\cL=26\), that is \(\cm=0\) and \(\kappa=8/3\) in \((0,4]\).
Thus the Fisher metric of the welding family is Kähler--Einstein in the
regularized sense, \(\mathrm{Ric}^{\rm reg}=-\omega_I\), if and only if the matter
charge vanishes, which is the case of exact conformal restriction
(Werner's measure \cite{WangAIF,CP}); the potential of this metric is then
the relative entropy of the SLE\(_{8/3}\) welding laws.

(c) \emph{General \(\kappa\).} Corollary~\ref{cor:KE} reads
\[
 \mathrm{Ric}^{\rm reg}=-\omega_I-\cm\,\omega_{1,0},\qquad
 \omega_{1,0}=i\partial\bar\partial\bigl(E/24\bigr)\ge0,
\]
a twisted Kähler--Einstein equation whose twisting form has potential
\(\frac{\cm}{24}I^L\), the Onsager--Machlup functional of
\cite{CarfagniniWang}: Liouville, matter and ghost (\(-26\)) charges add to
zero. The twist changes sign at \(\kappa=8/3\) (\(\cm<0\) for \(\kappa<8/3\),
\(0<\cm\le1\) for \(8/3<\kappa\le4\)).

(d) \emph{A heuristic.} For a negatively curved Kähler--Einstein metric,
\(\log\det(g_{k\bar l})\) is a Kähler potential up to a pluriharmonic
function. Formally, at \(\kappa=8/3\) the relative entropy
\(\operatorname{Ent}(E_h|P)\) is therefore the logarithm of the regularized volume
density of the Fisher metric. We do not make this precise.

To turn (b)--(c) into theorems on \(\mathcal M\) one would have to extend the
regularized Ricci form from \(S^1\backslash\operatorname{Diff}_+(S^1)\) to the
Weil--Petersson class. 

\end{remark}
\section{On the value of the central charge \(\cL\)}\label{sec:open}

\begin{remark}[An independent test of the central charge]\label{rem:cLtest}
We conclude with a note about the value of the central charge \(\cL=26-\cm\).
This value depends on the normalization of the restriction
covariance (Lemma~\ref{lem:normalization}). The small-\(\kappa\) comparison of
Proposition~\ref{prop:kappa0} is consistent with it, but it does not exclude
the alternative \(26-2\cm\). Neither do the explicit densities of
\cite[Corollary~1.3]{BJweld} and \cite[Theorems~1.9, 1.10]{FanCM}, which contain
the same \(\cL\): both take the matter charge from the restriction covariance
in the convention \(e^{\frac c2\Lambda^*}\), through
\cite[Theorem~4.1]{BJloop} (see \cite[(6.4)]{FanCM}, \cite[(1.1), (1.10)]{BJweld}).
Under the alternative convention they would change together with
\eqref{eq:W}, and the comparison of Proposition~\ref{prop:FanComparison} would
remain valid (Remark~\ref{rem:FanTests}). The Gaussian multiplicative chaos
inputs for \(\kappa<4\) \cite{FS}, \cite[Lemma~1.2]{FanCM} concern only the measure
class. A test that discriminates the two would be valuable:
\begin{itemize}
\item an upper bound \(\limsup_{\kappa\to0}\kappa\operatorname{Ent}(P_a|P)\le I^L(a)\);
\item a derivation of the coefficient of \(I^L\) from the Liouville side, as
\(1+6Q^2\) with \(Q=2/\sqrt\kappa+\sqrt\kappa/2\), through a Weil--Petersson
deformation of the quantum-disk welding \cite{AHSloop,FS} (note
\(1+6Q^2=26-\cm\neq26-2\cm\));
\item a comparison of a score identity from \eqref{eq:W} with an exactly
solvable observable of SLE loops;
\item a numerical estimate at small \(\kappa\) of the variance
\(\E|\beta_1|^2=1/(C-3)\) of Proposition~\ref{prop:moments}, which differs
between the two conventions by about \(31\%\) at \(\kappa=1\) (the fourth
moment \(\E|\alpha_1|^4\) differs only by a few per cent).
\end{itemize}

\end{remark}
\appendix
\section{The M\"obius likelihood and coefficient moments}
\label{app:mobius}\label{app:moments}

These are consistency checks of the normalization; nothing else depends
on them. By Corollary~\ref{cor:entMobius}, \(M=e^{2K}P\) is invariant under
left M\"obius composition (\(M\) is infinite, but its weighted pushforward
is the probability in question), so for \(m_t(z)=(z+\bar t)/(1+tz)\),
\(|t|<1\), \(r_{m_t}(h)=dP_{m_t}/dP(h)=e^{2K(h)-2K(m_t^{-1}\circ h)}\). The
normalizer \(W_t(w)=w/(1-w/g(-1/t))\) fixes \(0\) with derivative one, and
the exterior map \(W_t\circ g\circ m_t^{-1}\) has leading coefficient
\(b_t=t^2g(-1/t)^2/((1-|t|^2)g'(-1/t))\), by
\(m_t^{-1}(z)=-1/t-(1-|t|^2)/(t^2z)+O(z^{-2})\). With \(m_t^{-1}=m_{-t}\) and
\(G'(t)=b\,g'(1/t)/(t^2g(1/t)^2)\),
\begin{equation}\label{eq:mobiuscheck}
 r_{m_t}(h)=(1-|t|^2)^2|G_h'(t)|^2 .
\end{equation}
Schlicht distortion bounds fixed derivatives uniformly over samples, so
\(\sqrt{r_{m_t}}=1+t\alpha_1+\bar t\bar\alpha_1+O_{L^\infty}(|t|^2)\), with complex
score \(2\alpha_1=-2\beta_0=s_1\). Normalization gives
\(\E|G'(t)|^2=(1-|t|^2)^{-2}\), and comparing series in \(t,\bar t\),
\(\E[\alpha_j\bar\alpha_k]=\delta_{jk}/(j+1)\) and
\(\E[G'(z)\overline{G'(w)}]=(1-z\bar w)^{-2}\). Hence \(\E|s_1|^2=2\) and
\(H^2(P_{m_t},P)=|t|^2+O(|t|^3)\); the coefficient pattern has precedent in
\cite[Theorem 6.1]{CP}.

\begin{proposition}\label{prop:moments}
For \(C=C_2=4+\cL/2\),
\(\E|\alpha_1|^2=\tfrac12\), \(\E|\alpha_2|^2=\E[\alpha_1^2\bar\alpha_2]=\tfrac13\) and
\(\E|\alpha_1|^4=C/(3(C-3))\). Equivalently,
\(\E\bigl[(\beta_0^2,\beta_1)^{\mathsf T}\,\overline{(\beta_0^2,\beta_1)}\bigr]
=\frac1{C-3}\bigl(\begin{smallmatrix}C/3&1\\1&1\end{smallmatrix}\bigr)\).
\end{proposition}
\begin{proof}
The second moments were obtained above. For the normalized first-mode
variation of the reciprocal exterior map, \(D_1G=1-G'+2\alpha_1G\) and
\(\bar D_1G=z^2G'\), so \(D_1\alpha_1=2\alpha_1^2-3\alpha_2\) and \(\bar D_1\alpha_1=1\).
Differentiating the explicit density \eqref{eq:mobiuscheck} gives
\(\E\bar D_1F=\E[2\bar\alpha_1F]\) for bounded coefficient polynomials; with
\(F=\alpha_1^2\bar\alpha_1\), \(Y=\E[\alpha_1^2\bar\alpha_2]\) and \(X=\E|\alpha_1|^4\), this reads
\(2\E|\alpha_1|^2+2X-3Y=2X\), so \(Y=\frac13\). Finally
\(s_2=C\alpha_2+(3-C)\alpha_1^2\) and \(\E|s_2|^2=C\) give
\(C=\frac{C^2}3+\frac{2C(3-C)}3+(C-3)^2X\), and \(\alpha_1=-\beta_0\),
\(\alpha_2=\beta_0^2-\beta_1\) give the Gram matrix.
\end{proof}

At \(\kappa=8/3\), \(\cm=0\), \(C=17\), and \(\E|\alpha_1|^4=17/42\), the Werner
moment of \cite[Proposition 6.4]{CP}.
\section{The conformal potential at the round circle}\label{sec:round}

The K\"ahler-potential theorem is due to Takhtajan--Teo
\cite[Theorems 4.8 and 5.1]{TT}; see also \cite{Wang}. The following
calculation fixes the metric convention at the round circle; it does
not evaluate a conformal energy on an SLE sample.

\begin{proposition}[Round-circle comparison]\label{prop:potential}
For smooth welding maps let \(S_1\) denote the pullback to the
marked curve of the universal Liouville action on the M\"obius
quotient. In the indicated welding normalization, write
\[
 K=\log|g'(\infty)/f'(0)|,\qquad
 S_1=\int_\D|f''/f'|^2\dd
     +\int_\Ds|g''/g'|^2\dd-4\pi K,
 \qquad V=2K+\frac{\cL}{24\pi}S_1.
\]
On any fixed finite Fourier tangent space at the round circle, the
real Hessian of \(V\) equals the Fisher matrix of
Theorem~\ref{thm:finite}.
\end{proposition}
\begin{proof}
Take \(f_{\mathbf t}(z)=z-\sum_{n\le N}t_nz^{n+1}\), univalent on a
fixed larger disc for small \(\mathbf t\), a holomorphic path in the
standard chart. For a single mode the boundary welding and
the exterior map have first jets
\[
 h_t(z)=z-tz^{n+1}+\bar t z^{1-n},\qquad
 g_t(z)=z-\bar t z^{1-n},
\]
with remainders of higher order. The interior area divided by \(\pi\)
is \(1+(n+1)|t|^2\). For the exterior Laurent coefficients,
the same area divided by \(\pi\) is
\(b_t^2-\sum_{k\ge1}k|b_{-k}(t)|^2\), where \(b_t>0\).
For \(n\ge2\), the first jet has \(b_{-(n-1)}=-\bar t\);
for \(n=1\) it is a constant translation with zero area weight.
Analytic collar dependence bounds the remaining weighted coefficient
sum by \(o(|t|^2)\). Thus
\(g_t'(\infty)^2=b_t^2=1+2n|t|^2+o(|t|^2)\).
Hence \(K=n|t|^2+o(|t|^2)\). The two pre-Schwarzian energy
coefficients are \(\pi n(n+1)^2\) and \(\pi n(n-1)^2\),
so \(S_1=2\pi(n^3-n)|t|^2+o(|t|^2)\).
Analytic collar dependence justifies these expansions in the energy
norms. In a finite sum of modes, angular orthogonality removes all
cross terms. Thus
\[
 V(t)=\sum_{n=1}^NC_n|t_n|^2+o(|t|^2).
\]
Its real Hessian has blocks \(2C_nI_2\), as does the Fisher matrix.
With \(X=\partial_{\Rea t_n}\), \(JX=\partial_{\Ima t_n}\), one has
\(X^{1,0}=\partial_{t_n}\); hence
\(\partial_{t_n}\partial_{\bar t_n}V=C_n\), and twice its real
part is \(2C_n\). This fixes the conversion to the mixed-Hessian
convention used in Theorem~\ref{thm:globalKahler}.
The first derivatives of these potentials vanish at the round circle,
so only the first jet of the chosen welding path matters.
\end{proof}
\section{The constant of the restriction covariance}\label{app:normalization}

This appendix verifies, from Zhan's proof, the constant \(\frac c2\) in the
conformal restriction covariance that is used in the literature
\cite[(2.3)]{CarfagniniWang}, \cite[Definition~2.1]{GQW},
\cite[Definition~1.8]{BJloop}, and in Proposition~\ref{prop:Wp-loopWard};
see Section~\ref{sec:conventionsLit}, (C4).

\begin{lemma}[Normalization of the restriction covariance]\label{lem:normalization}
Let \(0<\kappa\le4\) and \(c=c(\kappa)=\cm\). Let \(\mu^{\rm loop}\) be the
Brownian loop measure, normalized as in \cite[Section~4.1, (15)]{LawlerWerner}
and \cite[Section~2.1.3]{FieldLawler} (area \(\times\frac{dt}{2\pi t^2}\times\)
Brownian bridge; the measure \(\mu\) of Appendix~\ref{app:caiInput}), and let
\(\mathcal B\) and \(\Lambda^*\) be as there, \(\Lambda^*\) being the
normalization of \cite[Theorem~1.1, (1.6)]{FieldLawler}.
For a domain \(D\subset\widehat{\mathbb C}\) with non-polar complement put
\[
 \mu^0_D=\mathbf 1_{\{\gamma\subset D\}}\,
 e^{\frac c2\Lambda^*(\gamma,\widehat{\mathbb C}\setminus D)}\,\mu^0 .
\]
\begin{enumerate}
\item If \(\psi:D\to D'\) is conformal, then \(\psi_*\mu^0_D=\mu^0_{D'}\).
\item If \(\varphi:D\to D'\) is conformal and \(\Psi\ge0\) is a Borel function on
loops in \(D'\), then
\[
 \int_{\gamma\subset D}\Psi(\varphi\gamma)\,d\mu^0(\gamma)
 =\int_{\gamma'\subset D'}\Psi(\gamma')\,
 e^{\frac c2[\Lambda^*(\gamma',\widehat{\mathbb C}\setminus D')
       -\Lambda^*(\varphi^{-1}\gamma',\widehat{\mathbb C}\setminus D)]}\,d\mu^0(\gamma').
\]
\end{enumerate}
Thus \(\mu^0\) has the properties (RC) and (CI) of
\cite[Definition~2.1, (1)]{GQW} with charge \(c\), and (ii) is
\cite[Lemma~2.2]{GQW}.
\end{lemma}

\begin{proof}
(ii) follows from (i): for \(\gamma\subset D\) the sets \(\gamma\) and
\(\widehat{\mathbb C}\setminus D\) are disjoint non-polar compacta, so
\(\Lambda^*(\gamma,\widehat{\mathbb C}\setminus D)\in\mathbb R\), and
\[
 \int_{\gamma\subset D}\Psi(\varphi\gamma)\,d\mu^0
 =\int\Psi(\varphi\gamma)e^{-\frac c2\Lambda^*(\gamma,\widehat{\mathbb C}\setminus D)}d\mu^0_D(\gamma)
 =\int\Psi(\gamma')e^{-\frac c2\Lambda^*(\varphi^{-1}\gamma',\widehat{\mathbb C}\setminus D)}d\mu^0_{D'}(\gamma').
\]

(i) We show that the proof of \cite[Theorem~5.1, (5.2)]{ZhanLoop} gives the
exponent \(\frac c2\) for the Field--Lawler \(\Lambda^*\).

\emph{Step 1 (where the loop measure enters).} For \(\kappa\le4\) the proof of
\cite[Theorem~5.1]{ZhanLoop} combines the conformal Markov property of the
rooted measure \(\mu^1_0\), the distortion of whole-plane SLE\(_\kappa(2)\)
under a conformal map \cite[Lemmas~3.2--3.3]{ZhanLoop}, and the chordal
covariance \cite[Proposition~A.3]{ZhanLoop}, transferred to \(\Ds\) in
\cite[Lemma~A.5]{ZhanLoop}; the passage from rooted to unrooted measures
\cite[(4.6)]{ZhanLoop} does not involve loops. The Brownian loop measure
\(\mu^{\rm lp}\) of \cite{ZhanLoop} enters in three ways only.
(a) Through \eqref{eq:LambdaNest} and the conformal invariance of
\(\mathcal B\); both are linear in the loop measure.
(b) Through the chordal measures
\(\mathbf 1_{\{\gamma\cap(D\setminus U)=\emptyset\}}e^{c\mu^{\rm lp}(L_D(\gamma,D\setminus U))}\mu^\#_{D;a\to b}\)
of \cite[(A.3)]{ZhanLoop}, whose covariance \cite[Proposition~A.3]{ZhanLoop} is
quoted for \(\kappa\le4\) from Lawler's work. In Lawler's formulation this is
the boundary perturbation rule
\(\frac{d\mu_1}{d\mu}(\gamma)=\mathbf 1_{\{\gamma\subset D_1\}}\exp\{\frac c2\Lambda(\gamma,D\setminus D_1;D)\}\)
\cite[(1.4)]{FieldLawler} (see also \cite[Section~2, (2)]{LawlerMC} and
\cite{LawlerPF}), where \(\Lambda(\gamma,D\setminus D_1;D)\) is the
\(\mu^{\rm loop}\)-mass of the loops in \(D\) that meet \(\gamma\) and
\(D\setminus D_1\), for the measure \(\mu^{\rm loop}\) of
\cite[Section~2.1.3]{FieldLawler}, which is that of \cite[(15)]{LawlerWerner}.
So the exponent is \(\frac c2\mu^{\rm loop}\). (We do not use the loop-soup
forms of the restriction formula \cite[Theorem~7.3]{LSW},
\cite[Theorem~1]{LawlerWerner}, whose constant depends on a factor-two
orientation convention for loops \cite[Section~1.3]{SheffieldWerner}.)
(c) Through \cite[Lemma~3.4]{ZhanLoop}, which turns the exponent
\(cE(\tau)\),
\(E(\tau)=\frac16\int_{-\infty}^\tau\mathcal S(\widetilde W_t)(\lambda_t)\,dt
+\frac1{12}\int_{-\infty}^\tau(\widetilde W_t'(\lambda_t)^2-1)\,dt\), of Zhan's
martingale into \(c[\Lambda^*(\beta_{u(\tau)},V^c)-\Lambda^*(\gamma_\tau,U^c)]\)
(see \cite[(5.6)]{ZhanLoop}). Its proof writes, in \cite[(3.35)--(3.36)]{ZhanLoop},
\[
 \mathcal B\bigl(\gamma[a,\tau],U^c;\widehat{\mathbb C}\setminus\gamma_a\bigr)
 =k\int_a^\tau\mu^{\rm bub}_{\lambda_t}\bigl(L(\Ds\setminus U_t)\bigr)\,dt ,
\]
and the same for \(\beta\), with \(k=1\). Here \(\mu^{\rm bub}_x\) is the
Brownian bubble measure in \(\Ds\) rooted at \(e^{ix}\) \cite{LawlerWerner}.
Its normalization is the Poisson-kernel formula displayed (unnumbered) right
after \cite[(3.36)]{ZhanLoop} in the proof of \cite[Lemma~3.4]{ZhanLoop}:
for a subdomain \(U\subset\Ds\) containing a neighbourhood of \(S^1\) in
\(\Ds\),
\begin{equation}\label{eq:bubnorm}
 \mu^{\rm bub}_x\bigl(L(\Ds\setminus U)\bigr)
 =\lim_{U\ni z\to e^{ix}}\frac1{|z-e^{ix}|^{2}}
 \Bigl(1-\frac{P^U_x(z)}{P^{\Ds}_x(z)}\Bigr),
 \qquad
 P^{\Ds}_x(z)=\frac1{2\pi}\Rea\frac{z+e^{ix}}{z-e^{ix}},
\end{equation}
where \(L(\Ds\setminus U)\) is the set of bubbles that meet \(\Ds\setminus U\)
and \(P^U_x\) is the Poisson kernel of \(U\) with pole at \(e^{ix}\).
This is the normalization used in Step~2. The decomposition theorem of
\cite{LawlerWerner} gives this representation with a constant \(k\) that does
not depend on the curve, and the remaining Poisson-kernel computation in that
proof then gives \(\Lambda^*(\beta_{u(\tau)},V^c)-\Lambda^*(\gamma_\tau,U^c)=kE(\tau)\).

\emph{Step 2 (\(k=2\)).} Take the constant driving function \(\lambda\equiv0\).
Then \(g_t^{-1}(w)=e^t(w+1)^2/w\), so \(\gamma_t=[0,4e^t]\) and
\(\lim_{t\to-\infty}e^tg_t(z)=z\), as in \cite[(3.1)]{ZhanLoop} and the display before it. Take
\(U=\D\), so that \(U^c=\widehat{\mathbb C}\setminus\D\) and
\(\Ds\setminus U_t=g_t(U^c)\) for \(4e^t<1\), and put
\(G(x)=\Lambda^*(\overline{B(0,x)},U^c)\) for \(0<x<1\).

(\(\alpha\)) \(G(x)+\log\log(1/x)\to0\) as \(x\downarrow0\). Indeed, for \(0<r<x\),
\eqref{eq:LambdaNest} with \(K=U^c\), \(D=\widehat{\mathbb C}\setminus\overline{B(0,r)}\),
\(D'=\widehat{\mathbb C}\setminus\overline{B(0,x)}\) gives
\(G(x)=G(r)+\mathcal B(U^c,\overline{B(0,x)};\widehat{\mathbb C}\setminus\overline{B(0,r)})\),
since loops avoiding \(\overline{B(0,r)}\) that meet \(\overline{B(0,x)}\) meet
\(\overline{B(0,x)}\setminus\overline{B(0,r)}\). By \cite[Theorem~1.1, (1.6)]{FieldLawler},
\(\Lambda^*(V_1,V_2)=\lim_{r\downarrow0}[\Lambda(V_1,V_2;\mathcal O_r)-\log\log(1/r)]\),
\(\mathcal O_r=\{|z|>r\}\), where \(\Lambda(V_1,V_2;\mathcal O_r)\) is the
\(\mu^{\rm loop}\)-mass of loops in \(\mathcal O_r\) meeting \(V_1\) and \(V_2\).
That theorem is stated for all disjoint non-polar closed sets
\(V_1,V_2\subset\widehat{\mathbb C}\), with the base point \(0\) and no
condition on where \(0\) lies; so it applies to \(V_1=U^c\) and
\(V_2=\overline{B(0,x)}\ni0\). In this case, for \(r<x\),
\(\Lambda(U^c,\overline{B(0,x)};\mathcal O_r)=\mathcal B(U^c,\overline{B(0,x)};\widehat{\mathbb C}\setminus\overline{B(0,r)})\),
since \(\mu^{\rm loop}\)-almost no loop passes through \(\infty\). The same
definition with an arbitrary base point \(z_0\in\mathbb C\) is
\cite[(2.1)]{CarfagniniWang}. Its value does not depend on \(z_0\)
\cite[(1.2)]{ZhanLoop}, which also follows from the M\"obius invariance in
\cite[Theorem~1.1]{FieldLawler} (see the remark after it). Hence
\(\mathcal B(U^c,\overline{B(0,x)};\widehat{\mathbb C}\setminus\overline{B(0,r)})-\log\log(1/r)\to G(x)\)
as \(r\downarrow0\), and therefore \(G(r)+\log\log(1/r)\to0\).

(\(\beta\)) \(\Lambda^*(\gamma_a,U^c)+\log(-a)\to0\) as \(a\to-\infty\). Since
\((|w|-1)^2/|w|\le|(w+1)^2/w|\le(|w|+1)^2/|w|\), we have
\(\{|w|\ge R_+\}\subset g_a(U^c)\subset\{|w|\ge R_-\}\), where \(R_\pm>1\) solve
\((R_+-1)^2/R_+=e^{-a}=(R_-+1)^2/R_-\); thus \(\log R_\pm=-a+O(e^a)\). Also
\(g_a(\overline{B(0,4e^a)}\setminus\gamma_a)\subset\{|w|\le6\}\). Since
\(\gamma_a\subset\overline{B(0,4e^a)}\), \eqref{eq:LambdaNest} with \(K=U^c\),
\(D=\widehat{\mathbb C}\setminus\gamma_a\), \(D'=\widehat{\mathbb C}\setminus\overline{B(0,4e^a)}\),
and conformal invariance give
\[
 0\le G(4e^a)-\Lambda^*(\gamma_a,U^c)
 =\mathcal B\bigl(g_a(U^c),g_a(\overline{B(0,4e^a)}\setminus\gamma_a);\Ds\bigr)
 \le\mathcal B\bigl(\{|w|\ge R_-\},\{|w|\le6\};\Ds\bigr).
\]
After the inversion \(w\mapsto1/w\), the right side is
\(\mathcal B(\overline{B(0,1/R_-)},\{1/6\le|w|<1\};\D)=G(6/R_-)-G(1/R_-)\), by
\eqref{eq:LambdaNest} (\(K=\overline{B(0,1/R_-)}\), \(D=\D\), \(D'=B(0,1/6)\)) and
the dilation invariance of \(\Lambda^*\). By (\(\alpha\)) it tends to \(0\), and
\(G(4e^a)=-\log(-a)+o(1)\).

(\(\gamma\)) \(\mu^{\rm bub}_0(L(\Ds\setminus U_t))=\frac1{2|t|}+O(t^{-2})\) as
\(t\to-\infty\). For \(R>1\) let \(U_R=\{1<|w|<R\}\) and
\(m(R)=\mu^{\rm bub}_0(L(\Ds\setminus U_R))\). By \eqref{eq:bubnorm} with
\(x=0\) and \(z=r\downarrow1\),
\[
 m(R)=\lim_{r\downarrow1}\frac{1}{(r-1)^2}\,
 \frac{P^{\Ds}_0(r)-P^{U_R}_0(r)}{P^{\Ds}_0(r)},
 \qquad P^{\Ds}_0(r)=\frac{r+1}{2\pi(r-1)} .
\]
For \(1<r<R\), by the strong Markov
property, \(P^{\Ds}_0(r)-P^{U_R}_0(r)=\E^r[P^{\Ds}_0(B_T);T<\sigma]\), where \(T\) and
\(\sigma\) are the hitting times of \(\{|w|=R\}\) and \(S^1\). Here
\(\mathbb P^r(T<\sigma)=\log r/\log R\), and
\(P^{\Ds}_0\) lies between \(\frac{R-1}{2\pi(R+1)}\) and \(\frac{R+1}{2\pi(R-1)}\) on
\(\{|w|=R\}\). Since \(\log r/((r-1)(r+1))\to\frac12\) as \(r\downarrow1\), this
gives \(m(R)=(1+O(1/R))/(2\log R)\), which is all we use. (The Fourier
expansion of \(P^{U_R}_0\) in the annulus gives the exact value
\(m(R)=\frac1{2\log R}+2\sum_{n\ge1}\frac{nR^{-n}}{R^n-R^{-n}}=\frac1{2\log R}+O(R^{-2})\).)
Since \(\log R_\pm=|t|+O(e^{t})\), either form gives the error \(O(t^{-2})\)
in (\(\gamma\)). The
bubble mass increases with the removed set, because \(P^U\) decreases as \(U\)
decreases; so the inclusions in (\(\beta\)) at time \(t\) give
\(m(R_+)\le\mu^{\rm bub}_0(L(\Ds\setminus U_t))\le m(R_-)\), and (\(\gamma\)) follows.

Now fix \(\tau<-\log4\). By \eqref{eq:LambdaNest} and (\(\beta\)),
\(\mathcal B(\gamma[a,\tau],U^c;\widehat{\mathbb C}\setminus\gamma_a)
=\Lambda^*(\gamma_\tau,U^c)-\Lambda^*(\gamma_a,U^c)=\log(-a)+O(1)\), whereas by
(\(\gamma\)) \(k\int_a^\tau\mu^{\rm bub}_0\,dt=\frac k2\log(-a)+O(1)\) as
\(a\to-\infty\). Hence \(k=2\). With \(k=1\) the right side of
\cite[(3.35)]{ZhanLoop} would tend to \(-\infty\). So
\cite[(3.35)--(3.36)]{ZhanLoop}, as displayed there, hold for the loop
measure \(\frac12\mu^{\rm loop}\) (Step~3).

\emph{Step 3 (conclusion).} Put \(\hat\mu=\frac12\mu^{\rm loop}\),
\(\hat{\mathcal B}=\frac12\mathcal B\) and \(\hat\Lambda=\frac12\Lambda^*\). By
Step~2, \cite[(3.35)--(3.36)]{ZhanLoop} hold for \(\hat{\mathcal B}\) and
\(\hat\Lambda\) as displayed, with \(\log(-a)\) replaced by \(\frac12\log(-a)\),
a term that cancels in the difference; so \cite[Lemma~3.4]{ZhanLoop} holds
verbatim for \(\hat\Lambda\). By (b), \cite[(A.3), Proposition~A.3]{ZhanLoop} hold
verbatim for \(\hat\mu\), and (a) holds for \(\hat\mu\), \(\hat{\mathcal B}\),
\(\hat\Lambda\). No other property of the loop measure is used, so the proof of
\cite[Theorem~5.1]{ZhanLoop} gives (5.1)--(5.2) there with
\(e^{c\hat\Lambda}=e^{\frac c2\Lambda^*}\). This is (i).
\end{proof}

\begin{remark}[Normalization]\label{rem:Wp-normalization}
(a) The constant \(\frac{\cm}6\) in Proposition~\ref{prop:Wp-loopWard} is the
product of \(\frac c2\) and the coefficient \(\frac13\) of
Lemma~\ref{lem:GQWanalytic}; by \eqref{eq:LambdaIL} it is \(\frac c{24}\) times
\(\frac{d}{dh}I^L(\phi_h\zeta)=4\Rea(\lambda P_p)\). This product is independent
of the normalization of the loop measure (rescaling \(\mu^{\rm loop}\) changes
the two factors inversely), and in this form it is the coefficient
\(\frac c{24}\) of the Onsager--Machlup functional
\cite[Theorem~1.1]{CarfagniniWang} and of the formal density
\(e^{\frac c{24}I^L}\) \cite[(1.9)]{BJloop}, \cite[(1.10)]{BJweld}, and the
factor \(\frac{\cm}{12}\) in \cite[Theorem~4.1]{BJloop}. These sources use the
same convention, so the agreement is not an independent check.
(b) No identity internal to this paper fixes the constant independently.
The M\"obius invariance (P1) gives only the case \(p=1\) of
Proposition~\ref{prop:Wp-loopWard} (\(P_1=0\)) and \eqref{eq:Wp-rot}, hence
only the first-mode moments of Appendix~\ref{app:moments}; every quantity
depending on \(\cL\) (the \(C_n\), \(n\ge2\), and \(\E|\alpha_1|^4\)) passes through
Lemma~\ref{lem:normalization}. The Werner-moment check is at \(\kappa=8/3\),
where \(c=0\).
\end{remark}

We end this appendix with a small-\(\kappa\) consistency check. The next proposition
records the \(\kappa\to0\) behaviour of Theorem~\ref{thm:entropyWP}; it is
compared in Remark~\ref{rem:kappa0} with an exact chordal identity, which
involves neither the loop measure nor the central charge.

\begin{proposition}[Small-\(\kappa\) limit]\label{prop:kappa0}
For \(0<\kappa\le4\) and every WP homeomorphism \(a\),
\[
 \kappa\operatorname{Ent}(P_a|P)=\Bigl(1+\frac{13\kappa}{24}+\frac{\kappa^2}{16}\Bigr)I^L(a)+2\kappa K(a).
\]
In particular \(\kappa\operatorname{Ent}(P_a|P)\to I^L(a)\) as \(\kappa\to0\),
uniformly on sets where \(I^L\) and \(K\) are bounded. If the exponent
\(\frac c2\Lambda^*\) in Lemma~\ref{lem:normalization} is replaced by
\(c\Lambda^*\), the proofs of Sections~\ref{app:Wproof}
and~\ref{sec:entropy} give instead
\(\kappa\operatorname{Ent}(P_a|P)=(2+\frac{\kappa^2}8)I^L(a)+2\kappa K(a)\),
with limit \(2I^L(a)\).
\end{proposition}

\begin{proof}
Since \(\kappa\cm=-\frac12(3\kappa^2-26\kappa+48)\), we have
\(\kappa\cL/24=(3\kappa^2+26\kappa+48)/48\), and the identity follows from
Theorem~\ref{thm:entropyWP}. Now suppose the exponent is \(c\Lambda^*\). Then
the factor \(\frac{\cm}2\) in the proof of Proposition~\ref{prop:Wp-loopWard}
becomes \(\cm\), and the constant in \eqref{eq:Wp-loopWard} becomes
\(-\frac{\cm}3\). In the proof of Theorem~\ref{thm:Wproof} the term
\(-\frac{\cm}{12}M_{n,p}P_p\) becomes \(-\frac{\cm}6M_{n,p}P_p\), and by
Lemmas~\ref{lem:Wp-basis} and~\ref{lem:Wp-ghost} the last display there
becomes
\(-[z^{n-2}]\{\frac{26-2\cm}{12}\mathcal Sf+\frac{f'^2}{f^2}-\frac1{z^2}\}\).
The parameter \(\cL\) enters Sections~\ref{sec:finiteinput}
and~\ref{sec:entropy} only through the stress coefficients
\eqref{eq:stress}, that is, through \eqref{eq:W}. So
Theorem~\ref{thm:entropyWP} holds with \(26-2\cm\) in place of \(\cL\), and
\(\kappa(26-2\cm)/24=(3\kappa^2+48)/24\).
\end{proof}

\begin{remark}[Chordal comparison; a heuristic]\label{rem:kappa0}
For chordal SLE\(_\kappa\) driven by \(\sqrt\kappa B+\lambda\), \(\lambda(0)=0\),
\(I(\lambda)=\frac12\int_0^\infty\dot\lambda^2dt<\infty\), the Cameron--Martin theorem and
the injectivity of the Loewner map, which has the Borel inverse
\(W_t=\lim_{s\downarrow0}g_t(\gamma(t+s))\) \cite{RohdeSchramm}, give the exact identity
\(\operatorname{Ent}(Q^\lambda_\kappa|Q^0_\kappa)=I(\lambda)/\kappa\) for the laws of the traces. Its rate
\(I(\lambda)\), the chordal Loewner energy, is the large-deviation rate function
of SLE\(_{0+}\) \cite[Proposition~3.2]{WangJEMS}, \cite[Theorem~1.5]{PeltolaWang},
and any such large-deviation upper bound gives, by data processing on a
neighbourhood of the limit point, the lower bound
\(\liminf\kappa\operatorname{Ent}\ge I\) for deformations converging to a point.
Proposition~\ref{prop:kappa0} shows that \(\frac c2\Lambda^*\) gives the same
leading coefficient, whereas \(c\Lambda^*\) would give \(2I^L(a)\). This is a
heuristic consistency check, not a calibration: the deformations differ;
no large-deviation principle for \(P\) is used or known to us; and even
granted one, the lower bound \(\liminf\kappa\operatorname{Ent}(P_a|P)\ge I^L(a)\) holds for
both leading coefficients. The limit sees only the leading coefficient of
\(\cL/24\), not the ghost term \(\frac{13}{24}\) (Lemma~\ref{lem:Wp-ghost}).
\end{remark}
\section{The first variation of \(\Lambda^*\) and the score covariance}\label{app:caiInput}

This appendix proves the first-variation formula used in
Section~\ref{app:Wproof} (Lemmas~\ref{lem:GQWanalytic} and~\ref{lem:GQWext})
and derives the score covariance from \eqref{eq:W}. Here \(f=f_h=z(1+\sum_{m\ge1}a_mz^m)\) is the
normalized interior map, \(\mathfrak t_n\) its stress coefficients
(Lemma~\ref{lem:CaiDictionary}), and \(c_{\rm m}=26-\cL\).

We begin with the first variation of \(\Lambda^*\) along Jordan curves.
Let \(\mu\) be the Brownian loop measure on \(\widehat{\mathbb C}\) and,
for disjoint non-polar compact sets \(V_1,V_2\), let
\(\Lambda^*(V_1,V_2)\in\mathbb R\) be its Field--Lawler normalization
\cite{FieldLawler}, \cite[Section~2.1]{GQW}. For a domain \(D\) and sets
\(A,B\) let \(\mathcal B(A,B;D)\) be the \(\mu\)-mass of loops contained in
\(D\) that meet \(A\) and \(B\). We use two facts. First, if
\(D'\subset D\) are domains with non-polar complements and \(K\subset D'\) is compact and
non-polar, then
\begin{equation}\label{eq:LambdaNest}
 \Lambda^*(K,\widehat{\mathbb C}\setminus D')
 =\Lambda^*(K,\widehat{\mathbb C}\setminus D)+\mathcal B(K,D\setminus D';D)
\end{equation}
\cite[Lemma~2.1]{CarfagniniWang}. (This follows from the definition: at
each cutoff the loops avoiding a small ball split into those meeting
\(\widehat{\mathbb C}\setminus D\) and those contained in \(D\) that meet
\(D\setminus D'\); the \(\log\log\) terms cancel and the second family
increases to its limit, since \(\mu\)-almost no loop passes through a
given point.) Second, if \(\psi:D\to\psi(D)\) is conformal, then \(\psi\)
pushes \(\mu\) restricted to loops in \(D\) to \(\mu\) restricted to loops
in \(\psi(D)\), modulo time change (\cite{LawlerWerner}; on
\(\widehat{\mathbb C}\) this follows from the planar statement by an
inversion, since \(\mu\)-almost no loop passes through \(\infty\)); hence
\(\mathcal B(\psi A,\psi B;\psi D)=\mathcal B(A,B;D)\).

Fix \(p\ge2\), \(\lambda\in\{1,i\}\) and \(\delta>0\). Let \(\phi_t\) be the
flow of the holomorphic field \(-\lambda z^{1-p}\partial_z\); this is the flow
\(\varphi^{p,\lambda}_t\) (or \(\varphi_t\)) of Section~\ref{app:Wproof}, written
\(\phi_t\) here, while \(\varphi_{-p,t}\) below is the notation of \cite{GQW}. Along a
trajectory \(w(t)=\phi_t(z)\) one has \(w(t)^p=z^p-p\lambda t\). Hence, if
\(2p\varepsilon\le(\delta/4)^p-(\delta/8)^p\), then for \(|t|\le2\varepsilon\) the
map \(\phi_t\) is defined and conformal on
\(U_0=\{|z|>\delta/4\}\cup\{\infty\}\), with \(\phi_t(\infty)=\infty\). The
flow rule \(\phi_h\circ\phi_s=\phi_{s+h}\) holds on \(\{|z|>\delta/2\}\) for
\(|s|,|h|\le\varepsilon\), and \(\phi_s(D_\delta)\subset D_{\delta/2}\), where
\(D_\rho=\{|z|>\rho\}\cup\{\infty\}\). For \(\lambda=1\), resp.\
\(\lambda=i\), \(\phi_t\) is the flow \(\varphi_{-p,t}\), resp.\
\(\psi_{-p,t}\), of \cite[Section~3.3]{GQW}.
Put
\[
 K_t(\rho)=\widehat{\mathbb C}\setminus\phi_t(D_\rho),\qquad
 A_t(\eta)=\Lambda^*\bigl(\phi_t\eta,K_t(\delta)\bigr)-\Lambda^*\bigl(\eta,K_0(\delta)\bigr).
\]
(For the relation of \(K_t(\delta)\) to \cite[Corollary~3.9]{GQW} see
Section~\ref{sec:conventionsLit}, (C3).)
For \(\rho>0\) let \(\mathcal J_\rho\) be the set of Jordan curves \(\eta\subset\mathbb C\)
with \(\overline{B(0,\rho)}\subset\operatorname{int}\eta\). For
\(\eta\in\mathcal J_{\delta/2}\) let \(F_\eta:\D\to\operatorname{int}\eta\) be conformal with
\(F_\eta(0)=0\). Define
\[
 P_p(\eta)=\frac1{2\pi i}\oint_{|w|=\delta/4}w^{1-p}\,
            \mathcal S[F_\eta^{-1}](w)\,dw ,
\]
which is the coefficient of \(z^p\) in \(z^2\mathcal S[F_\eta^{-1}](z)\)
\cite[Definition~3.6]{GQW}.
It does not depend on the choice of \(F_\eta\), since
\(\mathcal S[M\circ g]=\mathcal S[g]\) for rotations \(M\). By
\cite[Lemma~3.7]{GQW}, it is a polynomial in the ratios \(b_k/b_1\) of finitely many Taylor
coefficients \(b_k\) of \(F_\eta^{-1}\) at \(0\); by the Kraus--Nehari bound
\(|\mathcal S[F_\eta^{-1}]|\) is bounded on \(|w|=\delta/4\) uniformly in
\(\eta\in\mathcal J_{\delta/2}\), so \(\sup_\eta|P_p(\eta)|<\infty\).

\begin{lemma}[First variation of \(\Lambda^*\) for analytic curves]\label{lem:GQWanalytic}
Let \(\zeta\in\mathcal J_{\delta/2}\) be an analytic Jordan curve. Then, for
\(|h|\le2\varepsilon\),
\begin{equation}\label{eq:LambdaIL}
 \Lambda^*\bigl(\phi_h\zeta,K_h(\tfrac\delta2)\bigr)-\Lambda^*\bigl(\zeta,K_0(\tfrac\delta2)\bigr)
 =-\tfrac1{12}\bigl(I^L(\phi_h\zeta)-I^L(\zeta)\bigr),
\end{equation}
where \(I^L\) is the Loewner energy \cite{WangInv}, and
\(\frac{d}{dh}\big|_{h=0}I^L(\phi_h\zeta)=4\Rea(\lambda P_p(\zeta))\). Consequently,
with \(\rho=\delta/2\),
\begin{equation}\label{eq:GQWanalytic}
 \frac{d}{dh}\Big|_{h=0}
 \Bigl[\Lambda^*\bigl(\phi_h\zeta,K_h(\rho)\bigr)-\Lambda^*\bigl(\zeta,K_0(\rho)\bigr)\Bigr]
 =-\tfrac13\Rea\bigl(\lambda P_p(\zeta)\bigr).
\end{equation}
\end{lemma}

The proof uses neither \cite[Proposition~3.1]{GQW} nor its Beltrami flows.

\begin{proof}
Write \(\Omega=\operatorname{int}\zeta\), \(\Omega^*=\operatorname{ext}\zeta\),
\(f=f_\zeta:\D\to\Omega\) with \(f(0)=0\), \(f'(0)>0\), \(g:\Ds\to\Omega^*\) with
\(g(\infty)=\infty\), \(g'(\infty)>0\), \(F=f^{-1}\), \(G=g^{-1}\),
\(u(\xi)=-\lambda\xi^{1-p}\), and \(A_H=H''/H'\) for locally univalent \(H\). Since
\(\zeta\) is analytic, \(f\) and \(g\) extend conformally across \(S^1\)
\cite[Proposition~3.1]{PommerenkeBB}. So \(F,A_F\) are holomorphic near
\(\overline\Omega\), \(G,A_G\) near \(\overline{\Omega^*}\setminus\{\infty\}\), with
\(A_G(\xi)=O(|\xi|^{-3})\) and \(G'/G=\xi^{-1}+O(|\xi|^{-2})\) at \(\infty\), and
\(\zeta\) is a Weil--Petersson quasicircle, so \(I^L(\zeta)<\infty\)
\cite[Theorem~1.4]{WangInv}. The same holds for the analytic curves
\(\zeta_h=\phi_h\zeta\).

\emph{Step 1 (reduction to the Loewner energy).} For a domain
\(D\subset U_0\) containing \(\zeta\) put
\(X(D)=\Lambda^*(\zeta_h,\widehat{\mathbb C}\setminus\phi_hD)-\Lambda^*(\zeta,\widehat{\mathbb C}\setminus D)\).
If \(D_1\subset D_2\), then by \eqref{eq:LambdaNest} (twice) and conformal
invariance,
\(X(D_1)-X(D_2)=\mathcal B(\zeta_h,\phi_hD_2\setminus\phi_hD_1;\phi_hD_2)-\mathcal B(\zeta,D_2\setminus D_1;D_2)=0\).
Choose \(\delta/2<\rho_1<\min_\zeta|\xi|\), \(R>\max_\zeta|\xi|\), and
\(A=\{\rho_1<|\xi|<R\}\subset D_{\delta/2}\). The left side of \eqref{eq:LambdaIL}
is \(X(D_{\delta/2})=X(A)\). The curve \(\zeta\) separates the boundary
components of \(A\), and \(\phi_h\) is conformal on \(A\). By
\cite[Corollary~2.6]{CarfagniniWang},
\(X(A)=\mathcal W(\zeta_h,\widehat{\mathbb C}\setminus\phi_hA)-\mathcal W(\zeta,\widehat{\mathbb C}\setminus A)\),
where \(\mathcal W\) is Werner's measure normalized as the outer boundaries of
Brownian loops \cite[(2.4)]{CarfagniniWang}; and since \(I^L(\zeta)<\infty\),
Wang's restriction formula \cite[Theorem~4.1]{WangAIF}, restated as
\cite[Theorem~2.3]{CarfagniniWang} for an annular neighbourhood \(A\), gives
\(12\,[\mathcal W(\zeta,\widehat{\mathbb C}\setminus A)-\mathcal W(\zeta_h,\widehat{\mathbb C}\setminus\phi_hA)]
=I^L(\zeta_h)-I^L(\zeta)\). This is \eqref{eq:LambdaIL}.

\emph{Step 2 (the two Riemann maps).} By \cite[Theorem~1.4]{WangInv}
(\cite[Theorem~2.2]{CarfagniniWang}),
\[
 I^L(\zeta_h)=\frac1\pi\int_\D|A_{f_h}|^2\dd+\frac1\pi\int_{\Ds}|A_{g_h}|^2\dd
 +4\log f_h'(0)-4\log g_h'(\infty),
\]
where \(f_h=f_{\zeta_h}\) and \(g_h\) is the exterior map of \(\zeta_h\). Since
\(\phi_h\) is conformal on \(U_0\supset\overline{\Omega^*}\) and
\(\phi_h(\xi)=\xi+O(|\xi|^{1-p})\) at \(\infty\), \(g_h=\phi_h\circ g\) and
\(g_h'(\infty)=g'(\infty)\). For \(f_h\) we use the proof of
Proposition~\ref{prop:Wp-var} with \(\gamma=\zeta\): on \(\{r_*<|z|<1\}\),
\(f_h=\phi_h\circ f\circ\widetilde E_h^{-1}\), where \(\widetilde E_h\) is conformal on
an annulus \(\{r_*\le|z|\le1/r_*\}\),
\(\widetilde E_h^{-1}=\mathrm{id}+hw+O(h^2)\) uniformly on a smaller annulus, \(w\) is a
Laurent polynomial with \(w(z)/(iz)\in\mathbb R\) on \(S^1\), and
\(\delta f=f'Y=u\circ f+f'w\) on \(\{r_*<|z|<1\}\). Taking \(r_*\) closer to \(1\),
we may assume that \(f\) is conformal on \(\{|z|<1/r_*\}\) with
\(f(\{r_*\le|z|\le1/r_*\})\subset U_0\). Fix \(r_*<r_2<r_3<1\). For small \(|h|\)
the function \(\phi_h\circ f\circ\widetilde E_h^{-1}\) is holomorphic on
\(\{r_2<|z|<1/r_2\}\) and equals \(f_h\) on \(\{r_2<|z|<1\}\), so \(f_h\) extends
holomorphically to \(B(0,1/r_2)\); likewise \(\delta f\) extends by
\(u\circ f+f'w\). On \(\{|z|=1/r_3\}\) we have
\(f_h=f+h\,\delta f+O(h^2)\) uniformly, hence on \(\overline{B(0,1/r_3)}\) by the
maximum principle, and by Cauchy's estimates
\(f_h^{(j)}=f^{(j)}+h\,\delta f^{(j)}+O(h^2)\) uniformly on \(\overline\D\) for
\(j\le3\). Put \(V=\delta f\circ F\), holomorphic near \(\overline\Omega\), and
\(W=V-u\), holomorphic near \(\zeta\). At \(\xi=f(z)\in\zeta\),
\(W(\xi)=f'(z)w(z)=\frac{w(z)}{iz}\cdot izf'(z)\), and \(izf'(z)\) is tangent to
\(\zeta\); so \(W=T\tau\) on \(\zeta\), with \(T\) real and \(\tau\) the positive unit
tangent. Finally \(\frac{d}{dh}\big|_0\log f_h'(0)=\delta f'(0)/f'(0)=V'(0)\), which
is real because \(f_h'(0)>0\).

\emph{Step 3 (differentiation).} On \(\Ds\),
\(A_{g_h}=A_g+(\omega_h\circ g)g'\) with
\(\omega_h=\phi_h''/\phi_h'=hu''/(1-p\lambda h\xi^{-p})\), and
\(|A_g(z)|\le C|z|^{-3}\), \(|(\omega_h\circ g)g'|\le C|h||z|^{-p-1}\). By dominated
convergence, the change of variables \(\xi=g(z)\) and \(A_g=-(A_G\circ g)g'\),
\[
 \frac{d}{dh}\Big|_0\frac1\pi\int_{\Ds}|A_{g_h}|^2\dd
 =\frac2\pi\Rea\int_{\Ds}\overline{A_g}\,(u''\circ g)g'\dd
 =-\frac2\pi\Rea\int_{\Omega^*}\overline{A_G}\,u''\,d^2\xi
 =\frac1\pi\Rea\Bigl(i\oint_\zeta\overline{A_G}\,u'\,d\bar\xi\Bigr),
\]
by Stokes' formula \(\int_R\partial_\xi\Phi\,d^2\xi=\frac i2\oint_{\partial R}\Phi\,d\bar\xi\)
for \(\Phi=\overline{A_G}u'\) on \(\Omega^*\cap B(0,M)\), \(M\to\infty\)
(\(\Phi=O(|\xi|^{-p-3})\)), \(\zeta\) being oriented positively with respect to
\(\Omega\). Likewise, by Step~2, \(A_{f_h}=A_f+h(\delta f'/f')'+O(h^2)\) uniformly
on \(\overline\D\), and \(A_f=-(A_F\circ f)f'\), \(\delta f'/f'=V'\circ f\), so
\[
 \frac{d}{dh}\Big|_0\frac1\pi\int_\D|A_{f_h}|^2\dd
 =-\frac2\pi\Rea\int_\Omega\overline{A_F}\,V''\,d^2\xi
 =-\frac1\pi\Rea\Bigl(i\oint_\zeta\overline{A_F}\,V'\,d\bar\xi\Bigr).
\]
Hence
\begin{equation}\label{eq:ILprime}
 \frac{d}{dh}\Big|_0I^L(\zeta_h)
 =\frac1\pi\Rea\Bigl(i\oint_\zeta\bigl(\overline{A_G}\,u'-\overline{A_F}\,V'\bigr)d\bar\xi\Bigr)+4V'(0).
\end{equation}

\emph{Step 4 (boundary identities).} Parametrize \(\zeta\) by arclength \(s\),
positively, and let \(\tau=d\xi/ds\) and \(\tau_s=i\kappa_\zeta\tau\), with
curvature \(\kappa_\zeta\). Let \(H\in\{F,G\}\). Since \(|H|=1\) on \(\zeta\),
differentiating \(\overline{H}=1/H\) along \(\zeta\) gives
\(\overline{H'}\,\bar\tau=-H'\tau/H^2\), i.e.\ \(\overline{H'}=-H'\tau^2/H^2\);
differentiating the logarithm once more,
\begin{equation}\label{eq:conjA}
 \overline{A_H}\,d\bar\xi=A_H\,d\xi+2\,\frac{d\tau}{\tau}-2\,\frac{dH}{H}
 \qquad\text{on }\zeta .
\end{equation}
Insert \eqref{eq:conjA} into \eqref{eq:ILprime}. The integrals
\(\oint_\zeta A_Gu'\,d\xi\) and \(\oint_\zeta u'G'/G\,d\xi\) vanish, since the
integrands are holomorphic on \(\Omega^*\) and \(O(|\xi|^{-p-1})\) at \(\infty\);
\(\oint_\zeta A_FV'\,d\xi=0\), and \(\oint_\zeta V'F'/F\,d\xi=2\pi iV'(0)\), since
\(F\) has a simple zero at \(0\). Thus
\[
 \frac{d}{dh}\Big|_0I^L(\zeta_h)
 =\frac1\pi\Rea\Bigl(-2i\oint_\zeta W'\,\frac{d\tau}{\tau}\Bigr)-4V'(0)+4V'(0)
 =-\frac2\pi\Rea\Bigl(i\oint_\zeta W'\,\frac{d\tau}{\tau}\Bigr).
\]
Since \(W=T\tau\), \(W'=\tau^{-1}dW/ds=T_s+i\kappa_\zeta T\) and
\(d\tau/\tau=i\kappa_\zeta\,ds\), so
\(\Rea(iW'\,d\tau/\tau)=-\kappa_\zeta T_s\,ds\), and an integration by parts gives
\(\frac{d}{dh}|_0I^L(\zeta_h)=-\frac2\pi\oint_\zeta T\,\partial_s\kappa_\zeta\,ds\).

\emph{Step 5 (the Schwarzian).} Put \(q=A_F\tau\) on \(\zeta\). Since
\(dF/F=i|F'|\,ds\) on \(\zeta\), \eqref{eq:conjA} for \(H=F\) gives
\(\Ima q=|F'|-\kappa_\zeta\), and \(\Rea q=\partial_s\log|F'|\). From
\(q_s=A_F'\tau^2+i\kappa_\zeta q\) we get
\(\tau^2\mathcal S[F]=q_s-i\kappa_\zeta q-\frac12q^2\), hence
\[
 \Ima\bigl(\tau^2\mathcal S[F]\bigr)=\partial_s(|F'|-\kappa_\zeta)-\kappa_\zeta\Rea q-\Rea q\,\Ima q
 =\partial_s|F'|-\partial_s\kappa_\zeta-|F'|\Rea q=-\partial_s\kappa_\zeta .
\]
Therefore \(\Rea(i\oint_\zeta W\mathcal S[F]\,d\xi)=-\oint_\zeta T\Ima(\tau^2\mathcal S[F])\,ds
=\oint_\zeta T\,\partial_s\kappa_\zeta\,ds\), and, since \(\oint_\zeta V\mathcal S[F]\,d\xi=0\),
\[
 \frac{d}{dh}\Big|_0I^L(\zeta_h)=-\frac2\pi\Rea\Bigl(i\oint_\zeta(V-u)\mathcal S[F]\,d\xi\Bigr)
 =\frac2\pi\Rea\Bigl(-i\lambda\oint_{|\xi|=\delta/4}\xi^{1-p}\mathcal S[F](\xi)\,d\xi\Bigr)
 =4\Rea\bigl(\lambda P_p(\zeta)\bigr),
\]
because \(\mathcal S[F]\) is holomorphic on \(\Omega\supset\overline{B(0,\delta/2)}\) and
\(F=F_\zeta^{-1}\). With \eqref{eq:LambdaIL} this gives \eqref{eq:GQWanalytic}.
\end{proof}

We also checked \(\frac{d}{dh}I^L(\phi_h\zeta)=4\Rea(\lambda P_p(\zeta))\)
numerically (script \texttt{check\_IL.py}) for \(\zeta=f(S^1)\),
\(f(z)=z+(0.12+0.05i)z^2+(-0.04+0.03i)z^3\), \(p=2,3,4\), \(\lambda\in\{1,i\}\),
\(\phi_h(\xi)=\xi(1-p\lambda h\xi^{-p})^{1/p}\): both Riemann maps of the
star-shaped curve \(\phi_h\zeta\) are computed by Theodorsen's method on
\(N=1024\) nodes, \(I^L\) by \cite[Theorem~1.4]{WangInv} from the Fourier
coefficients of \(\log f_h'\), \(\log g_h'\), and \(\frac d{dh}\) by a central
difference with step \(10^{-4}\); \(P_p(\zeta)\) is exact. In all six cases the
relative deviation is below \(10^{-6}\) (largest for \(p=4\), \(\lambda=1\):
\(2.2257922\) against \(2.2257907\)), consistent with the \(O(h^2)\) error.

\begin{lemma}[First variation of \(\Lambda^*\) for Jordan curves]
\label{lem:GQWext}
For every Jordan curve \(\gamma\in\mathcal J_\delta\) the function
\(s\mapsto P_p(\phi_s\gamma)\) is continuous on \([-\varepsilon,\varepsilon]\), and
\begin{equation}\label{eq:GQWintegrated}
 A_t(\gamma)=-\frac13\int_0^t\Rea\bigl(\lambda P_p(\phi_s\gamma)\bigr)\,ds,
 \qquad |t|\le\varepsilon .
\end{equation}
In particular \(t\mapsto A_t(\gamma)\) is \(C^1\) with derivative
\(-\frac13\Rea(\lambda P_p(\phi_t\gamma))\) at every \(t\in(-\varepsilon,\varepsilon)\).
\end{lemma}

\begin{proof}
\emph{Step 1 (independence of the domain).} Let \(D_1\subset D_2\subset U_0\)
be domains containing \(\infty\), let \(\zeta\subset D_1\) be a Jordan curve,
and let \(|h|\le2\varepsilon\). Then
\[
 \Lambda^*\bigl(\phi_h\zeta,\widehat{\mathbb C}\setminus\phi_hD_1\bigr)
 -\Lambda^*\bigl(\zeta,\widehat{\mathbb C}\setminus D_1\bigr)
\]
is unchanged when \(D_1\) is replaced by \(D_2\). Indeed, by
\eqref{eq:LambdaNest} and the symmetry of \(\Lambda^*\), the difference
between the two expressions is
\(\mathcal B(\zeta,D_2\setminus D_1;D_2)-\mathcal B(\phi_h\zeta,\phi_hD_2\setminus\phi_hD_1;\phi_hD_2)\),
which vanishes by conformal invariance.

\emph{Step 2 (analytic curves).} Let \(\eta\in\mathcal J_\delta\) be analytic,
let \(|s|\le\varepsilon\), and put \(\zeta=\phi_s\eta\). This curve is analytic.
Moreover \(\operatorname{int}\zeta=K_s(\delta)\cup\phi_s(\operatorname{int}\eta\setminus K_0(\delta))\)
and \(K_s(\delta)\supset\overline{B(0,\delta/2)}\), so
\(\zeta\in\mathcal J_{\delta/2}\). By the flow rule,
\(A_{s+h}(\eta)-A_s(\eta)\) is the difference in Step~1 with
\(D_1=\phi_s(D_\delta)\). Step~1 with \(D_2=D_{\delta/2}\) and
\eqref{eq:GQWanalytic} with \(\rho=\delta/2\) show that \(s\mapsto A_s(\eta)\)
is differentiable everywhere (Lemma~\ref{lem:GQWanalytic} applies to \(\zeta\)), with derivative
\(-\frac13\Rea(\lambda P_p(\phi_s\eta))\). This derivative is continuous by
Step~4, so \eqref{eq:GQWintegrated} holds for \(\eta\).

\emph{Step 3 (continuity of \(\Lambda^*\) under interior approximation).} Let
\(\gamma\in\mathcal J_\delta\) and \(\Omega=\operatorname{int}\gamma\), and let
\(f:\D\to\Omega\) be conformal with \(f(0)=0\). By Carath\'eodory's theorem
\cite[Theorem~2.6]{PommerenkeBB}, \(f\) extends to a homeomorphism
\(\overline\D\to\overline\Omega\). Put \(\gamma_r=f(rS^1)\) and \(\gamma_1=\gamma\).
Choose \(r_1<1\) with \(\overline{B(0,\delta)}\subset f(r_1\D)\). Then
\(\gamma_r\in\mathcal J_\delta\) is analytic for \(r_1\le r<1\), and
\(\gamma_r\subset X:=\overline\Omega\setminus f(r_1\D)\), a compact subset of
\(U_0\). Fix \(|t|\le\varepsilon\), and write \(K=K_t(\delta)\),
\(\eta_r=\phi_t\gamma_r\) and \(\eta=\phi_t\gamma\). Since \(\phi_t\) and \(f\)
are uniformly continuous on \(X\) and \(\overline\D\),
\begin{equation}\label{eq:unifconv}
 \sup_\theta\bigl|\phi_t(f(re^{i\theta}))-\phi_t(f(e^{i\theta}))\bigr|\to0
 \qquad(r\uparrow1).
\end{equation}
For \(r<1\) we have \(K\subset\operatorname{int}\eta_r\subset\operatorname{int}\eta\).
By \eqref{eq:LambdaNest} with \(D=\operatorname{int}\eta\) and
\(D'=\operatorname{int}\eta_r\), using that a loop meets
\(\widehat{\mathbb C}\setminus\operatorname{int}\eta\) and \(K\) if and only if it
meets \(\eta\) and \(K\), and likewise for \(\eta_r\),
\[
 \Lambda^*(\eta_r,K)-\Lambda^*(\eta,K)=\mu(E_r),\qquad
 E_r=\{\ell\subset\operatorname{int}\eta:\ \ell\cap K\ne\emptyset,\
       \ell\cap\eta_r\ne\emptyset\}.
\]
Every loop \(\ell\subset\operatorname{int}\eta\) is compact, so
\(\operatorname{dist}(\ell,\eta)>0\), and by \eqref{eq:unifconv}
\(\ell\notin E_{r}\) for \(r\) close to \(1\); loops not contained in
\(\operatorname{int}\eta\) lie in no \(E_r\). Thus
\(\mathbf 1_{E_r}\to0\) pointwise as \(r\uparrow1\). For \(r_1\le r<1\), a loop in \(E_r\) joins
\(K\subset\operatorname{int}\eta_{r_1}\) to
\(\eta_r\subset\widehat{\mathbb C}\setminus\overline{\operatorname{int}\eta_{r_1}}\)
(when \(r>r_1\)), so it meets \(\eta_{r_1}\). Hence \(E_r\subset E_{r_1}\), and
\(\mu(E_{r_1})<\infty\) is a difference of two finite values of \(\Lambda^*\).
Dominated convergence gives \(\Lambda^*(\eta_r,K)\to\Lambda^*(\eta,K)\). Taking
\(t\) and \(0\) gives \(A_t(\gamma_r)\to A_t(\gamma)\).

\emph{Step 4 (continuity of \(P_p\)).} We claim that
\(\Pi(s,r)=P_p(\phi_s\gamma_r)\) is continuous on
\([-\varepsilon,\varepsilon]\times[r_1,1]\). Let \((s_n,r_n)\to(s,r)\). The
parametrizations \(J_n(\theta)=\phi_{s_n}(f(r_ne^{i\theta}))\) converge
uniformly to \(J(\theta)=\phi_s(f(re^{i\theta}))\), by joint uniform continuity
of \((s,z)\mapsto\phi_s(z)\) on \([-2\varepsilon,2\varepsilon]\times X\). Put
\(G_n=\operatorname{int}J_n\) and \(G=\operatorname{int}J\). All of these
contain \(\overline{B(0,\delta/2)}\). We show that \(G_n\to G\) in the sense of
kernels with respect to \(0\).
(i) A compact set \(L\subset G\) lies in a compact connected
\(L'\subset G\) containing \(0\). Once \(\|J_n-J\|_\infty<\operatorname{dist}(L',J)\),
the set \(L'\) misses \(J_n\) and contains \(0\in G_n\), so \(L'\subset G_n\).
(ii) Let \(H\ni0\) be the kernel of a subsequence. If \(H\not\subset G\), then,
since \(H\) is open and \(\operatorname{int}\overline G=G\), there is
\(w\in H\setminus\overline G\). Once \(\|J_n-J\|_\infty<\operatorname{dist}(w,J)\),
the straight-line homotopy shows that the winding number of \(J_n\) about
\(w\) is \(0\), so \(w\notin G_n\). This is a contradiction.
By the Carath\'eodory kernel theorem \cite[Theorem~3.1]{Duren}, the maps
\(F_{G_n}\) normalized by \(F(0)=0\), \(F'(0)>0\) converge to \(F_G\) locally
uniformly on \(\D\). So their Taylor coefficients converge, with
\(F_G'(0)>0\). By Lagrange inversion, so do the Taylor coefficients of
\(F_{G_n}^{-1}\) at \(0\). By polynomiality, \(\Pi(s_n,r_n)\to\Pi(s,r)\).
Since \(\Pi\) is continuous on a compact set, it is bounded (as also follows
from \cite[Lemma~3.7]{GQW}) and uniformly continuous. Hence
\(\sup_{|s|\le\varepsilon}|\Pi(s,r)-\Pi(s,1)|\to0\) as \(r\uparrow1\).

\emph{Step 5.} Apply Step~2 to \(\gamma_r\) and let \(r\uparrow1\). The left
side converges by Step~3 and the right side by Step~4, which proves
\eqref{eq:GQWintegrated}. The \(C^1\) statement follows from the
fundamental theorem of calculus and the continuity of
\(s\mapsto\Pi(s,1)\).
\end{proof}

\phantomsection\label{app:covariance}%
We now derive the score covariance from \eqref{eq:W}.
By Corollary~\ref{cor:coeffTm}(ii),
\[
 L_n\mathfrak t_m=\delta_{mn}C_n+(m+n)\,\mathfrak t_{m-n}\mathbf 1_{m>n}.
\]
First, \eqref{eq:W} extends from \(\mathcal C\) to every bounded Borel
\(F\) on \(\mathcal W\) whose orbits \(t\mapsto F(h\circ\psi^v_t)\) are
differentiable with derivative bounded uniformly in \(h\) and \(|t|\le1\),
in particular to bounded coefficient polynomials
(Corollary~\ref{cor:coeffBounded}). Indeed, by Lemma~\ref{lem:measureODE} the image laws
\(Q_t=(T_t)_*P\) have densities \(r_t\) with \(\|r_t-1-t\rho_v\|_\infty=O(t^2)\);
hence \(\int F\circ T_t\,dP=\int F\,r_t\,dP\), and differentiating at
\(t=0\), by dominated convergence on the left, gives
\(\int D_vF\,dP=\int F\rho_v\,dP\). For \(n\ge1\) and holomorphic
coefficient polynomials \(\varphi,\psi\), Lemma~\ref{lem:coeffResp} gives
\(\mathfrak d_n(\varphi\bar\psi)=i(L_n\varphi)\bar\psi-i\varphi\,\overline{L_{-n}\psi}\), and
\(\int\mathfrak d_nF\,dP=i\int F\bar{\mathfrak t}_n\,dP\) (Theorem~\ref{thm:Wproof},
Remark~\ref{rem:Wp-checks}(c)) then gives
\(\E[(L_n\varphi)\bar\psi]=\E[\varphi\,\overline{(L_{-n}\psi+\mathfrak t_n\psi)}]\),
with \(L_{-n}1=0\). Then \(\varphi=\psi=1\) gives \(\E\mathfrak t_n=0\),
\(\varphi=\mathfrak t_m\), \(\psi=1\), \(m\le n\), gives
\(\E[\mathfrak t_m\bar{\mathfrak t}_n]=\delta_{mn}C_n\), and rotation
covariance \(\mathfrak t_n(h\circ R_\alpha)=e^{in\alpha}\mathfrak t_n(h)\)
gives \(\E[\mathfrak t_n\mathfrak t_m]=0\). Since
\(s_n=\overline{\mathfrak t_n\circ\mathsf J}\) and \(\mathsf J_*P=P\),
this is \eqref{eq:allcov} for \(m,n\ge2\). We also checked the rule for
\(L_n\mathfrak t_m\) for \(m,n\le7\) symbolically
(Remark~\ref{rem:coeffRespChecks}).
\section{Deterministic lemmas for Section~\ref{sec:entropy}}\label{sec:entropyDet}

We first fix conventions. For a Jordan curve \(\gamma\subset\mathbb C\) with
\(0\in\Omega=\operatorname{int}\gamma\), \(\Omega^*=\operatorname{ext}\gamma\), a
\emph{pair} for \(\gamma\) consists of conformal maps \(f:\D\to\Omega\),
\(g:\Ds\to\Omega^*\) with \(f(0)=0\), \(g(\infty)=\infty\); both extend to
homeomorphisms of the closures \cite[Theorem~2.6]{PommerenkeBB}, and
\(y=g^{-1}\circ f|_{S^1}\in X\) is a \emph{welding of \(\gamma\) with interior
map \(f\)}. Rotating \(g\) composes \(y\) with a left rotation. A compact set
\(E\) is (conformally) removable if every homeomorphism of \(\widehat{\mathbb C}\)
that is conformal on \(\widehat{\mathbb C}\setminus E\) is M\"obius;
\(\mathcal W\) is the set of weldings of removable Jordan curves, as in
Section~\ref{sec:finiteinput}. Analytic Jordan curves are removable (by
Morera's theorem) and are quasicircles \cite[Chapter~5]{PommerenkeBB}.
By Lemma~\ref{lem:entR}(iv), for \(y\in\mathcal W\) the
normalized interior map \(f_y\) (\(f_y(0)=0\), \(f_y'(0)=1\)) is determined by
\(y\) and is unchanged under \(y\mapsto R_\alpha\circ y\); \(a_m(y)\),
\(\mathfrak t_n(y)\) and \(\rho_v(y)\) are its coefficients, stress coefficients
and stress score, bounded by de Branges' theorem \cite{deBranges}, and on
\(\mathcal W\) they are the Borel versions fixed in Section~\ref{sec:finiteinput}.
For a \(C^1\) orientation-preserving diffeomorphism \(\Psi(e^{i\theta})=e^{i\psi(\theta)}\)
(\(\psi'>0\)), its \emph{radial extension} is
\(\Lambda_\Psi(re^{i\theta})=re^{i\psi(\theta)}\), \(0\le r\le\infty\). It is a
homeomorphism of \(\widehat{\mathbb C}\) preserving \(\D\), \(\Ds\), \(0\), \(\infty\),
\(\Lambda_\Psi^{-1}=\Lambda_{\Psi^{-1}}\), and a direct computation gives
\begin{equation}\label{eq:entRadialBeltrami}
 \mu_{\Lambda_\Psi}(re^{i\theta})=e^{2i\theta}\,\frac{1-\psi'(\theta)}{1+\psi'(\theta)},
\end{equation}
continuous off \(\{0,\infty\}\) with \(\|\mu_{\Lambda_\Psi}\|_\infty<1\); so
\(\Lambda_\Psi\) is quasiconformal (points are removable for quasiconformal
maps \cite{LehtoVirtanen}). We use the measurable Riemann mapping theorem,
the composition rule for Beltrami coefficients, Weyl's lemma, and the
absolute continuity of quasiconformal maps
\cite{AhlforsLectures,AhlforsBers,LehtoVirtanen}. For a conformal \(F\) and a
qc map \(A\), \(\mu_{A\circ F^{-1}}(F(z))=\mu_A(z)F'(z)/\overline{F'(z)}\).

\begin{lemma}[Removability]\label{lem:entR}
Let \(\gamma\) be a removable Jordan curve.
\begin{enumerate}
\item \(\gamma\) has zero area.
\item A homeomorphism \(H\) of \(\widehat{\mathbb C}\) which is quasiconformal on
\(\widehat{\mathbb C}\setminus\gamma\) is quasiconformal on \(\widehat{\mathbb C}\).
\item For every quasiconformal \(E\) of \(\widehat{\mathbb C}\), \(E(\gamma)\) is removable.
\item If \((f_1,g_1)\) is a pair for \(\gamma\), \((f_2,g_2)\) a pair for a Jordan
curve \(\gamma_2\), and \(g_1^{-1}\circ f_1=g_2^{-1}\circ f_2\) on \(S^1\), then
\(f_2=M\circ f_1\), \(g_2=M\circ g_1\) with \(M(w)=cw\), \(c\ne0\).
\end{enumerate}
\end{lemma}
\begin{proof}
(i) If \(|\gamma|>0\), the normalized solution of the Beltrami equation with
coefficient \(\frac12\mathbf 1_\gamma\) is a homeomorphism conformal off \(\gamma\)
which is not M\"obius. (ii) Let \(\nu=\mu_{H^{-1}}\) on
\(\widehat{\mathbb C}\setminus H(\gamma)\), \(\nu=0\) on \(H(\gamma)\), and let \(Q\) solve
the Beltrami equation with coefficient \(\nu\). On \(\widehat{\mathbb C}\setminus\gamma\)
the maps \(Q\) and \(H^{-1}\) have equal coefficients on the image, so
\(Q\circ H\) is a homeomorphism conformal off \(\gamma\), hence M\"obius, and
\(H=Q^{-1}\circ(Q\circ H)\) is quasiconformal. (iii) If \(H\) is a homeomorphism
conformal off \(E(\gamma)\), then \(H\circ E\) is quasiconformal off \(\gamma\),
hence on \(\widehat{\mathbb C}\) by (ii), so \(H=(H\circ E)\circ E^{-1}\) is
quasiconformal. Its coefficient vanishes off the null set \(E(\gamma)\)
(by (i) and absolute continuity), so \(H\) is M\"obius by Weyl's lemma.
(iv) Put \(H=f_2\circ f_1^{-1}\) on \(\overline\Omega_1\) and
\(H=g_2\circ g_1^{-1}\) on \(\overline{\Omega_1^*}\). For \(\xi=f_1(z)=g_1(y(z))\in\gamma\)
both definitions give \(f_2(z)=g_2(y(z))\), so \(H\) is a homeomorphism of
\(\widehat{\mathbb C}\), conformal off \(\gamma\), hence M\"obius; it fixes \(0\) and
\(\infty\).
\end{proof}

\begin{lemma}[Right composition]\label{lem:entQ}
Let \(y\in\mathcal W\) have the pair \((f,g)\) for \(\gamma\), and let \(\Psi\in X\) have
a quasiconformal extension \(\Lambda:\overline\D\to\overline\D\) with \(\Lambda(0)=0\)
(e.g.\ \(\Lambda=\Lambda_\Psi\) for a \(C^1\) diffeomorphism). Let \(\sigma\) be the
coefficient of \(\Lambda^{-1}\circ f^{-1}\) on \(\Omega\), \(\sigma=0\) off \(\Omega\), and
\(E\) the principal solution of \(\partial_{\bar z}E=\sigma\,\partial_zE\). Then
\(E\circ f\circ\Lambda\) and \(E\circ g\) are conformal on \(\D\), \(\Ds\), they weld
\(y\circ\Psi\), and \(E(\gamma)\) is removable. Hence \(y\circ\Psi\in\mathcal W\) and
\(f_{y\circ\Psi}=(E\circ f\circ\Lambda-E(0))/(E\circ f\circ\Lambda)'(0)\). In particular
\(\mathcal W\) is invariant under right composition by \(C^1\) diffeomorphisms,
and every \(C^1\) diffeomorphism lies in \(\mathcal W\) (take \(y=\mathrm{id}\)).
\end{lemma}
\begin{proof}
\(E\) and \(\Lambda^{-1}\circ f^{-1}\) have the same coefficient on \(\Omega\), so
\(E\circ f\circ\Lambda=E\circ(\Lambda^{-1}\circ f^{-1})^{-1}\) is conformal on \(\D\);
\(E\) is conformal on \(\Omega^*\). On \(S^1\),
\((E\circ g)^{-1}\circ(E\circ f\circ\Lambda)=g^{-1}\circ f\circ\Psi\). The curve \(E(\gamma)\) is
removable by Lemma~\ref{lem:entR}(iii), and Lemma~\ref{lem:entR}(iv) gives the
formula for \(f_{y\circ\Psi}\).
\end{proof}

\begin{lemma}[Continuity]\label{lem:entC}
Let \(F_k,F\) be conformal maps of \(\D\) onto Jordan domains, normalized at \(0\),
with \(F_k\to F\) uniformly on \(\overline\D\), and with removable boundary curves.
Let \(y_k,y\) be weldings with interior maps \(F_k,F\), and let \(\Psi_k,\Psi\) be
\(C^1\) diffeomorphisms with lifts \(\psi_k\to\psi\) in \(C^1\). Then
\(a_m(y_k\circ\Psi_k)\to a_m(y\circ\Psi)\) for every \(m\).
\end{lemma}
\begin{proof}
Let \(\Lambda_k=\Lambda_{\Psi_k}\), \(\Lambda=\Lambda_\Psi\), and let \(\beta_k,\beta\) be the
coefficients of \(\Lambda_k^{-1},\Lambda^{-1}\). By \eqref{eq:entRadialBeltrami} applied to
\(\Psi_k^{-1}\), \(\|\beta_k\|_\infty\le k_0<1\) and \(\beta_k\to\beta\) uniformly on \(\D\),
\(\beta\) being continuous on \(\D\setminus\{0\}\). Let \(\sigma_k,\sigma\) be as in
Lemma~\ref{lem:entQ}: \(\sigma_k(F_k(z))=\beta_k(z)F_k'(z)/\overline{F_k'(z)}\) on
\(\Omega_k=F_k(\D)\), \(\sigma_k=0\) elsewhere, all supported in a fixed disc
\(B(0,R)\). We claim \(\sigma_k\to\sigma\) almost everywhere. For \(w\in\Omega\setminus\{0\}\),
Hurwitz's theorem gives \(F_k^{-1}(w)\to F^{-1}(w)\) and
\(F_k'(F_k^{-1}(w))\to F'(F^{-1}(w))\), so \(\sigma_k(w)\to\sigma(w)\). For
\(w\notin\overline\Omega\), once \(\|F_k-F\|_{S^1}<\operatorname{dist}(w,\gamma)\) the
straight-line homotopy shows that the winding number of \(F_k(S^1)\) about
\(w\) is \(0\), so \(\sigma_k(w)=0=\sigma(w)\). Finally \(\gamma=F(S^1)\) is null by
Lemma~\ref{lem:entR}(i).

As in the proof of Lemma~\ref{lem:Wp-AB}, the principal solutions are
\(E_k=z+\mathsf C[\omega_k]\), \(E=z+\mathsf C[\omega]\), with
\(\omega_k=(I-\sigma_kT)^{-1}\sigma_k\) in \(L^p\), \(p>2\), \(k_0\|T\|_p<1\). Then
\(\omega_k-\omega=(I-\sigma_kT)^{-1}[(\sigma_k-\sigma)(1+T\omega)]\to0\) in \(L^p\) by
dominated convergence on \(B(0,R)\), and \(E_k\to E\) uniformly on
\(\mathbb C\). The conformal maps \(\Xi_k=E_k\circ F_k\circ\Lambda_k\) of
Lemma~\ref{lem:entQ} converge to \(\Xi=E\circ F\circ\Lambda\) uniformly on
\(\overline\D\), since \(F_k\circ\Lambda_k\to F\circ\Lambda\) uniformly and \(E\) is uniformly
continuous on compacts. Their Taylor coefficients at \(0\) converge, with
\(\Xi'(0)\ne0\), and Lemma~\ref{lem:entQ} identifies
\((\Xi_k-\Xi_k(0))/\Xi_k'(0)\) with \(f_{y_k\circ\Psi_k}\).
\end{proof}

\begin{lemma}[Holomorphy]\label{lem:entH}
Let \(U\subset\mathbb C\) be open and let \(F_\zeta\), \(\zeta\in U\), be normalized
maps such that every \(\zeta_0\in U\) has a neighbourhood \(U_0\) and a radius
\(\varrho>1\) with \(F_\zeta\) univalent on \(\{|z|<\varrho\}\) for \(\zeta\in U_0\) and
\((\zeta,z)\mapsto F_\zeta(z)\) holomorphic on \(U_0\times\{|z|<\varrho\}\). Let \(y_\zeta\) be a
welding with interior map \(F_\zeta\), and \(\Psi\) a \(C^1\) diffeomorphism. Then
\(\zeta\mapsto a_m(y_\zeta\circ\Psi)\) is holomorphic on \(U\) for every \(m\).
\end{lemma}
\begin{proof}
Fix \(\zeta_0\), and let \(W_0\) be a quasiconformal map of \(\widehat{\mathbb C}\) with
\(W_0=F_{\zeta_0}\) on \(\overline\D\), \(W_0(\infty)=\infty\) (it exists since
\(F_{\zeta_0}(S^1)\) is an analytic curve, hence a quasicircle
\cite[Chapter~5]{PommerenkeBB}; composing with a quasiconformal self-map of the
exterior domain equal to the identity near the curve, we may assume
\(W_0(\infty)=\infty\)). Choose
\(1<\varrho'<\varrho\), compact neighbourhoods \(N_1\subset N_2^\circ\) of
\(F_{\zeta_0}(\overline\D)\) with \(N_2\subset F_{\zeta_0}(\{|z|<\varrho'\})\), and
\(\chi\in C_c^\infty(N_2^\circ)\) with \(\chi=1\) on \(N_1\). Put
\(T_\zeta=F_\zeta\circ F_{\zeta_0}^{-1}\) on \(F_{\zeta_0}(\{|z|<\varrho'\})\) and
\(\widetilde T_\zeta=\mathrm{id}+\chi(T_\zeta-\mathrm{id})\). For \(|\zeta-\zeta_0|<\delta\) small,
\(\widetilde T_\zeta\) is \(C^1\)-close to the identity, hence a diffeomorphism of
\(\mathbb C\), and its coefficient
\[
 \mu_{\widetilde T_\zeta}=\frac{(\partial_{\bar w}\chi)(T_\zeta-w)}
 {1+(\partial_w\chi)(T_\zeta-w)+\chi(T_\zeta'-1)}
\]
is holomorphic in \(\zeta\) at every point, with \(|\mu_{\widetilde T_\zeta}|\le\frac12\); no
conjugate of \(T_\zeta\) occurs because \(T_\zeta\) is holomorphic in \(w\). Let
\(W_\zeta=\widetilde T_\zeta\circ W_0\); then \(W_\zeta=F_\zeta\) on \(\overline\D\). With
\(\Lambda^{\rm e}=\Lambda_\Psi|_{\overline{\Ds}}\), let \(\nu_\zeta\) be the coefficient of
\(W_\zeta\circ\Lambda^{\rm e}\) on \(\Ds\) and \(\nu_\zeta=0\) on \(\D\). By the composition rule,
\(\nu_\zeta(z)\) is a M\"obius function of \(\mu_{\widetilde T_\zeta}(W_0(\Lambda^{\rm e}(z)))\) with
\(\zeta\)-independent coefficients, so it is holomorphic in \(\zeta\) for almost
every \(z\) and bounded by some \(k<1\); by Cauchy's formula in \(\zeta\),
\(\zeta\mapsto\nu_\zeta\in L^\infty\) is holomorphic on a smaller disc. Let \(E_\zeta\) be the
solution fixing \(0,1,\infty\). By the Ahlfors--Bers theorem
\cite{AhlforsBers}, \cite[Chapter~V]{AhlforsLectures}, \(E_\zeta(z)\) is holomorphic in
\(\zeta\) for each \(z\); it is holomorphic in \(z\in\D\), hence jointly holomorphic
by Hartogs' theorem. If \(g_\zeta\) is the exterior map of the pair, the map
\(E_\zeta\circ(W_\zeta\circ\Lambda^{\rm e})^{-1}\circ g_\zeta\) is conformal on \(\Ds\), and on \(S^1\)
\[
 \bigl(E_\zeta\circ(W_\zeta\circ\Lambda^{\rm e})^{-1}\circ g_\zeta\bigr)^{-1}\circ E_\zeta
 =g_\zeta^{-1}\circ F_\zeta\circ\Psi=y_\zeta\circ\Psi .
\]
Since \(E_\zeta(S^1)\) is a quasicircle, Lemma~\ref{lem:entR}(iv) gives
\(f_{y_\zeta\circ\Psi}=E_\zeta|_\D/E_\zeta'(0)\). The Taylor coefficients of
\(E_\zeta\) at \(0\) are holomorphic in \(\zeta\) (Cauchy integrals over \(|z|=\frac12\)),
and \(E_\zeta'(0)\ne0\).
\end{proof}

The next lemma is the first-variation formula of the Liouville--capacity
action along right composition, at analytic weldings. It is the
circle-field form of the first-variation formula of Takhtajan--Teo for
\(S_1\) \cite[Theorem~4.8]{TT} and of the corresponding formula for the
Velling--Kirillov potential \(K\) (cf.\ \cite[Theorem~5.3]{TTI}); compare
\cite[(1.13)]{BJweld}. Since several
sign conventions enter (left versus right action, complex versus real
parameters), we give a direct proof.

\begin{lemma}[First variation]\label{lem:entV}
Let \(\hat h\) be a real-analytic diffeomorphism of \(S^1\) and
\(v(\theta)=c_0+\sum_{n=1}^N\{x_n\cos(n\theta)+y_n\sin(n\theta)\}\). Then
\(t\mapsto H_\kappa(\hat h\circ\psi^v_t)\) is differentiable at \(t=0\), and
\[
 \frac{d}{dt}\Big|_{t=0}H_\kappa(\hat h\circ\psi^v_t)=\rho_v(\hat h).
\]
More precisely, with \(f=f_{\hat h}\) and \((x_n+iy_n)\) as above,
\[
 \frac{d}{dt}\Big|_{0}S_1(\hat h\circ\psi^v_t)=2\pi\Ima\sum_{n=2}^N(x_n+iy_n)[z^{n-2}]\mathcal Sf,
 \qquad
 \frac{d}{dt}\Big|_{0}2K(\hat h\circ\psi^v_t)=\Ima\sum_{n=1}^N(x_n+iy_n)\,k_n,
\]
where \(k_n=[z^{n-2}]\bigl(f'^2/f^2-z^{-2}\bigr)\). The formula for \(2K\) holds at every
\(\hat h\in\mathcal W\).
\end{lemma}

\begin{proof}
\emph{Step 0 (regularity).} By Lemma~\ref{lem:entQ}, \(\hat h\in\mathcal W\); let
\((f,g)\) be its normalized pair, \(\gamma=f(S^1)\). The lift of \(\hat h\) is real
analytic, so \(\hat h\) extends to a conformal map \(Y\) of an annulus around \(S^1\)
with \(Y(S^1)=S^1\), mapping the outer side into \(\Ds\). The function equal to
\(f\) on \(\overline\D\) and to \(g\circ Y\) just outside \(S^1\) is continuous and
holomorphic off \(S^1\), hence holomorphic, and it is injective near
\(\overline\D\). Thus \(f\) extends conformally to \(\{|z|<R\}\), likewise \(g\) to
\(\{|z|>1/R\}\), for some \(R>1\); \(\gamma\) is analytic, and \(F=f^{-1}\),
\(G=g^{-1}\) are conformal on neighbourhoods of \(\overline\Omega\) and
\(\overline{\Omega^*}\). Write \(N_h=h''/h'\).

\emph{Step 1 (representation).} Let \(L\) be the Laurent polynomial with
\(L(e^{i\theta})=ie^{i\theta}v(\theta)\), namely
\(L(z)=iz\{c_0+\sum_n(\frac{x_n-iy_n}2z^n+\frac{x_n+iy_n}2z^{-n})\}\). Fix
\(\chi\in C^\infty([0,\infty))\) with \(\chi=0\) on \([0,\frac14]\), \(\chi=1\) on
\([\frac12,\infty)\), and let \(\widehat\Psi_t\) be the flow of \(\widetilde V=\chi(|z|)L(z)\).
Since \(\widetilde V\) is tangent to \(S^1\), \(\widehat\Psi_t\) preserves \(\overline\D\), and
\(\widehat\Psi_t|_{S^1}=\psi^v_t\). For \(|t|\le t_0\): \(\widehat\Psi_t=\mathrm{id}\) on
\(\{|z|<\frac15\}\); \(\widehat\Psi_t\) is the holomorphic flow of \(L\) on an annulus
\(A=\{r_1<|z|<R_1\}\), \(\frac12<r_1<1<R_1<R\), with
\(\widehat\Psi_t=\mathrm{id}+tL+O(t^2)\) uniformly there; and the coefficient \(\mu_t\) of
\(\widehat\Psi_t^{-1}=\widehat\Psi_{-t}\) is supported in \(\{\frac15\le|z|\le r_1'\}\),
\(r_1'<r_1\), with \(\mu_t=-t\,\partial_{\bar z}\widetilde V+O(t^2)\) uniformly. Let \(\sigma_t\)
be the coefficient of \(\widehat\Psi_t^{-1}\circ F\) on \(\Omega\) (zero elsewhere), supported
in the compact \(\mathcal K=f(\{\frac15\le|z|\le r_1'\})\), and \(E_t\) the principal
solution. Since \(\|\sigma_t-t\dot\sigma\|_\infty=O(t^2)\), \(\dot\sigma(f(z))=-\partial_{\bar z}\widetilde V(z)f'(z)/\overline{f'(z)}\),
Lemma~\ref{lem:Wp-AB} gives \(E_t=\mathrm{id}+t\phi+O(t^2)\) uniformly on \(\mathbb C\), with
\(\phi=\mathsf C[\dot\sigma]\) holomorphic off \(\mathcal K\), \(\phi(\xi)=O(1/\xi)\) at \(\infty\);
as \(E_t\) is holomorphic off \(\mathcal K\), Cauchy estimates give the same
expansion for all derivatives on compact subsets of \(\mathbb C\setminus\mathcal K\), and
\((E_t-\mathrm{id}-t\phi)^{(j)}=O(t^2|\xi|^{-1-j})\) near \(\infty\). By
Lemma~\ref{lem:entQ} and the normalization \(A_t(\xi)=(\xi-E_t(0))/E_t'(0)\), the
normalized pair of \(\hat h\circ\psi^v_t\) is
\[
 f_t=A_t\circ E_t\circ f\circ\widehat\Psi_t,\qquad g_t=A_t\circ E_t\circ g,
\]
(\(f_t'(0)=1\) since \(\widehat\Psi_t=\mathrm{id}\) near \(0\)). Put \(\Theta_t=f_t\circ F\) on \(\Omega\). On
the neighbourhood \(f(A)\) of \(\gamma\), \(\Theta_t=A_t\circ E_t\circ(f\circ\widehat\Psi_t\circ F)\) is a
composition of holomorphic maps, so \(\Theta_t\) extends holomorphically to
\(\Omega\cup f(A)\), and there
\begin{equation}\label{eq:entThetaVel}
 \Theta_t=\mathrm{id}+t\widetilde\phi+O(t^2),\qquad
 \widetilde\phi=\phi+\alpha+f_*L,\quad \alpha(\xi)=-\phi(0)-\phi'(0)\xi,\quad
 f_*L=(f'L)\circ F .
\end{equation}
By the maximum principle on \(f(\{|z|<\varrho\})\), \(1<\varrho<R_1\), the expansion
holds uniformly on this domain with \(\widetilde\phi\) holomorphic there, and by
Cauchy estimates in \(C^2\) on \(\overline\Omega\). Since \(\Theta_t(0)=0\),
\(\Theta_t'(0)=1\), we have \(\widetilde\phi(0)=\widetilde\phi'(0)=0\).

\emph{Step 2 (the capacity).} \(K(\hat h\circ\psi_t^v)=\log|g_t'(\infty)|=\log|g'(\infty)|-\log|E_t'(0)|\),
so its derivative is \(-\Rea\phi'(0)\). Differentiating the Cauchy transform
at \(0\notin\mathcal K\) and changing variables \(\xi=f(z)\),
\[
 \phi'(0)=-\frac1\pi\int\frac{\dot\sigma(\xi)}{\xi^2}\,d^2\xi
 =\frac1\pi\int_\D\partial_{\bar z}\widetilde V\,\frac{f'^2}{f^2}\,\dd .
\]
Let \(q=f'^2/f^2-z^{-2}=\sum_{n\ge1}k_nz^{n-2}\), holomorphic on
\(\{0<|z|<1\}\). Since \(\widetilde V\) vanishes near \(0\) and equals \(L\) on \(\{|z|\ge\frac12\}\),
Stokes' formula \(\int_{|z|<r}\partial_{\bar z}u\,\dd=\frac1{2i}\oint_{|z|=r}u\,dz\),
\(\frac12<r<1\), and the residue theorem give
\[
 \int_\D\partial_{\bar z}\widetilde V\,z^{-2}\dd=\frac{2\pi i}{2i}[z^1]L=i\pi c_0\in i\mathbb R,\qquad
 \int_\D\partial_{\bar z}\widetilde V\,q\,\dd=\frac{2\pi i}{2i}\,[z^{-1}](Lq)
 =\frac{i\pi}2\varkappa,
\]
where \(\varkappa=\sum_nk_n(x_n+iy_n)\): only the terms
\(\frac i2(x_n+iy_n)z^{1-n}\) of \(L\) contribute to \([z^{-1}](Lq)\). Hence
\(-\Rea\phi'(0)=-\frac1\pi\Rea\frac{i\pi\varkappa}2=\frac12\Ima\varkappa\), so the
derivative of \(2K\) is \(\Ima\varkappa\). Step~1 up to the formula for
\((f_t,g_t)\), and Step~2, use only \(\hat h\in\mathcal W\) (Lemma~\ref{lem:entQ}); so the
formula for \(2K\) holds at every \(\hat h\in\mathcal W\).

\emph{Step 3 (reduction of \(S_1\) to the curve).} Since \(f_t'(0)=1\),
\(S_1(\hat h\circ\psi^v_t)=\int_\D|N_{f_t}|^2\dd+\int_{\Ds}|N_{g_t}|^2\dd-4\pi\log|g_t'(\infty)|\).
Now \(N_{f_t}=(N_{\Theta_t}\circ f)f'+N_f\) and \(N_{g_t}=(N_{E_t}\circ g)g'+N_g\), with
\(N_{\Theta_t}=t\widetilde\phi''+O(t^2)\) uniformly on \(\overline\Omega\) and
\(N_{E_t}=t\phi''+O(t^2(1+|\xi|)^{-3})\) uniformly on \(\overline{\Omega^*}\); since \(f'\),
\(g'\) are bounded and \(N_g=O(|z|^{-3})\), we may differentiate under the
integrals. Changing variables to \(\xi=f(z)\) and \(\xi=g(z)\), with
\(N_f=-(N_F\circ f)f'\), \(N_g=-(N_G\circ g)g'\),
\[
 \frac{d}{dt}\Big|_0S_1(\hat h\circ\psi^v_t)
 =-2\Rea\int_\Omega\overline{N_F}\,\widetilde\phi''\,d^2\xi
  -2\Rea\int_{\Omega^*}\overline{N_G}\,\phi''\,d^2\xi+4\pi\Rea\phi'(0).
\]
Since \(\overline{N_F}\,\widetilde\phi''=\partial_\xi(\overline{N_F}\,\widetilde\phi')\), and likewise on
\(\Omega^*\) (where the integrand is \(O(|\xi|^{-5})\)), Stokes' formula
\(\int_U\partial_\xi u\,d^2\xi=-\frac1{2i}\oint_{\partial U}u\,d\bar\xi\) gives, with \(\gamma\)
positively oriented,
\begin{equation}\label{eq:entS1curve}
 \frac{d}{dt}\Big|_0S_1(\hat h\circ\psi^v_t)
 =\Ima\oint_\gamma\bigl(\overline{N_F}\,\widetilde\phi'-\overline{N_G}\,\phi'\bigr)\,d\bar\xi
 +4\pi\Rea\phi'(0).
\end{equation}

\emph{Step 4 (reflection on \(\gamma\)).} Parametrize \(\gamma\) by arclength, so that
\(\xi_s=e^{i\tau}\) with \(\tau\) the tangent angle. Since \(|F|=1\) on \(\gamma\),
\(\overline{F'(\xi)\xi_s}=-F'(\xi)\xi_s/F(\xi)^2\); taking the logarithmic derivative in
\(s\) gives \(\overline{N_F\xi_s}-i\tau_s=(N_F-2F'/F)\xi_s+i\tau_s\), that is,
\[
 \overline{N_F}\,d\bar\xi=(N_F-2F'/F)\,d\xi+2i\,d\tau\quad\text{on }\gamma,
\]
and the same identity holds for \(G\). The resulting holomorphic integrals
vanish: \((N_F-2F'/F)\widetilde\phi'\) is holomorphic near \(\overline\Omega\) except for a
simple pole at \(0\) with residue \(-2\widetilde\phi'(0)=0\), and
\((N_G-2G'/G)\phi'\) is holomorphic near \(\overline{\Omega^*}\cap\mathbb C\) and
\(O(|\xi|^{-3})\) at \(\infty\). Hence \eqref{eq:entS1curve} becomes
\(2\Rea\oint_\gamma(\widetilde\phi'-\phi')\,d\tau+4\pi\Rea\phi'(0)\). By
\eqref{eq:entThetaVel}, \(\widetilde\phi'-\phi'=-\phi'(0)+(f_*L)'\) near \(\gamma\), and
\(\oint_\gamma d\tau=2\pi\) (see below), so the terms in \(\phi'(0)\) cancel:
\[
 \frac{d}{dt}\Big|_0S_1(\hat h\circ\psi^v_t)=2\Rea\oint_\gamma(f_*L)'\,d\tau .
\]

\emph{Step 5 (back to the circle).} On \(\gamma\), \(\xi=f(e^{i\theta})\),
\((f_*L)'(\xi)=(N_fL+L')(z)\) and \(\tau=\arg(izf'(z))\), so
\(d\tau=(1+\Rea w)\,d\theta\) with \(w=zN_f(z)\) (and
\(\int\Rea w\,d\theta=\Rea\frac1i\oint N_f\,dz=0\), which gives \(\oint d\tau=2\pi\)).
On \(S^1\), \(L=izv\) and \(L'=iv+v'\), so \(N_fL+L'=v'+iv(1+w)\) and
\[
 \frac{d}{dt}\Big|_0S_1(\hat h\circ\psi^v_t)
 =2\int_0^{2\pi}(1+\Rea w)(v'-v\Ima w)\,d\theta
 =-2\int_0^{2\pi}v\,\bigl\{\partial_\theta\Rea w+(1+\Rea w)\Ima w\bigr\}\,d\theta ,
\]
after integrating \(v'\) by parts. On the other hand
\(\partial_\theta N_f(e^{i\theta})=izN_f'\) gives \(z^2N_f'=-i\partial_\theta w-w\), hence
\(z^2\mathcal Sf=-i\partial_\theta w-w-\frac12w^2\) and
\(\Ima(z^2\mathcal Sf)=-\partial_\theta\Rea w-(1+\Rea w)\Ima w\). Therefore
\[
 \frac{d}{dt}\Big|_0S_1(\hat h\circ\psi^v_t)=2\int_0^{2\pi}v\,\Ima\bigl(z^2\mathcal Sf\bigr)\,d\theta
 =2\pi\Ima\sum_{n\ge2}(x_n+iy_n)[z^{n-2}]\mathcal Sf ,
\]
using \(\int_0^{2\pi}e^{in\theta}\cos(m\theta)\,d\theta=\pi\delta_{nm}\),
\(\int_0^{2\pi}e^{in\theta}\sin(m\theta)\,d\theta=i\pi\delta_{nm}\) and
\(\int_0^{2\pi}e^{in\theta}d\theta=0\) for \(n,m\ge1\).

\emph{Step 6.} Combining,
\(\frac d{dt}|_0H_\kappa(\hat h\circ\psi^v_t)=\Ima\sum_n(x_n+iy_n)\{\frac{\cL}{12}[z^{n-2}]\mathcal Sf+k_n\}
=\Ima\sum_n(x_n+iy_n)\mathfrak t_n=\rho_v(\hat h)\), since
\(\mathfrak t_n=[z^{n-2}]T_f\) with \(T_f\) as in \eqref{eq:stress} and
\(\Ima((x+iy)\mathfrak t)=x\Ima\mathfrak t+y\Rea\mathfrak t\).
\end{proof}

\begin{remark}[Checks of Lemma~\ref{lem:entV}]\label{rem:entVchecks}
(a) Along the M\"obius flow \(m_{\tanh s}\), \(m_x(z)=(z+x)/(1+xz)\), the field is
\(v=-2\sin\theta\), the pair is \((A\circ m_x,A)\) with \(A(w)=(w-x)/(1-x^2)\), so
\(a_1=-x\), \(k_1=2a_1\) and the lemma gives \(\frac{d}{ds}K=-2\Rea a_1=2\tanh s\), which is
the derivative of the exact value \(K(m_x)=-\log(1-x^2)\) of
Corollary~\ref{cor:entropyConsequences}(ii); and \(S_1\) is constant, as
\(\mathcal Sf=0\). (b) At \(\psi^{\cos n\theta}_s\) the interior map is
\(z+\frac{is}2z^{n+1}+O(s^2)\) (Proposition~\ref{prop:potential}), and the lemma gives
\(\frac d{ds}S_1=\pi(n^3-n)s+O(s^2)\), \(\frac d{ds}K=\frac{n}2s+O(s^2)\), the
derivatives of the expansions \(S_1=\frac\pi2(n^3-n)s^2\), \(K=\frac n4s^2\) of
Proposition~\ref{prop:potential}.
\end{remark}

\begin{lemma}[Density of finite-mode words]\label{lem:entD}
Let \(\mathcal G_{\rm fin}\) be the group generated by the flows \(\psi^v_t\) of real
trigonometric polynomials \(v\). For every orientation-preserving
\(C^\infty\) diffeomorphism \(\Phi\) of \(S^1\) there are \(\Phi_j\in\mathcal G_{\rm fin}\) with
\(\Phi_j\to\Phi\) in \(C^\infty\).
\end{lemma}
\begin{proof}
A smooth diffeomorphism need not be the time-one map of an autonomous
field, so we use an isotopy. Let \(\widetilde\Phi\) be a lift, and
\(\widetilde\Phi_\tau=(1-\tau)\mathrm{id}+\tau\widetilde\Phi\), \(\tau\in[0,1]\); each \(\widetilde\Phi_\tau\) is the lift of a
smooth diffeomorphism (\(\widetilde\Phi_\tau'\ge\min(1,\min\widetilde\Phi')>0\)), and
\(\partial_\tau\widetilde\Phi_\tau=u(\tau,\widetilde\Phi_\tau)\) with the smooth time-dependent field
\(u(\tau,\cdot)=(\widetilde\Phi-\mathrm{id})\circ\widetilde\Phi_\tau^{-1}\), \(2\pi\)-periodic. Let \(\sigma_N\) be the
Fejér mean of order \(N\) and \(u_N(\tau,\cdot)=\sigma_Nu(\lfloor N\tau\rfloor/N,\cdot)\). Since
\(\partial^j_\theta\sigma_N=\sigma_N\partial^j_\theta\), Fejér means are sup-norm contractions converging
uniformly on equicontinuous families, and \(\tau\mapsto u(\tau)\) is continuous into
every \(C^k\), we get \(\sup_\tau\|u_N(\tau)-u(\tau)\|_{C^k}\to0\) for every \(k\). The
time-one map of \(\dot\theta=u_N(\tau,\theta)\) is
\(\psi^{v_{N-1}}_{1/N}\circ\cdots\circ\psi^{v_0}_{1/N}\in\mathcal G_{\rm fin}\),
\(v_j=\sigma_Nu(j/N)\). Gronwall's inequality applied to the equation and to its
variational equations up to order \(k\) (whose coefficients involve
\(\partial^i_\theta u\), \(i\le k\), and which require \(u\in C^{k+1}\)) shows that these
time-one maps converge to \(\Phi\) in \(C^k\). A diagonal choice gives
\(C^\infty\) convergence.
\end{proof}
\section{Coefficient responses}\label{app:coeffResp}

We prove the coefficient responses used in Sections~\ref{sec:finiteinput}
and~\ref{app:Wproof} and Appendix~\ref{app:covariance}. The proof uses only
Lemma~\ref{lem:entQ}, Lemma~\ref{lem:Wp-AB}, Koebe's distortion theorem and de
Branges' theorem, and treats a general real trigonometric field at every point
of \(\mathcal W\), without approximation by smooth weldings. It is the interior
counterpart of Proposition~\ref{prop:Wp-var}; the two agree by
Remark~\ref{rem:Wp-checks}(d).

The notation is as follows. Here \(f=z(1+\sum_{m\ge1}a_mz^m)\) is a normalized univalent
map of \(\D\), \(u=f/z\), and \([z^k]\) is a Laurent coefficient at \(0\). For
\(n\ge0\) and \(0\le k\le n\) put
\[
 \mathfrak m_{n,k}=[w^{n-2}]\bigl(f'(w)^2f(w)^{k-2}\bigr)=[w^{n-k}]\bigl(f'^2u^{k-2}\bigr)
 \in\mathbb C[a_1,\ldots,a_{n-k}],
\]
so \(\mathfrak m_{n,n}=1\) and \(\mathfrak m_{n,k}=M_{n,k}\) of \eqref{eq:Wp-E} for \(k\ge1\). A real
trigonometric field of degree \(\le N\) is written
\[
 v(\theta)=c_0+\sum_{n=1}^N\{x_n\cos(n\theta)+y_n\sin(n\theta)\}
 =\sum_{|j|\le N}\hat v_je^{ij\theta},\qquad
 \hat v_0=c_0,\ \hat v_{\pm n}=\tfrac12(x_n\mp iy_n).
\]
\(\psi^v_t\) is its circle flow (\(\dot\theta=v(\theta)\)),
\(D_vF(h)=\frac d{dt}|_{t=0}F(h\circ\psi^v_t)\), and, once \(D_vF\) is known to be real-linear
in \(v\), \(\mathfrak d_j=D_{e^{ij\theta}}\) is its complex-linear extension, as in
Section~\ref{app:Wproof}; coefficient polynomials are evaluated at \(h\in\mathcal W\) through \(f_h\).

\begin{lemma}[Coefficient responses]\label{lem:coeffResp}
\begin{enumerate}
\item For \(j\in\mathbb Z\), \(z\in\D\) and \(|z|<r<1\), the function
\[
 \mathbf L_jf(z)=\frac{f(z)^2}{2\pi i}\oint_{|w|=r}
 \frac{w^{1+j}f'(w)^2}{f(w)^2\,(f(w)-f(z))}\,dw
\]
does not depend on \(r\). It is holomorphic on \(\D\), and \(\mathbf L_jf=O(z^2)\). Moreover
\[
 \mathbf L_nf=z^{n+1}f'\ (n\ge1),\qquad \mathbf L_0f=zf'-f,\qquad
 \mathbf L_{-n}f=z^{1-n}f'-\sum_{k=0}^n\mathfrak m_{n,k}\,f^{1-k}\ (n\ge1).
\]
\item Let \(L_j\) be the derivation of \(\mathbb C[a_1,a_2,\ldots]\) with
\(L_ja_m=[z^{m+1}]\mathbf L_jf\). Then, for \(m,n\ge1\),
\[
 L_na_m=\delta_{mn}+(m-n+1)\,a_{m-n}\mathbf 1_{m>n},\qquad L_0a_m=ma_m,
\]
\[
 L_{-n}a_m=(m+n+1)\,a_{m+n}-\sum_{k=0}^n\mathfrak m_{n,k}\,[z^{m+k}]u^{1-k}
 \ \in\mathbb C[a_1,\ldots,a_{m+n}].
\]
Each \(L_ja_m\) is weighted homogeneous of weight \(m-j\), where \(a_k\) has weight \(k\).
\item Let \(h\in\mathcal W\), \(f=f_h\), and let \(v\) be a real trigonometric field of degree
\(\le N\). Then there is a constant \(C_v\), independent of \(h\), with
\begin{equation}\label{eq:coeffResp-f}
 \sup_{|z|\le1/20}\Bigl|f_{h\circ\psi^v_t}(z)-f(z)-t\,i\sum_{|j|\le N}\hat v_j\,\mathbf L_jf(z)\Bigr|
 \le C_vt^2\qquad(|t|\le1).
\end{equation}
Consequently, for every \(m\ge1\), \(t\mapsto a_m(h\circ\psi^v_t)\) is differentiable at
\(t=0\) and
\begin{equation}\label{eq:coeffResp-a}
 D_va_m(h)=i\sum_{|j|\le N}\hat v_j\,(L_ja_m)(h),
\end{equation}
which is real-linear in \(v\). In particular \(\mathfrak d_ja_m=iL_ja_m\) on \(\mathcal W\) for every
\(j\in\mathbb Z\):
\[
 D_{\cos n\theta}=\tfrac i2(L_n+L_{-n}),\qquad D_{\sin n\theta}=\tfrac12(L_n-L_{-n}),\qquad
 D_1=iL_0
\]
on \(\mathbb C[a_1,a_2,\ldots]\). For each \(m\), \(D_va_m\) is a polynomial in
\(a_1,\ldots,a_{m+N}\).
\end{enumerate}
\end{lemma}

\begin{proof}
\emph{(i)} Since \(f\) is univalent, for fixed \(z\) the integrand is meromorphic on
\(\{|w|<1\}\). Its poles lie at most at \(w=z\) and \(w=0\). Hence the integral does
not depend on \(r\in(|z|,1)\), and it is holomorphic in \(z\) on \(\{|z|<r\}\), with
the factor \(f(z)^2\) giving \(O(z^2)\). At \(w=z\) the pole is simple
(\(f'(z)\ne0\)). Its residue is \(z^{1+j}f'(z)/f(z)^2\), which contributes
\(z^{1+j}f'(z)\). For the residue at \(w=0\), integrate over \(|w|=\varepsilon\), where
\(\varepsilon\) is so small that \(|f(w)|<|f(z)|\) there. Then
\((f(w)-f(z))^{-1}=-\sum_{k\ge0}f(w)^kf(z)^{-k-1}\) uniformly on the circle, and
\[
 \operatorname{Res}_{w=0}\,w^{1+j}f'(w)^2f(w)^{k-2}=[w^{-2-j}]\bigl(f'^2f^{k-2}\bigr).
\]
This vanishes unless \(j\le0\) and \(k\le-j\), since \(f'^2f^{k-2}=O(w^{k-2})\). For
\(j=-n\le0\) it equals \(\mathfrak m_{n,k}\), and \(\mathfrak m_{0,0}=1\). Collecting terms gives the three
formulas.

\emph{(ii)} Write \(f'=1+\sum_m(m+1)a_mz^m\). Then
\([z^{m+1}](z^{n+1}f')=\delta_{mn}+(m-n+1)a_{m-n}\mathbf 1_{m>n}\),
\([z^{m+1}](zf'-f)=ma_m\), and \([z^{m+1}](z^{1-n}f')=(m+n+1)a_{m+n}\). Also
\([z^{m+1}]f^{1-k}=[z^{m+k}]u^{1-k}\in\mathbb C[a_1,\ldots,a_{m+k}]\), because \(u\) is a unit with
constant term \(1\), so every power \(u^r\), \(r\in\mathbb Z\), has coefficients in
\(\mathbb C[a_1,a_2,\ldots]\). Since \(k\le n\) and \(\mathfrak m_{n,k}\in\mathbb C[a_1,\ldots,a_n]\), the
formula for \(L_{-n}a_m\) follows. For weights, replace \(f\) by \(f_\beta(z)=e^{-i\beta}f(e^{i\beta}z)\),
which multiplies \(a_k\) by \(e^{ik\beta}\). The substitution \(\omega=e^{i\beta}w\) in the integral of (i)
gives \(\mathbf L_jf_\beta(z)=e^{-i(1+j)\beta}\,\mathbf L_jf(e^{i\beta}z)\). Taking
\([z^{m+1}]\) gives \(L_ja_m(f_\beta)=e^{i(m-j)\beta}L_ja_m(f)\). A polynomial with this
covariance for all \(\beta\) is a combination of monomials of weight \(m-j\).

\emph{(iii), Step 1: an extension of the flow.} Let
\(\mathcal L(z)=iz\sum_j\hat v_jz^j\), a Laurent polynomial with
\(\mathcal L(e^{i\theta})=ie^{i\theta}v(\theta)\). Fix \(\chi\in C^\infty([0,\infty))\) with \(\chi=0\) on
\([0,\frac14]\) and \(\chi=1\) on \([\frac12,\infty)\), and put \(\widetilde V(z)=\chi(|z|)\mathcal L(z)\). This
field is smooth on \(\overline\D\) and tangent to \(S^1\), so its flow \(\widehat\Psi_t\) is defined for all
\(t\). Each \(\widehat\Psi_t\) is a \(C^\infty\) diffeomorphism of \(\overline\D\) preserving \(S^1\), smooth in
\((t,z)\), with \(\widehat\Psi_t|_{S^1}=\psi^v_t\). It is the identity on \(\{|z|\le\frac14\}\), where
\(\widetilde V=0\). Let \(c_v=\sup_{\overline\D}|\widetilde V|\) and \(t_0=\min(1,1/(10c_v))\). For
\(|t|\le t_0\), a trajectory starting in \(\{|z|\ge\frac35\}\) stays in \(\{|z|>\frac12\}\) up to time \(|t|\). There
\(\widetilde V=\mathcal L\) is holomorphic, and \(\partial_{\bar z}\widehat\Psi_t\) solves a homogeneous linear equation
along it, so \(\widehat\Psi_t\) is holomorphic on \(\{\frac35<|z|<1\}\). Hence the Beltrami coefficient
\(\mu_t\) of \(\widehat\Psi_{-t}=\widehat\Psi_t^{-1}\) is supported in
\(\overline A=\{\frac14\le|z|\le\frac35\}\). Since \(\widehat\Psi_s=\mathrm{id}+s\widetilde V+O(s^2)\) in \(C^1(\overline\D)\),
\[
 \|\mu_t+t\,\partial_{\bar z}\widetilde V\|_\infty\le M_vt^2,\qquad
 \|\mu_t\|_\infty\le\tfrac12\qquad(|t|\le t_0),
\]
after decreasing \(t_0\). Here \(t_0\) and \(M_v\) depend only on \(v\) and \(\chi\).

\emph{Step 2: the deformed map.} Let \((f,g)\) be a pair for \(h\) with \(f=f_h\), and let
\(\Omega=f(\D)\). Apply Lemma~\ref{lem:entQ} with \(\Psi=\psi^v_t\) and \(\Lambda=\widehat\Psi_t\). The
coefficient of \(\widehat\Psi_{-t}\circ f^{-1}\) on \(\Omega\) is
\(\sigma_t(f(z))=\mu_t(z)f'(z)/\overline{f'(z)}\), and it vanishes off the compact set
\(\mathcal K=f(\overline A)\). Let \(E_t\) be the principal solution with coefficient \(\sigma_t\). Then
\(E_t\) is conformal off \(\mathcal K\), in particular near \(0\), and
\[
 f_{h\circ\psi^v_t}=\frac{E_t\circ f\circ\widehat\Psi_t-E_t(0)}{E_t'(0)}
 =\frac{E_t\circ f-E_t(0)}{E_t'(0)}\quad\text{on }\{|z|<\tfrac14\}.
\]
Since \(\widehat\Psi_t=\mathrm{id}\) near \(0\), we have \((E_t\circ f\circ\widehat\Psi_t)'(0)=E_t'(0)\).

\emph{Step 3: first order.} Put \(\dot\sigma(f(z))=-\partial_{\bar z}\widetilde V(z)f'(z)/\overline{f'(z)}\) on \(\mathcal K\),
\(\dot\sigma=0\) elsewhere, and \(\phi=\mathsf C[\dot\sigma]\). Then
\(\|\sigma_t\|_\infty\le\frac12\) and \(\|\sigma_t-t\dot\sigma\|_\infty=\|\mu_t+t\partial_{\bar z}\widetilde V\|_\infty\le M_vt^2\).
By Koebe's distortion theorem, \(\mathcal K\subset\{\frac4{25}\le|\xi|\le4\}\) and
\(f(\{|z|\le\frac1{20}\})\subset\{|\xi|\le\frac1{18}\}\), for every \(h\). By
Lemma~\ref{lem:Wp-AB},
\(\sup_{\mathbb C}|E_t-\mathrm{id}-t\phi|\le C t^2\). Inspection of its proof shows that \(C\)
depends only on \(R=4\), \(k=\frac12\) and \(M_v\), and the same proof bounds \(\sup|\phi|\). The
function \(E_t-\mathrm{id}-t\phi\) is holomorphic on \(B(0,\frac18)\subset\mathbb C\setminus\mathcal K\). By Cauchy's
estimates, \(|E_t(0)-t\phi(0)|\le Ct^2\) and \(|E_t'(0)-1-t\phi'(0)|\le8Ct^2\). Inserting
these into Step~2 gives, uniformly on \(\{|z|\le\frac1{20}\}\) and in \(h\),
\[
 f_{h\circ\psi^v_t}=f+t\dot f+O(t^2),\qquad \dot f=\phi\circ f-\phi(0)-\phi'(0)f
\]
for \(|t|\le t_0\), and trivially for \(t_0\le|t|\le1\), by Koebe's bound.

\emph{Step 4: the contour formula.} Let \(\xi=f(z)\) with \(|z|<\frac14\). Change variables
\(\zeta=f(w)\), \(dA(\zeta)=|f'(w)|^2dA(w)\). Then
\(\dot\sigma(f(w))|f'(w)|^2=-\partial_{\bar w}\widetilde V(w)f'(w)^2\), so
\[
 \phi(f(z))=-\frac1\pi\int\frac{\dot\sigma(\zeta)\,dA(\zeta)}{\zeta-f(z)}
 =\frac1\pi\int_\D\partial_{\bar w}\widetilde V(w)\,\frac{f'(w)^2}{f(w)-f(z)}\,dA(w).
\]
The function \(k_z(w)=f'(w)^2/(f(w)-f(z))\) is holomorphic on \(\D\setminus\{z\}\), and
\(\widetilde V=0\) on \(\{|w|<\frac14\}\ni z\). So \(\widetilde Vk_z\) is smooth on \(\{|w|\le r\}\), \(\frac35<r<1\), with
\(\partial_{\bar w}(\widetilde Vk_z)=(\partial_{\bar w}\widetilde V)k_z\). Moreover \(\widetilde V=\mathcal L\) on \(|w|=r\). Stokes' formula
\(\int_{|w|<r}\partial_{\bar w}F\,dA=\frac1{2i}\oint_{|w|=r}F\,dw\) gives
\[
 \phi(\xi)=\frac1{2\pi i}\oint_{|w|=r}\frac{\mathcal L(w)f'(w)^2}{f(w)-\xi}\,dw
 \qquad(\xi\in f(\{|z|<\tfrac14\})).
\]
Differentiating at \(\xi=0\) and using
\(\frac1{\zeta-\xi}-\frac1\zeta-\frac\xi{\zeta^2}=\frac{\xi^2}{\zeta^2(\zeta-\xi)}\),
\[
 \dot f(z)=\frac{f(z)^2}{2\pi i}\oint_{|w|=r}\frac{\mathcal L(w)f'(w)^2}{f(w)^2(f(w)-f(z))}\,dw
 =i\sum_{|j|\le N}\hat v_j\,\mathbf L_jf(z),
\]
since \(\mathcal L(w)=i\sum_j\hat v_jw^{1+j}\). This proves \eqref{eq:coeffResp-f}.

\emph{Step 5: coefficients.} By Cauchy's formula on \(|z|=\frac1{20}\),
\(|a_m(h\circ\psi^v_t)-a_m(h)-t[z^{m+1}]\dot f|\le20^{m+1}C_vt^2\). This proves
\eqref{eq:coeffResp-a}. The map \(v\mapsto(\hat v_j)_j\) is real-linear, so \(D_va_m\) is
real-linear in \(v\). Its complex-linear extension is \(\mathfrak d_ja_m=iL_ja_m\), and the real
forms follow from \(\cos n\theta=\frac12(e^{in\theta}+e^{-in\theta})\) and
\(\sin n\theta=\frac1{2i}(e^{in\theta}-e^{-in\theta})\). By (ii), \(L_ja_m\in\mathbb C[a_1,\ldots,a_{m+N}]\) for
\(|j|\le N\).
\end{proof}

\begin{corollary}[Iterated responses]\label{cor:coeffBounded}
Extend \(D_v\) to the derivation of \(\mathbb C[a_m,\bar a_m:m\ge1]\) with
\(D_v\bar a_m=\overline{D_va_m}\) and \(D_va_m\) as in \eqref{eq:coeffResp-a}. Let \(\Phi\) be a coefficient
polynomial, \(v\) a real trigonometric field, and \(h\in\mathcal W\). Then
\(t\mapsto\Phi(h\circ\psi^v_t)\) is \(C^\infty\) on \(\mathbb R\), with
\[
 \frac{d^j}{dt^j}\Phi(h\circ\psi^v_t)=(D_v^j\Phi)(h\circ\psi^v_t)\qquad(j\ge0).
\]
Each \(D_{v_1}\cdots D_{v_j}\Phi\) is a coefficient polynomial, real-multilinear in
\((v_1,\ldots,v_j)\). For each \(N\) there is \(C_{\Phi,j,N}\) with
\(\sup_{\mathcal W}|D_{v_1}\cdots D_{v_j}\Phi|\le C_{\Phi,j,N}\|v_1\|\cdots\|v_j\|\) for fields of degree
\(\le N\), for any fixed norm \(\|\cdot\|\). In particular, for the stress scores
\(\rho_w=\sum_n\{x_n\Ima\mathfrak t_n+y_n\Rea\mathfrak t_n\}\) of Section~\ref{sec:finiteinput},
\(D_v\rho_w\) is bilinear in \((v,w)\) and bounded on \(\mathcal W\) uniformly for \(v,w\) on the unit
sphere of fields of degree \(\le N\). Also
\(|\rho_w(h\circ\psi^v_t)-\rho_w(h)-tD_v\rho_w(h)|\le\frac{t^2}2\sup_{\mathcal W}|D_v^2\rho_w|\).
\end{corollary}

\begin{proof}
\(\mathcal W\) is invariant under \(\psi^v_t\) by Lemma~\ref{lem:entQ}, and
\(h\circ\psi^v_{t+s}=(h\circ\psi^v_t)\circ\psi^v_s\). Apply Lemma~\ref{lem:coeffResp}(iii) at
\(h\circ\psi^v_t\) to get differentiability of each \(a_m\), hence of \(\bar a_m\) (\(v\) is real), at
every \(t\). The chain rule for polynomials then gives the derivative formula for
\(j=1\), and induction gives it for all \(j\). Multilinearity follows from
\eqref{eq:coeffResp-a}. By de Branges' theorem, \(|a_m(h)|\le m+1\) on \(\mathcal W\), so a coefficient
polynomial is bounded on \(\mathcal W\) by its maximum on a compact polydisc. The
coefficients of \(D_{v_1}\cdots D_{v_j}\Phi\) are polynomials in the
\(\hat v_{i,k}\), \(|k|\le N\), so they are bounded on the unit sphere. The \(\rho_w\) are real
coefficient polynomials, linear in \(w\), by Corollary~\ref{cor:coeffTm}(i). The last bound is
Taylor's formula.
\end{proof}

\begin{corollary}[Stress coefficients]\label{cor:coeffTm}
Let \(q_f=\frac{\cL}{12}\mathcal Sf+(f'/f)^2-z^{-2}=\sum_{m\ge1}\mathfrak t_mz^{m-2}\).
\begin{enumerate}
\item \(\mathfrak t_m\in\mathbb C[a_1,\ldots,a_m]\) is weighted homogeneous of weight \(m\), and
\(\mathfrak t_m-C_ma_m\in\mathbb C[a_1,\ldots,a_{m-1}]\), where \(C_m=2m+\frac{\cL}{12}(m^3-m)\).
\item For \(m,n\ge1\),
\(L_n\mathfrak t_m=\delta_{mn}C_n+(m+n)\,\mathfrak t_{m-n}\mathbf 1_{m>n}\), and
\(L_0\mathfrak t_m=m\mathfrak t_m\).
\item Hence, at every \(h\in\mathcal W\),
\(\mathfrak d_n\mathfrak t_m=i\{\delta_{mn}C_n+(m+n)\mathfrak t_{m-n}\mathbf 1_{m>n}\}\) for \(n\ge1\).
\end{enumerate}
\end{corollary}

\begin{proof}
\emph{(i)} We have \(f'/f=z^{-1}+u'/u\), and \(1/f'\), \(1/u\) have coefficients in
\(\mathbb C[a_1,a_2,\ldots]\). So \(q_f\) is a Laurent series with coefficients in this ring, and
it has at most a simple pole at \(0\). The coefficient \([z^{m-2}]\mathcal Sf\) involves only
\([z^k]f'\) for \(k\le m\), and \([z^{m-2}](f'/f)^2\) only \([z^k](u'/u)\) for \(k\le m-1\); both lie in
\(\mathbb C[a_1,\ldots,a_m]\). For \(f_r(z)=f(rz)/r\) we have \(a_k(f_r)=r^ka_k\) and
\(q_{f_r}(z)=r^2q_f(rz)\), hence \(\mathfrak t_m(f_r)=r^m\mathfrak t_m(f)\), which is the weight statement. The
only monomial of weight \(m\) containing \(a_m\) is \(a_m\) itself. Its coefficient is
\(\frac d{d\varepsilon}\mathfrak t_m(z+\varepsilon z^{m+1})|_0=\frac{\cL}{12}(m^3-m)+2m=C_m\).

\emph{(ii)} Both sides are polynomials in \(a_1,\ldots,a_m\), since \(L_n\), \(n\ge1\), lowers
indices. It therefore suffices to prove the identity for every polynomial
\(f=z(1+\sum_{k\le m}a_kz^k)\). All functions below are then holomorphic near \(0\), where
\(f'\ne0\) and \(f/z\ne0\). The coefficient \([z^{m-2}]q_F\) equals \(\mathfrak t_m\) evaluated at the
coefficients of any \(F=z+O(z^2)\) holomorphic near \(0\). For \(n\ge1\), \(L_n\) is the derivative along \(\delta f=z^{n+1}f'\) (Lemma~\ref{lem:coeffResp}(i)):
\(L_n\Phi=\frac d{d\varepsilon}|_0\Phi(f+\varepsilon z^{n+1}f')\) for \(\Phi\in\mathbb C[a_1,a_2,\ldots]\). With
\(\phi_\varepsilon(z)=z+\varepsilon z^{n+1}\), we have \(f\circ\phi_\varepsilon=f+\varepsilon z^{n+1}f'+O(\varepsilon^2)\) coefficientwise,
and \(f\circ\phi_\varepsilon\) is normalized. The Schwarzian chain rule and
\(\bigl((f\circ\phi)'/(f\circ\phi)\bigr)^2=\bigl((f'/f)^2\circ\phi\bigr)\phi'^2\) give the cocycle
\[
 q_{f\circ\phi}=(q_f\circ\phi)\,\phi'^2+q_\phi .
\]
Here \(q_{\phi_\varepsilon}=\varepsilon\{\frac{\cL}{12}(n^3-n)+2n\}z^{n-2}+O(\varepsilon^2)=\varepsilon C_nz^{n-2}+O(\varepsilon^2)\). Moreover
\((q_f\circ\phi_\varepsilon)\phi_\varepsilon'^2=q_f+\varepsilon\bigl(z^{n+1}q_f'+2(n+1)z^nq_f\bigr)+O(\varepsilon^2)\), and
\(z^{n+1}q_f'+2(n+1)z^nq_f=\sum_k(k+2n)\mathfrak t_kz^{k+n-2}\). Taking \([z^{m-2}]\) gives the claim,
because \(\mathfrak t_k\) is defined only for \(k\ge1\). The identity \(L_0\mathfrak t_m=m\mathfrak t_m\) is (i).

\emph{(iii)} follows from (ii) and Lemma~\ref{lem:coeffResp}(iii).
\end{proof}

\begin{remark}[Checks]\label{rem:coeffRespChecks}
At \(f=z\), \(D_{\cos n\theta}f=\frac i2z^{n+1}\) (Proposition~\ref{prop:potential}); on the M\"obius orbit
of Remark~\ref{rem:entVchecks}(a) (\(v=-2\sin\theta\)), \eqref{eq:coeffResp-f} is the tangent of \(z/(1+\tanh(s)z)\). By computer algebra,
the Witt relations hold on \(a_1,\dots,a_4\) (\(|m|,|n|\le4\); not used),
and Corollary~\ref{cor:coeffTm}(ii) holds for \(m,n\le7\).
\end{remark}

\section*{Statements and declarations}
\noindent\textbf{Funding.} During the course of this project, V.O.H.~was partially supported by an AMS Simons Research Enhancement Grant for PUI
Faculty. 
The work of J. Yum was supported by the National Research Foundation of Korea(NRF) grant funded by the Korea government(MSIT) (RS-2026-25470773).\\
\textbf{Competing interests.} The authors have no competing interests to
declare.\\
\textbf{Data availability.} No datasets were generated or analysed. The
computer-algebra and numerical scripts used for the checks reported in
the paper (Lemmas~\ref{lem:Wp-basis}--\ref{lem:Wp-ghost},
Corollary~\ref{cor:coeffTm}, Lemma~\ref{lem:GQWanalytic}) are available from the authors upon request.\\
\textbf{Use of AI tools.} This paper is the result of work carried out by
the authors over the past two years (since 2024). More recently, the
authors used AI tools to help explore the literature; this led us to the
preprints of Cai and to the discovery of a significant overlap with that
work. The authors also used AI tools to help edit the text and to check
some computations, as noted above. The authors are solely responsible for
verifying all arguments and calculations and for the final content.
\end{document}